\documentclass[preprint,10pt]{article}
\usepackage[usenames]{color} 
\usepackage{amssymb} 
\usepackage{amsmath,graphicx} 
\usepackage[utf8]{inputenc} 
\usepackage{hyperref}
\usepackage{adjustbox}
\usepackage{wrapfig}
\usepackage{tikz}
\usetikzlibrary{decorations.pathreplacing}
\usetikzlibrary{fadings}
\usepackage{bbm}
\usepackage{bm}
\usepackage{algorithm,algorithmic}
\usepackage{subfigure}
\usepackage[nocompress]{cite}

\usepackage{footnote}
\makesavenoteenv{tabular}
\makesavenoteenv{table}

\usepackage[capitalise,noabbrev]{cleveref}
\def \bmu {\bm{\mu}}

\def \EFE {\mathcal{E}}
\def \ERB {\mathrm{E}}

\def \R {\mathbb{R}}
\def \V {\mathbb{V}}

\def \phiRB {\psi^i_{RB}}

\def \calN {\mathcal{N}}
\def \NFE {\calN}
\def \NRB {N}
\def \VFE {\V_\NFE}
\def \VRB {\V_{\NRB}}
\def \VRBp {\V_{\NRB+1}}
\def \uFE {u_{\NFE}}
\def \uRB {u_{\NRB}}
\def \uRBc {\mathrm{u}}
\def \vFE {v_{\NFE}}
\def \vRB {v_{\NRB}}
\def \vRBc {\mathrm{v}}
\def \phiRB {\psi}
\def \dimtrain {N_{train}}
\def \NEIM {N_{EIM}}
\def\dimEIM{\NEIM}
\def\dimRB{\NRB}

\def \RefDom {\mathcal{R}}
\def \EIMdof {\boldsymbol{\tau}^{EIM}}
\def \EIMprdof {\boldsymbol{\tau}^{EIM'}}
\def \EIMDof {\boldsymbol{\Sigma}^{EIM}}
\def \EIMfun {\boldsymbol{\rho}^{EIM}}
\def \EIMprfun {\boldsymbol{\rho}^{EIM'}}
\def \EIMFun {\boldsymbol{\mathcal{Q}}^{EIM}}
\def \EIMint {\mathcal{I}_{EIM}}
\def \proj {\Pi}
\def\EIM{\textrm{EIM}}
\def\POD{\textrm{POD}}
\def\RB{\textrm{RB}}

\def \II {\mathcal{I}}

\def \P {\mathcal{P}}
\DeclareMathOperator*{\argmax}{arg\,max}

\def \R {\mathbb{R}}

\def \P {\mathbb{P}}

\newcommand{\DD}{{\mathcal{ D}}}

\definecolor{lightgreen}{rgb}{0,1,0}

\newtheorem{oss}{Remark}
\newtheorem{ex}{Example}

\def\rbold{\mathbf{r}}

\def\vbold{\mathbf{v}}

\def\mubold{\boldsymbol{\mu}}

\def\0bold{\boldsymbol{0}}
\def\0{\mathbf{\0}}

\title{Model Reduction for Advection Dominated Hyperbolic Problems\\ in an ALE Framework: Offline and Online Phases} 
\author{Davide Torlo\footnote{I\lowercase{nstitut f\"ur} M\lowercase{athematik,} W\lowercase{interthurstrasse} 190, CH 8057  Z\lowercase{\"urich}, S\lowercase{witzerland}.} \footnote{ C\lowercase{orresponding author,  (\href{mailto:davide.torlo@math.uzh.ch}{davide.torlo@math.uzh.ch}).}}  }
\date{\today}

\begin{document}
\maketitle
\begin{abstract}
	Model order reduction (MOR) techniques have always struggled in compressing information for advection dominated problems. 
	Their linear nature does not allow to accelerate the slow decay of the Kolmogorov $N$--width of these problems.
	In the last years, new nonlinear algorithms obtained smaller reduced spaces. In these works only the \textit{offline} phase of these algorithms was shown.
	In this work, we study MOR algorithms for unsteady parametric advection dominated hyperbolic problems, giving a complete \textit{offline} and \textit{online} description and showing the time saving in the \textit{online }phase. 
	We propose an arbitrary Lagrangian--Eulerian approach that modifies both the \textit{offline} and \textit{online} phases of the MOR process.
	This allows to \textit{calibrate} the advected features on the same position and to strongly compress the reduced spaces.
	The basic MOR algorithms used are the classical Greedy, EIM and POD, while the \textit{calibration} map is learned through polynomial regression and artificial neural networks. 
	In the performed simulations we show how the new algorithm defeats the classical method on many equations with nonlinear fluxes and with different boundary conditions. Finally, we compare the results obtained with different \textit{calibration} maps.
\end{abstract}
\section{Introduction}
Advection dominated problems are very common in many engineering applications where particles or waves travel in space at fast speeds, for example in shallow water equations, gas dynamics simulations or explosions. They arise frequently in hyperbolic conservation laws, since their weak solutions develop shocks from smooth initial data and typical simulations consist of shock propagations. They can be found also in advection--diffusion--reaction problems when the advection parameter dominates the equation.

In this work we contribute to the development of model order reduction (MOR) algorithms for advection dominated problems. 
MOR techniques are extremely useful and sometimes of vital importance to reduce the computational time of parametric problems where many simulations or fast simulations are necessary. 

The extension of MOR algorithms to advection dominated problems is anyway not straightforward. 
It is well--known that these problems suffer of a slow decay of the Kolmogorov $N$--width.
Indeed, if one wants to compress the information of, for example, a traveling shock in a linear subspace, the number of basis functions needed is proportional to the number of degrees of freedom of the discretized domain. 
Hence, all the classical MOR techniques are not suited for these types of problem.

Nevertheless, much effort has been devoted in the last years in adapting and creating new MOR techniques for advection dominated problems. 
In particular, nonlinear approaches have been used and they can be distinguished in two big classes: the Eulerian methods, where the reduced structures are adapting to the moving frame, \textit{inter alia} \cite{sPOD,peherstorfer18adaptiveBases, amsallem08interpolation, rim18displacementInterpolation, zimmermann18onlineAdaptive}, and the Lagrangian methods, where the solutions are transformed before applying the MOR techniques \cite{RB_freezing,trasport_greedy,iollo2014advection, mojgani17aleRB, Cagniart2019,  taddei19registration, nair19transportedSnapshots,nonino2019overcoming}. 
Even being more complicated to construct, the Lagrangian methods allow to fully apply the classical MOR techniques in the new framework. 
All the previous works, even tackling the problematics of the topic, do not present a general approach for a parametric time--dependent problem. Namely, none of them proposed a complete MOR algorithm with an \textit{offline} and an \textit{online} phase, where the reduction in computational time is observable.
A different remark must be done for the work in \cite{Lye19DeepLearning}, where deep learning are used to predict the parameter--to--observable map. If considered as a MOR algorithm, it comprehend an \textit{online} and \textit{offline} phase, but it does not predict the whole solution of a parametric problem, but only few features, like physical relevant quantities, as drag and lift, or statistical averages.

The work that we present wants to solve a wide class of advection--dominated problems and shows the direction to develop more MOR techniques for many other problems.
It is based on a Lagrangian approach. In particular, we want to align all the solutions in order to minimize the effort of the MOR algorithms.
This leads to an arbitrary Lagrangian--Eulerian (ALE) formulation of the equations, that must be solved both in the \textit{offline} and \textit{online} phase of the algorithms.
Now, the classical \textit{speed} of the mesh of the ALE framework must be found with new techniques. 
We propose different regression maps to learn this speed, in order to use it cheaply also in the \textit{online} phase.
Here, regression algorithms, as polynomial least square regression or more recent neural networks, are compared.

The overall MOR algorithm makes use of this new ALE framework.
It also applies the classical Greedy algorithm in the parameter domain, the proper orthogonal decomposition on the time evolution and the empirical interpolation method to compress the nonlinearities of the fluxes.
The ALE framework exports these methods to a reference domain, allowing to compute a complete \textit{online} phase, where practical advantages in computational times are shown for different types of equations.

The work is organized as follows. In \cref{sec:preli} we introduce the classical MOR techniques that we use in the new algorithm and we highlight the difficulties in reducing advection--dominated problems with these algorithms.
In \cref{sec:ALE} we present the ALE framework where we can obtain aligned solutions that can be easily reduced.
In \cref{sec:transport_map} we present the regression maps that we use in the simulations and the learning algorithms that we need to optimize their parameters.
In \cref{sec:results} we demonstrate the quality of the algorithm on typical advection--dominated problems on which the classical MOR algorithms fail or struggle to compress information.
Finally, in \cref{sec:conclusion} we remark the new features of the proposed algorithm and we suggest new possible directions of investigation.

\section{Preliminaries}\label{sec:preli}
Before introducing the new methodology of this work, we present the standard techniques used to perform MOR on hyperbolic problems and some of the main troubles that we face when dealing with advection dominated problems. 
In particular, MOR techniques were first introduced to gain computational costs in diffusion dominated problems for elliptic and parabolic equations. The extension to hyperbolic problems lacks of tools, theorems and results that guarantee the quality of the obtained results. 
Nevertheless, many works tried to advance the study of MOR in the hyperbolic context, leading to interesting results, but still leaving many open problems.

\subsection{Model Order Reduction for Hyperbolic Problems}
As a benchmark MOR algorithm, we refer to Haasdonk et al.\cite{haasdonk_pod_greedy2008, Haasdonk2009, Drohmann2012} algorithms, which are suited for hyperbolic problems with different solvers. The most general one that we will consider consists of a PODEI--Greedy, which combines a proper orthogonal decomposition (POD) compression of the solutions in time, a Greedy in parameter space and an empirical interpolation method (EIM) to approximate the nonlinearities of fluxes and scheme operators.

Let us consider an interval $\Omega \subset \R$, the final time $t_f\in\R^+$, a parameter set $\P\subset \R^P$, a scalar unknown $u: \Omega \times [0,t_f] \times \P \to \R$ and the conservation law 
\begin{equation}\label{eq:cons_law}
\begin{cases}
\partial_t u (x,t,\bmu) + \frac{d}{dx} F(u, \bmu) =0, \quad &x\in \Omega,\, t\in [0,t_f], \, \bmu \in \P\subset \R^P, \\
\mathbf{B}(u, \bmu)  =0, &x \in \partial \Omega \\
u(x,t=0, \bmu) = u_0 (x,\bmu), & x\in \Omega  ,
\end{cases}
\end{equation}
where $\mathbf{B}$ is a boundary operator and $u_0$ some initial conditions. $F$ is a flux function $F:\R \times \P \to \R$, possibly nonlinear in all the inputs. More precisely, we aim to get the weak solution corresponding to \eqref{eq:cons_law}, but, in order to keep a simple notation, we will always write the equation for the strong solution.

Suppose that we have a discrete evolution scheme for an approximation $\uFE (t,\bmu)$ of the solution $u (\bmu,t)$ of \eqref{eq:cons_law}. Here, $\calN$ denotes the dimension of $\VFE \subset L^2(\Omega)$, a finite dimensional approximation space, where we search the full order model (FOM) solution $\uFE (t,\bmu)$. The space can be represented as a span of linearly independent basis functions $\lbrace\varphi_\sigma \rbrace_{\sigma\in\Sigma}$, where $\Sigma$ is the set of the degrees of freedom, and each solution in the space will be characterized by the coefficients $u_\sigma$, such that $\uFE(t,\bmu)=\sum_{\sigma \in \Sigma} u_\sigma(t,\bmu) \varphi_\sigma$.  After the time discretization $0=t^0, \dots , t^k, \dots , t^K=T$, the discrete (explicit) FOM, can be written in the following way:
\begin{equation}
\uFE^k(\bmu) = \uFE^{k-1}(\bmu) - \EFE(F(\uFE^{k-1}(\bmu),\bmu)),\quad k=1, \dots , K, \, \forall \bmu \in \P.
\end{equation}
Here, we denote with $\uFE^k$ the discretized variable at the time step $t^k$ and $\EFE: \VFE \to \VFE$ is an evolution operator given by a compact stencil scheme, e.\,g. finite volume, finite difference, finite element or residual distribution schemes \cite{leveque92FV,leveque07FD,quarteSaccoSaleri,abgrall2006residual}.  

Given a reduced space $\VRB= \left\lbrace \sum_{i=1}^\NRB \uRBc_i \phiRB^i : \uRBc_i \in \R \right\rbrace$ generated by selected basis functions $\lbrace \phiRB^i \rbrace_{i=1}^\NRB \subset \VFE$ provided by the PODEI-Greedy, we want to find a reduced solution $\uRB\in \VRB$ for the following reduced problem. The reduced basis space $\VRB$ is also coupled with a projection operator $\proj: \VFE \to \VRB$ given by the Galerkin projection. Given a parameter $\bmu\in \P$, solve
\begin{equation}\label{eq:RB_problem}
\uRB^{k+1}(\bmu) :=\sum_{i=1}^\NRB \uRBc^{k+1}_i(\bmu) \phiRB^i =\sum_{i=1}^\NRB \uRBc^{k}_i(\bmu) \phiRB^i - \sum_{i=1}^\NRB\ERB_i(F(\uRB^k(\bmu),\bmu))\phiRB^i.
\end{equation} 
Here, $\ERB: \VFE \to \R^\NRB$ is the reduced evolution operator. This will be provided by the EIM algorithm, through evaluation in few \textit{magic} points of the flux of the variable at the previous time step. To define the reduced basis space $\VRB$ and the operator $\ERB$, we will perform the PODEI--Greedy algorithm.
\subsubsection{Greedy Algorithm}
The key algorithm to generate the reduced space is the Greedy algorithm \cite{prudhomme2001,prud2002}. It selects iteratively the worst approximated snapshot and updates the reduced space with the information of the chosen snapshot. It requires the knowledge of an error measure or an estimation of it. The algorithm stops when the worst error is smaller than a prescribed tolerance. It can be summarized as in Algorithm \ref{algo:greedy}. In this framework, it requires the following procedures:
\begin{itemize}
	\item \textsc{InitBasis} which initializes the reduced basis (RB) space $\DD_N$;
	\item \textsc{ErrorEstimate} which estimates the error between the high--fidelity function and its projection on the RB space $\DD_N$;
	\item \textsc{UpdateBasis} which updates the RB space $\DD_N$, given a certain selected parameter.
\end{itemize}
\begin{algorithm}
	\fontsize{10pt}{10pt}\selectfont
	\caption{Greedy Algorithm} 
	\begin{algorithmic}[1]
		{\REQUIRE Training set $\mathcal{M}_{train} = \lbrace \bmu_i \rbrace_{i=1}^{\dimtrain}$, tolerance $\varepsilon^{tol}$ and $N_{max}$.
			\ENSURE RB space $\DD_N$
			\STATE Initialize an RB space of dimension $N$:\\
			$\DD_{N}$= \textsc{InitBasis}
			\WHILE{$\max_{\bmu \in \mathcal{M}_{train} }  $\textsc{ErrorEstimate}$(\DD_N,\uFE(\bmu)) \geq \varepsilon^{tol}$  AND $ N \leq  N_{max}$ } 
			\STATE Find the parameter of worst approximated snapshot:\\
			$\bmu_{max}=\argmax_{\bmu \in \mathcal{M}_{train}} \textsc{ErrorEstimate} (\uFE(\bmu), \DD_N)$
			\STATE Extend reduced basis $\DD_N$ with the found snapshot (adding the new snapshot to dictionary):\\
			$\DD_N, N=$\textsc{UpdateBasis}($\DD_N,\uFE (\bmu_{max})$)
			\ENDWHILE
		}
	\end{algorithmic}\label{algo:greedy}
\end{algorithm} 

\subsubsection{EIM Algorithm}
The empirical interpolation method (EIM or EI) \cite{barrault04} is a technique that allows to interpolate nonlinear functions into a set of \textit{magic} degrees of freedom (DoF) $\EIMDof :=\lbrace \EIMdof_m \rbrace_{m=1}^{\NEIM}\subset \VFE^*$ and basis functions 
$\EIMFun :=\lbrace \EIMfun_m \rbrace_{m=1}^{\NEIM}\subset \VFE$, 
in order to approximate the nonlinear functions with few evaluations. In particular, we need to use this procedure to interpolate the evolution operators for different times and parameters in the reduced setting, writing 
\begin{equation}
\EIMint [\EFE (F(\uFE,\bmu))] = \sum_{m=1}^{\NEIM} \EIMdof_m \left( \EFE (F(\uFE,\bmu)) \right) \EIMfun_m \approx  \EFE (F(\uFE,\bmu)) .
\end{equation}
\begin{algorithm}
	\fontsize{10pt}{10pt}\selectfont
	\caption{Empirical Interpolation Method($\EFE$)} 
	\textsc{EIM--InitBasis()}
	\begin{algorithmic}[1]
		{	\RETURN empty initial basis $\DD_{0} = \emptyset $
		}
	\end{algorithmic}
	\hrule
	\vspace{1mm}
	\textsc{EIM--ErrorEstimate}($(\EIMFun_{M},\EIMDof_{M}),\bmu, t^k, \EFE$ )
	\begin{algorithmic}[1]
		{	\STATE Compute the exact flux $\vbold_h=\EFE (F(\uRB^k(\bmu),\bmu))$\\
			\STATE Compute the interpolation coefficients $\sigma(\vbold_h):= (\sigma_j)_{j=1}^M\in \R^M$ by solving a linear system (upper triangular)
			$
			\sum\limits_{j=1}^M \sigma_j(\vbold_h)\EIMdof_i[\EIMfun_j]=\EIMdof_i[\vbold_h],$ for  $i=1,\dots,M
			$
			\RETURN approximation error $|| \vbold_h- \sum_{j=1}^{M}\sigma_j^M(\vbold_h)\EIMfun_j||_{\VFE}$
		}
	\end{algorithmic}
	\hrule
	\vspace{1mm}
	\textsc{EIM--UpdateBasis} ($(\EIMFun_{M},\EIMDof_{M}), \bmu_{max}, t^{k_{max}}, \EFE$)
	\begin{algorithmic}[1]
		{	\STATE Compute the exact flux $\vbold_h=\EFE (F(\uRB^k(\bmu_{max}),\bmu_{max}))$\\
			\STATE Compute the interpolation coefficients 	$\sigma(\vbold_h):= (\sigma_j)_{j=1}^M\in \R^M$ \\
			\STATE Compute the residual between the truth flux and its interpolant 
			$\rbold_M=  \vbold_h- \sum_{j=1}^{M}\sigma_j^M(\vbold_h)\EIMfun_j$\\
			\STATE Find the DoF that maximize the residual 
			$\EIMdof_{M+1} := \argmax_{\EIMdof \in \EIMDof_h} |\EIMdof (\rbold_M)|$\\
			\STATE Normalize the correspondent basis function
			$\EIMfun_{M+1}:=\EIMdof_{M+1}(\rbold_M)^{-1} \cdot \rbold_M$ \\
			\RETURN updated basis $\DD_{M+1}:=((\EIMfun_m)_{m=1}^{M+1}, (\EIMdof_{m})_{m=1}^{M+1}) $
		}
	\end{algorithmic}\label{algo:EIM}
\end{algorithm}

In this way, one can compute just the interpolation DoFs of the flux (which depend on the parameter $\bmu$) in an online phase and precompute the projection onto the reduced basis space of the interpolation functions $\EIMfun$ in an offline phase, i.\,e., $\proj(\EIMfun_m)$. Once stored these values, we can operate in a cheap way the computation of the reduced evolution operator, i.\,e.,
\begin{equation}\label{eq:EIM_interpolation}
\ERB_i(F(\uRB^k(\bmu),\bmu)) := \sum_{m=1}^{\NEIM} \EIMdof_m \left( \EFE (F(\uRB^k(\bmu),\bmu)) \right) \proj_i(\EIMfun_m).
\end{equation}
The error estimator in the EIM algorithm is an actual error computation between the interpolated functions and the original evolution operators.
The algorithm can be recast into a Greedy Algorithm \ref{algo:greedy}, through the specifications of Algorithm \ref{algo:EIM}.

Here, we can observe the importance of the compact stencil assumption on the evolution operator $\EFE$. Indeed, if, for example, we consider as DoFs $\EIMdof$ the point evaluations, the compactness of the stencil plays an important role to minimize the flux $F$ evaluations to compute the evolution operator.

More details on this algorithm can be found in \cite{barrault04,haasdonk_pod_greedy2008, Haasdonk2009, Drohmann2012, crisovan19uq}.
\subsubsection{POD--Greedy Algorithm}
The POD--Greedy algorithm is used to perform a reduction on a time-parameter domain. The Greedy algorithm is used to select the worst approximated time evolution snapshot $\lbrace \uFE^k(\bmu)\rbrace_{k=0}^K$ in the parameter domain $\P$.
The POD is used to extract information from a single time evolution snapshot, provided by the Greedy algorithm, to obtain few basis functions that summarize the time evolution. 
Then the new information is added to the previous RB and another POD is computed over it. 
In this algorithm it is crucial to have an inexpensive error indicator to minimize the computational cost during both the offline and the online phase. We will define it in \cref{sec:error}. 

The specification of the algorithm can be found in Algorithm \ref{algo:POD_greedy}.
\begin{algorithm}[h]
	\fontsize{10pt}{10pt}\selectfont
	\caption{POD--Greedy($\EFE$)} 
	\textsc{POD--Greedy--InitBasis()}
	\begin{algorithmic}[1]
		{	\STATE Pick a parameter $\bmu$ and compute the solution through all the time steps $t^k$: $ \lbrace \uFE^{k}(\bmu) \rbrace_{k=1}^K$\\
			\RETURN initial basis $\DD_{0} = \POD (\lbrace \uFE^{k} (\bmu) \rbrace_{k=1}^K )$
		}
	\end{algorithmic}
	\hrule
	\vspace{1mm}
	\textsc{POD--Greedy--ErrorEstimate}($\VRB, \mubold, t^k$ )
	\begin{algorithmic}[1]
		{	\RETURN error indicator $\eta^k_{\NRB ,\NEIM} (\bmu) \geq ||\uFE^{k}(\bmu)-\uRB^{k} (\bmu) ||_{\VFE}$
		}
	\end{algorithmic}
	\hrule
	\vspace{1mm}
	\textsc{POD--Greedy--UpdateBasis} ($\VRB, \bmu_{max} , \EFE$)
	\begin{algorithmic}[1]
		{	\STATE Compute the exact solution for all time steps $\lbrace \uFE^{k}(\bmu_{max}) \rbrace_{k=1}^K  $ with high--fidelity evolution operator $\EFE$  \\
			\STATE Compute the Galerkin projection of the solution onto the $\RB$ space $\proj[\uFE^{k}(\bmu_{max})]$
			\STATE Compute the POD over time steps of the orthogonal projection of the high--fidelity solution \\ 
			$\VRB^{add}=\POD\left( \lbrace \proj[\uFE^{k}(\bmu_{max})]-\uFE^{k}(\bmu_{max})\rbrace_{k=1}^K\right)$\\
			\STATE Compute a second POD to get rid of extra information \\ 
			$\VRB=\POD(\VRB^{add} \cup \VRB)$\\
			\RETURN updated basis $\VRB $
		}
	\end{algorithmic}\label{algo:POD_greedy}
\end{algorithm}
Hidden in this procedure is the EIM method, which must be used to compute the projection onto the reduced basis of the solutions as described in \eqref{eq:RB_problem} and \eqref{eq:EIM_interpolation}. 
More details on the algorithm can be found in \cite{haasdonk_pod_greedy2008, Haasdonk2009, Drohmann2012, RB_book_rozza, crisovan19uq}.

\subsubsection{PODEI--Greedy Algorithm}
The final algorithm presented in \cite{Drohmann2012} is a combination of the previous ones, where the EIM and the RB enrichment are synchronized to minimize the number of basis functions in both spaces, always reaching the desired precision. This choice is motivated also by other reasons specified in \cite{Drohmann2012, crisovan19uq}. 
The differences with respect to the POD--Greedy algorithm are the fact that also the EIM space is getting enriched every Greedy iteration and that the new RB space can be discarded if the error is increasing. 
This means that if the EIM space is not rich enough to obtain sufficiently reliable RB solutions, then, a refinement only in the EIM space is performed.

The final Algorithm can be rewritten into a Greedy Algorithm \ref{algo:greedy} with specifications of Algorithm \ref{algo:PODEIM_greedy}.

\begin{algorithm}
	\fontsize{10pt}{10pt}\selectfont
	\caption{PODEIM--Greedy($\EFE$)} 
	\textsc{PODEIM--Greedy--InitBasis($\mathcal{M}_{train},\EFE$)}
	\begin{algorithmic}[1]
		{	\STATE $D^{\EIM}_{\NEIM}$=$(\EIMFun_{M_{small}} , \EIMDof_{M_{small}}) = \textsc{EIM-Greedy}(\mathcal{M}_{train}, \varepsilon_{tol,small},\EFE)$
			\STATE Pick a parameter $\bmu$ and compute the solution through all the time steps $t^k$: $ \lbrace \uFE^{k}(\bmu) \rbrace_{k=1}^K$\\
			\STATE $\VRB = \POD (\lbrace \uFE^{k} (\bmu) \rbrace_{k=1}^K )$
			\RETURN  initial bases $\DD_{0}=(\VRB , D^{\EIM}_{\NEIM})$
		}
	\end{algorithmic}
	\hrule	
	\vspace{1mm}
	\textsc{PODEIM--Greedy--ErrorEstimate}($\DD_{S}, \bmu, t^k ,\EFE$ )
	\begin{algorithmic}[1]
		{	\RETURN error indicator $\eta^k_{\NRB ,\NEIM} (\bmu)$
		}
	\end{algorithmic}
	\hrule
	\vspace{1mm}
	\textsc{PODEIM--Greedy--UpdateBasis} ($\DD_{S}, \bmu_{max} ,\EFE$)
	\begin{algorithmic}[1]
		{	\STATE Extend EIM basis $D^{\EIM}_{\NEIM+1} = $ \textsc{EIM--UpdateBasis} $(D^{\EIM}_{\NEIM}, \bmu_{max},\EFE)$
			\STATE Extend RB basis $\VRBp =$ \textsc{POD--Greedy--UpdateBasis} $(\VRB, \bmu_{max},\EFE)$ \\
			\STATE Discard extended $\VRBp$ if error increases:
			\IF{ $\eta^k_{\NRB-1,\NEIM -1}(\bmu_{max})< \max_{\bmu_i \in \mathcal{M}_{train} } \eta^k_{\NRB,\NEIM }$} 
			\RETURN only {\EIM} updated basis: $\DD_{S+1} = (\VRB,D^{\EIM}_{\NEIM+1})$\\
			\ELSE 	\RETURN updated basis $\DD_{S+1} = (\VRBp ,D^{\EIM}_{\NEIM+1})$	\\
			\ENDIF
		}
	\end{algorithmic}\label{algo:PODEIM_greedy}
\end{algorithm}

\subsubsection{Error Indicator}\label{sec:error}
In order to perform a cheap error indicator $\eta_{\dimRB,\dimEIM,\dimEIM'}$ on the reduced solutions, we adopt the technique described in \cite{haasdonk_pod_greedy2008,crisovan19uq}, where $\dimEIM'$ extra basis functions of the EIM space are stored to give a measure of the error that the reduced fluxes have on this extra EIM space.
\begin{equation}\label{eq:error_estimator}
\lVert \uFE^K(\bmu)- \uRB^K(\bmu) \rVert_{\VFE} \leq \eta^K_{\dimRB,\dimEIM,\dimEIM'}(\mubold):= \sum_{k=1}^K C^{K-k} \!\! \left( \sum_{m=1}^{\dimEIM'} \!\! \Delta t \xi^k_m (\mubold)  \left\lVert  \EIMprfun_{m} \right\rVert_{\VFE}\!\!\!\!+ \!  ||\Delta t R^k(\mubold)||_{\VFE}\! \right),
\end{equation}
where $C$ is the Lipschitz continuity of $\EFE$,
\begin{equation}
\Delta t R^k(\mubold) := \uRB^{k}(\mubold) -\uRB^{k-1}(\mubold)+ \Delta t \II_{\dimEIM} [\EFE(F(\uRB^{k-1}(\mubold),\bmu)) ], \end{equation}
$\II_{\dimEIM}$ is the interpolation operator into basis functions and magic points of the EIM space 
and the coefficients $\xi^k_m(\mubold)$ are
\begin{equation}
\xi^k_m(\mubold)= \EIMprdof_m \left( \EFE (F(\uRB^{k-1}(\mubold),\bmu)) \right), \, \forall m= 1, \dots, \dimEIM'.
\end{equation}
The so--defined error estimator is an error bound under strict hypotheses, which are not always met in the simulations. Nevertheless, in practice, it shows a reliable behavior even with small $\NEIM'$ also when the hypotheses are not fulfilled \cite{Drohmann2012}. Hence, we choose the number of extra EIM bases $\NEIM'$ to be equal to 5 in all the simulations.

\subsubsection{Online Phase}
The \textit{online} of the PODEI-Greedy algorithm is given by the simple formula \eqref{eq:RB_problem}, where the reduced evolution operator is computed with the EIM algorithm as in \eqref{eq:EIM_interpolation} and the total number of flux evaluation computed are at most $\NEIM \cdot \NRB$ instead of $\NFE$.
This allow a computational saving when $$\NEIM ,\NRB \ll \NFE.$$

\subsection{Challenges in Model Order Reduction for Advection Dominated Problems}
This algorithm, but also all the other algorithms based on POD, shows difficulties in capturing the behavior of advection dominated problems. This is a common problem in many model order reduction techniques. Indeed, the classical MOR methods are better suited for diffusion dominated problems. Nevertheless, there are several works that tried to apply the classical techniques to advection dominated problems \cite{torlo18weighted, pacciarini14advectionDominated, crisovan19uq} resulting, most of the time, in very expensive online phase and using extra diffusion to smoothen out the advection phenomena.

\begin{figure}[ht]
\begin{center}
\begin{subfigure}[Snapshots at different times\label{fig:advection_wave_no_calib}]{\includegraphics[width=0.45\textwidth, trim={30 30 30 20},clip]{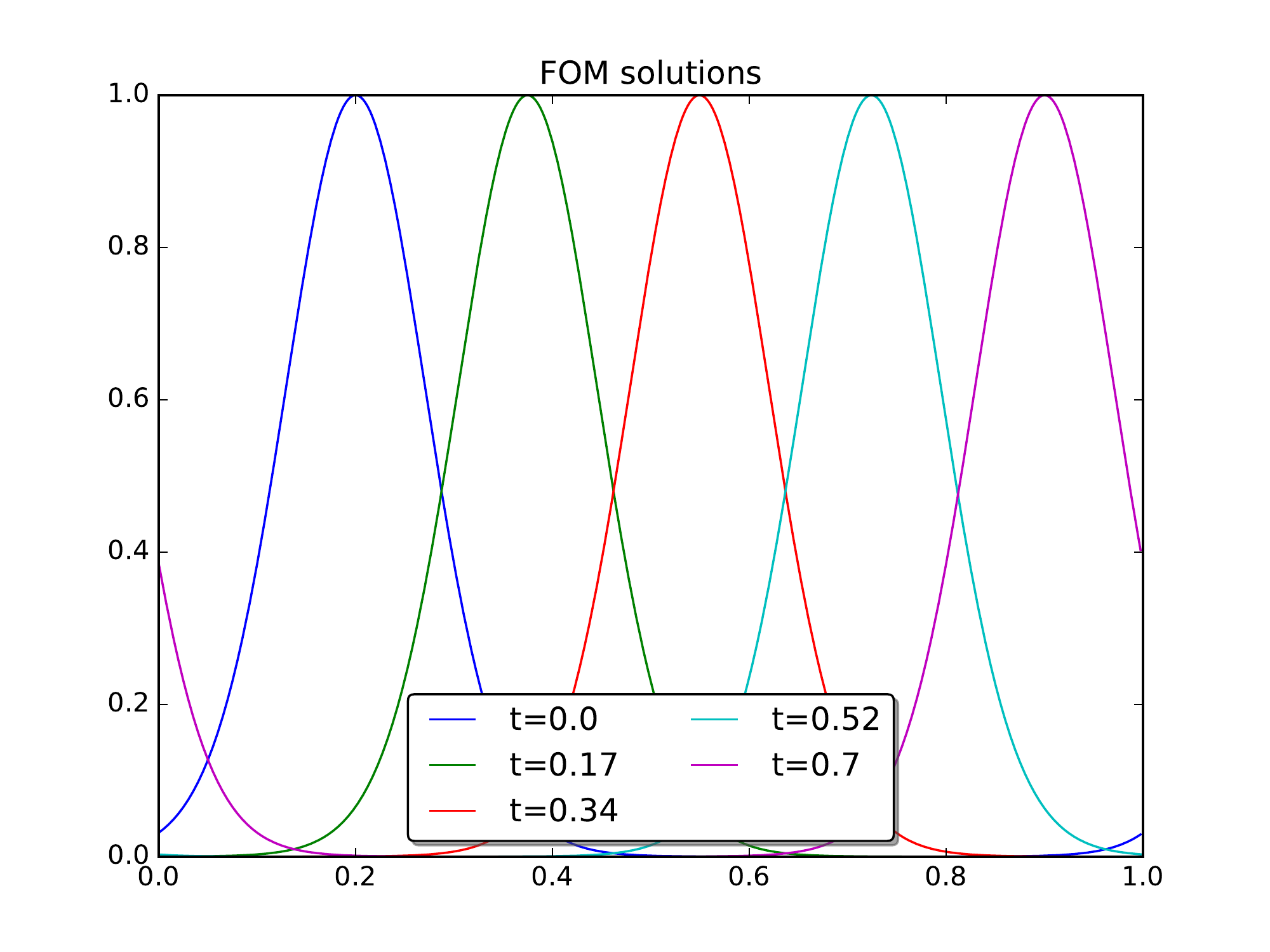}}
\end{subfigure}
\begin{subfigure}[POD modes\label{fig:POD_advection_wave_no_calib}]{\includegraphics[width=0.45\textwidth, trim={30 30 30 20},clip]{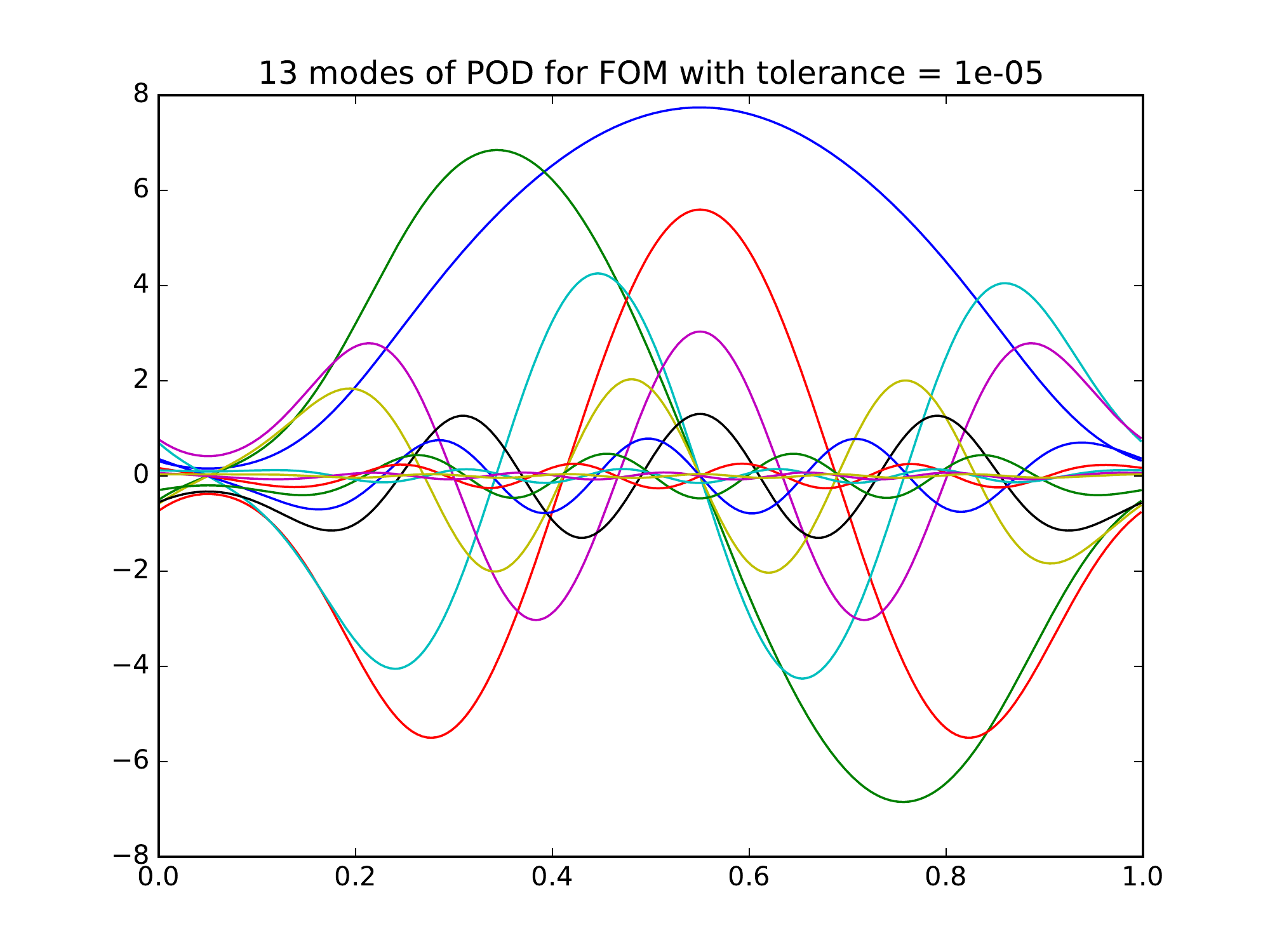}}
\end{subfigure}
\end{center}
\caption{Smooth solution for wave equation}
\end{figure}

We can point out very simple examples that make the POD reduction fails even just for one non parametric unsteady equation. Consider the 1D scalar wave equation 
\begin{equation}\label{eq:waveEq}
\partial_t u + \partial_x u =0, \quad x\in [0,1], \, t\in [0,t_f],
\end{equation}
with initial conditions $u_0(x) = \exp(-10\sin^2(\pi(x-0.2)))$ and periodic boundary conditions, which has exact solution $u(x,t)= u_0(x-t)$. This is a simple, smooth, one parameter (time--dependent) problem. Nevertheless, if we perform a POD with tolerance $10^{-5}$ on a set of these solutions computed on a space grid with $\NFE=1000$ equispaced intervals and with $K=1400$ time steps, we obtain 13 modes to describe a one-parameter problem. In \cref{fig:advection_wave_no_calib} we plot some snapshots and in \cref{fig:POD_advection_wave_no_calib} the POD basis functions.

\begin{figure}[ht]\begin{center}
\begin{subfigure}[Snapshots at different times\label{fig:advection_shock_no_calib}]{\includegraphics[width=0.45\textwidth, trim={30 30 30 20},clip]{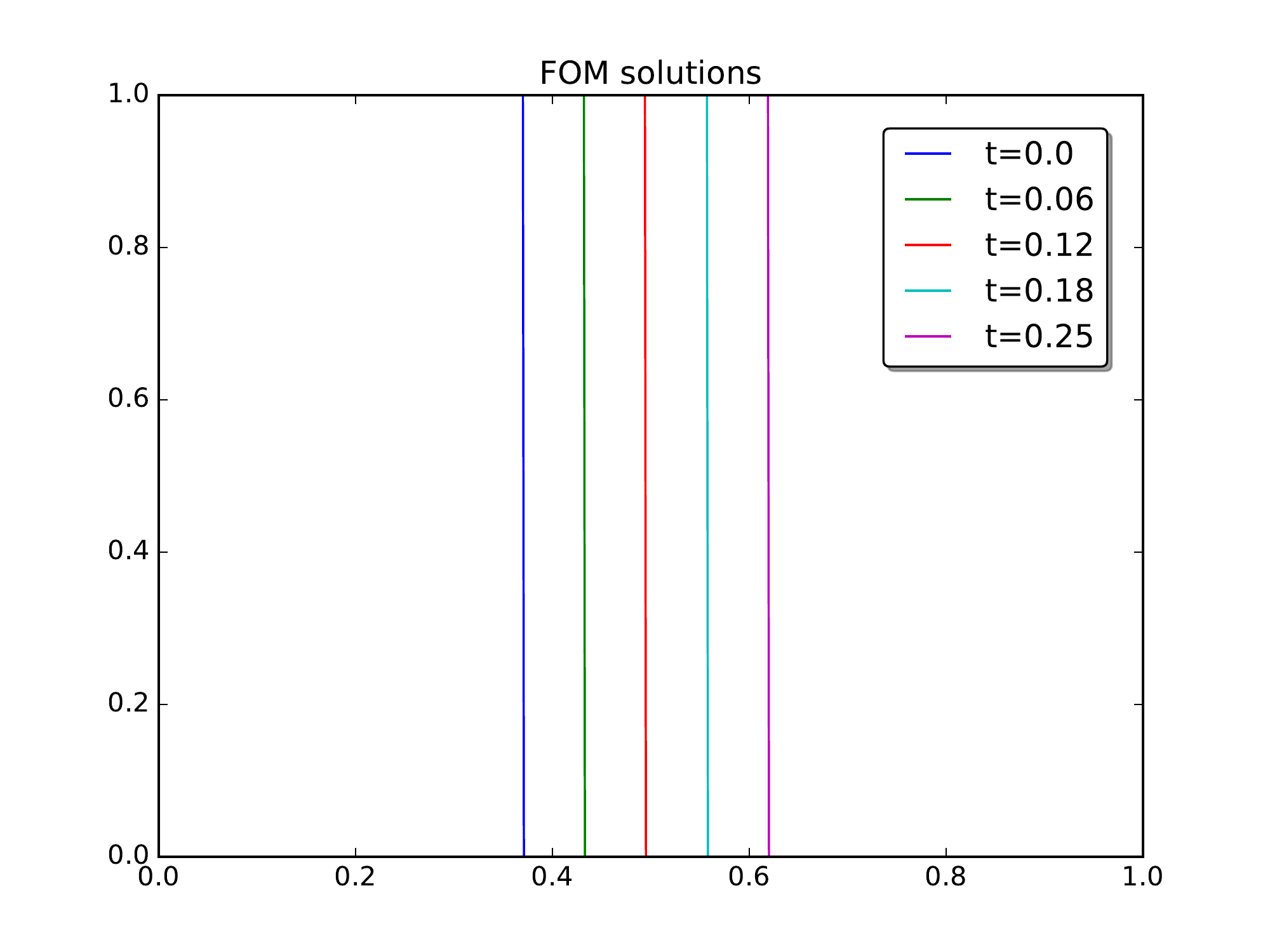}}
\end{subfigure}
\begin{subfigure}[POD modes\label{fig:POD_advection_shock_no_calib}]{\includegraphics[width=0.45\textwidth, trim={30 30 30 20},clip]{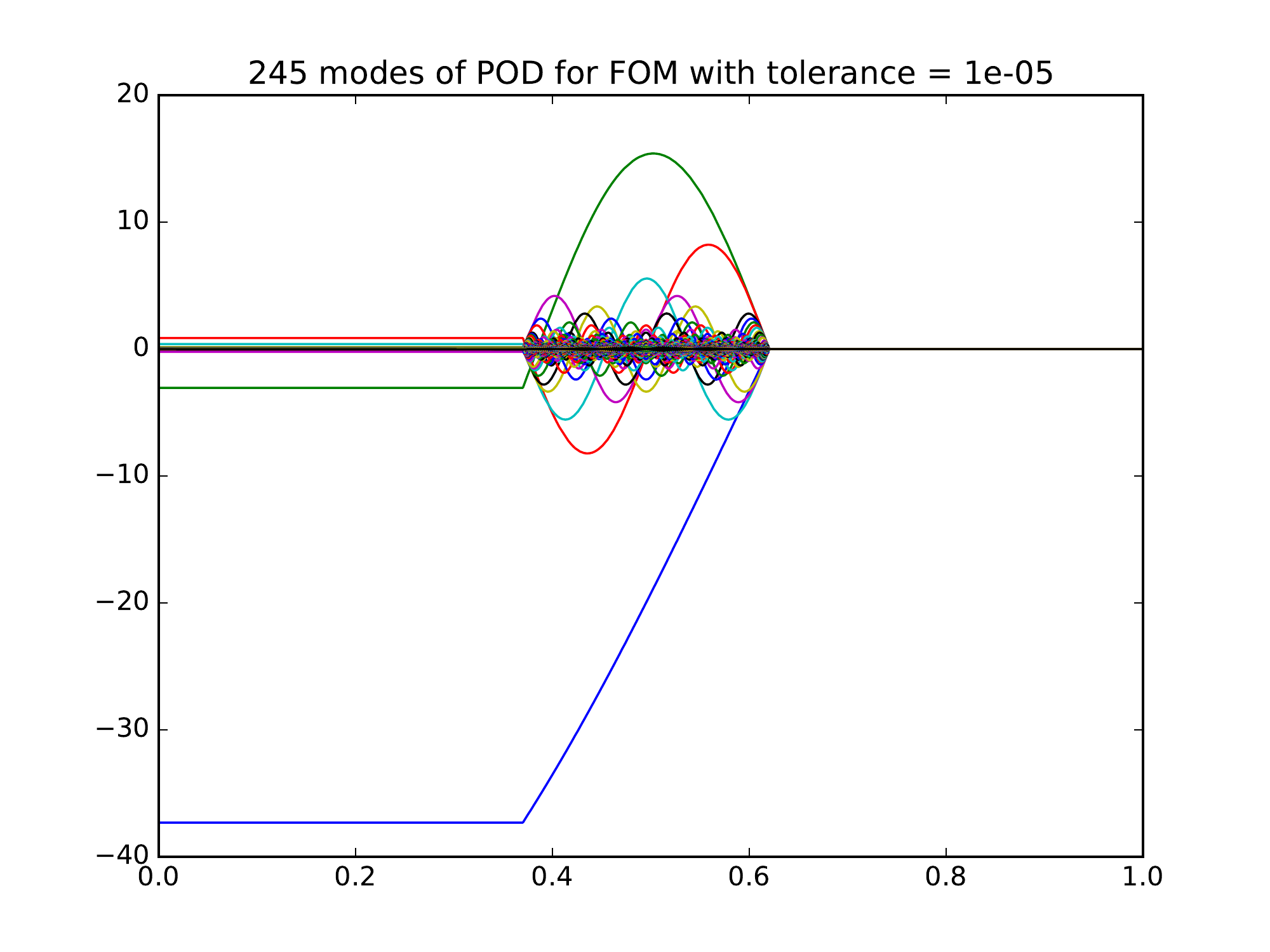}}
\end{subfigure}\\
\begin{subfigure}[Snapshots at different times, Rusanov scheme\label{fig:advection_shock_no_calib_LxF}]{\includegraphics[width=0.45\textwidth, trim={30 30 30 20},clip]{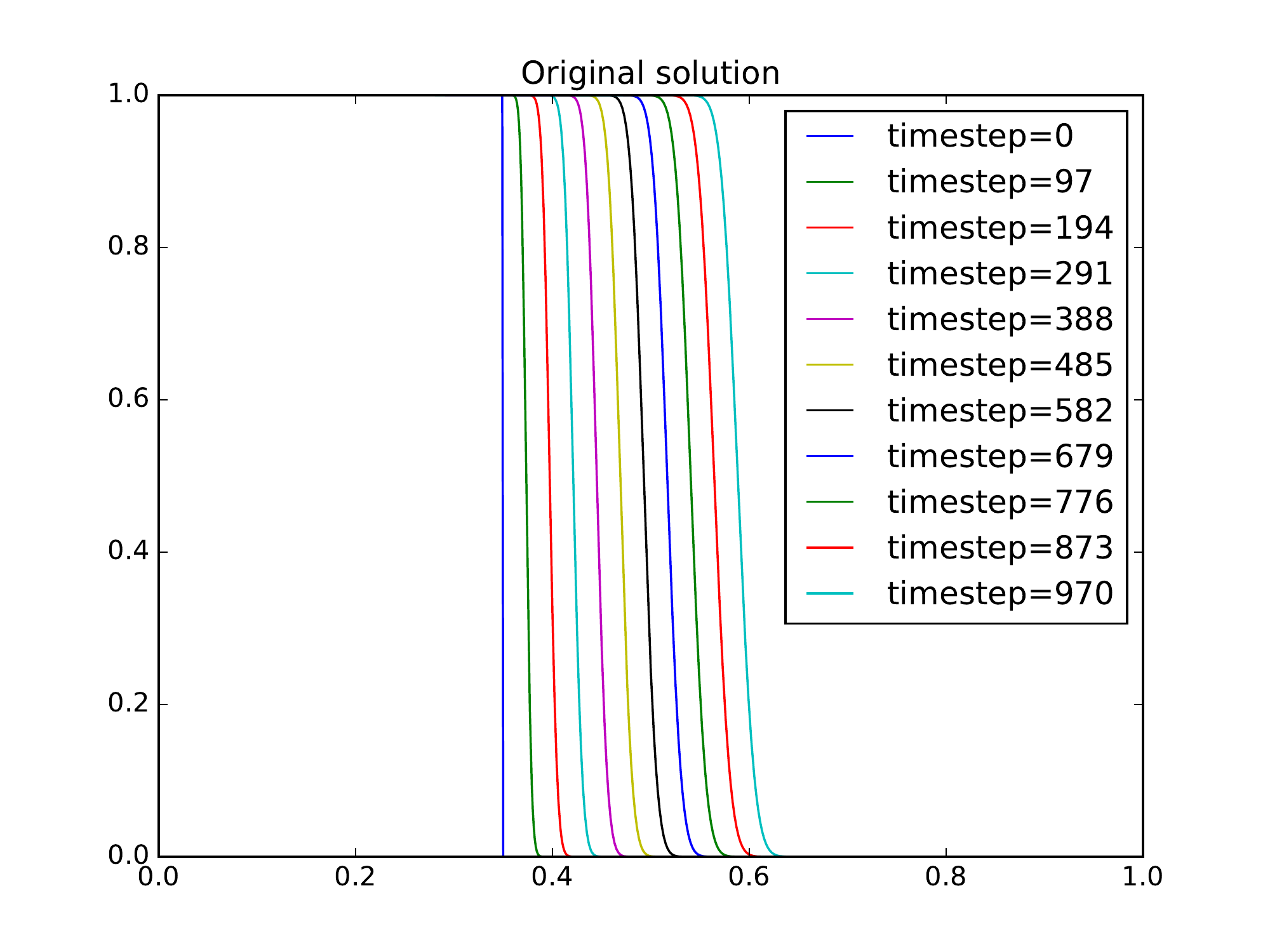}}
\end{subfigure}
\begin{subfigure}[POD modes on the Rusanov scheme\label{fig:POD_advection_shock_no_calib_LxF}]{\includegraphics[width=0.45\textwidth, trim={30 30 30 20},clip]{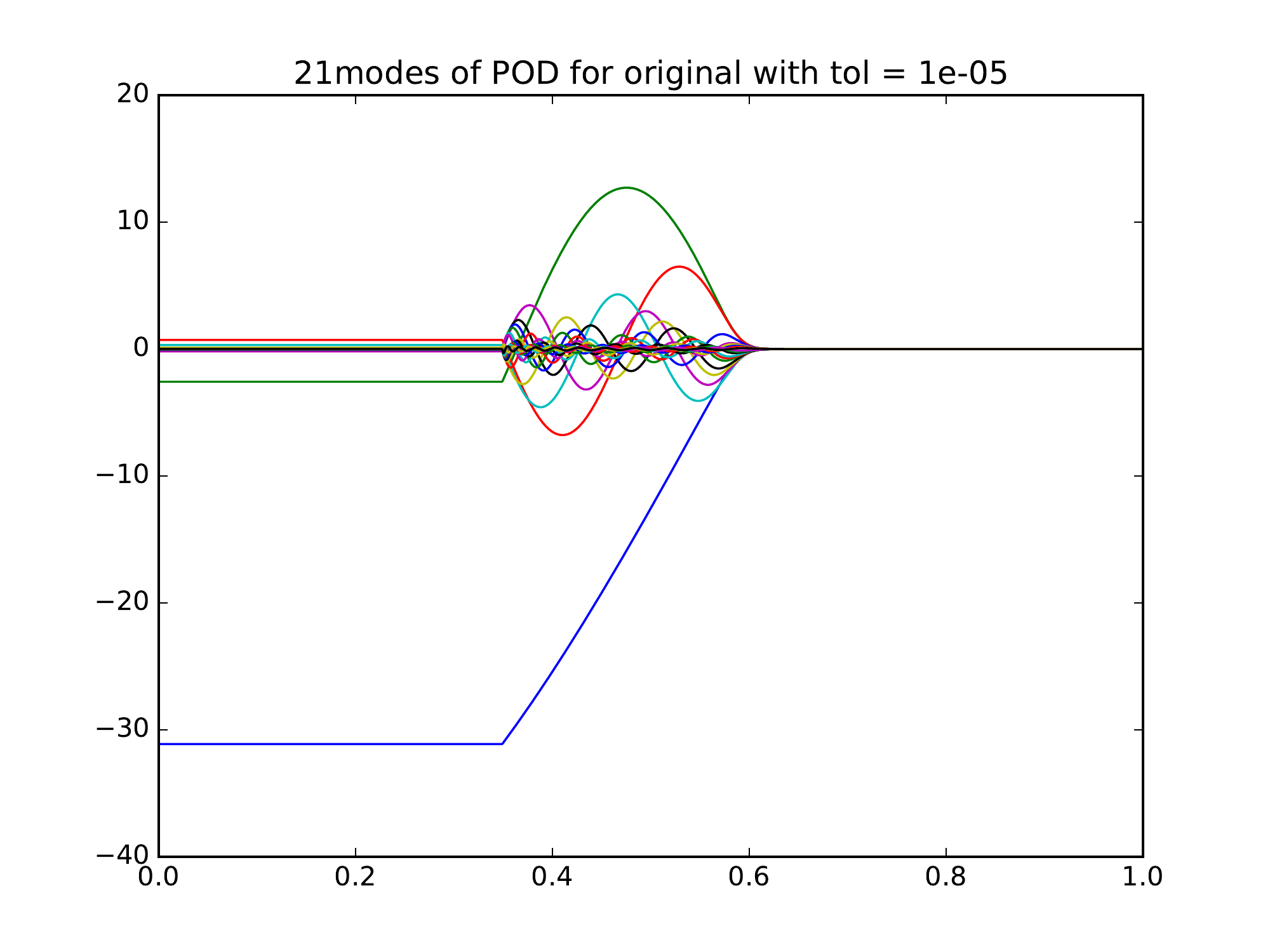}}
\end{subfigure}
\end{center}
\caption{Shock solution for wave equation}
\end{figure}

Another similar behavior can be observed with moving steep gradients or shocks, which are typical in the context of hyperbolic conservation laws, even starting with smooth initial conditions. 
Take the same wave equation \eqref{eq:waveEq} with the Riemann problem $u_0=\mathbbm{1}_{x<0.5}$ as initial condition and Dirichlet boundary conditions. 
The exact solution is a moving shock, as in \cref{fig:advection_shock_no_calib}, and, performing a POD, we get 245 modes to catch the evolution of this one-parameter problem, see \cref{fig:POD_advection_shock_no_calib}.
This number is proportional to the  size of the cells where the shock is located in any of the snapshots. This is very dangerous as soon as we increase the computational complexity.
We can observe in \cref{fig:advection_shock_no_calib_LxF} and in \cref{fig:POD_advection_shock_no_calib_LxF} that, as soon as we use a bit of diffusion in the solution scheme, in this case a Rusanov scheme, the problem becomes easier to be compressed. 
Still, 21 modes to compress information of a one-parameter problem are not what we would like to obtain.

More challenging examples can be constructed with systems of equations or nonlinear problems, where more shocks or waves moving at different speeds can meet or separate. As an example, we show the exact solution of a shallow water system in 1D, for a Riemann problem.
\begin{figure}[ht]\begin{center}
\begin{subfigure}[Snapshots at different times\label{fig:SW_no_calib}]{\includegraphics[width=0.45\textwidth, trim={30 30 30 20},clip]{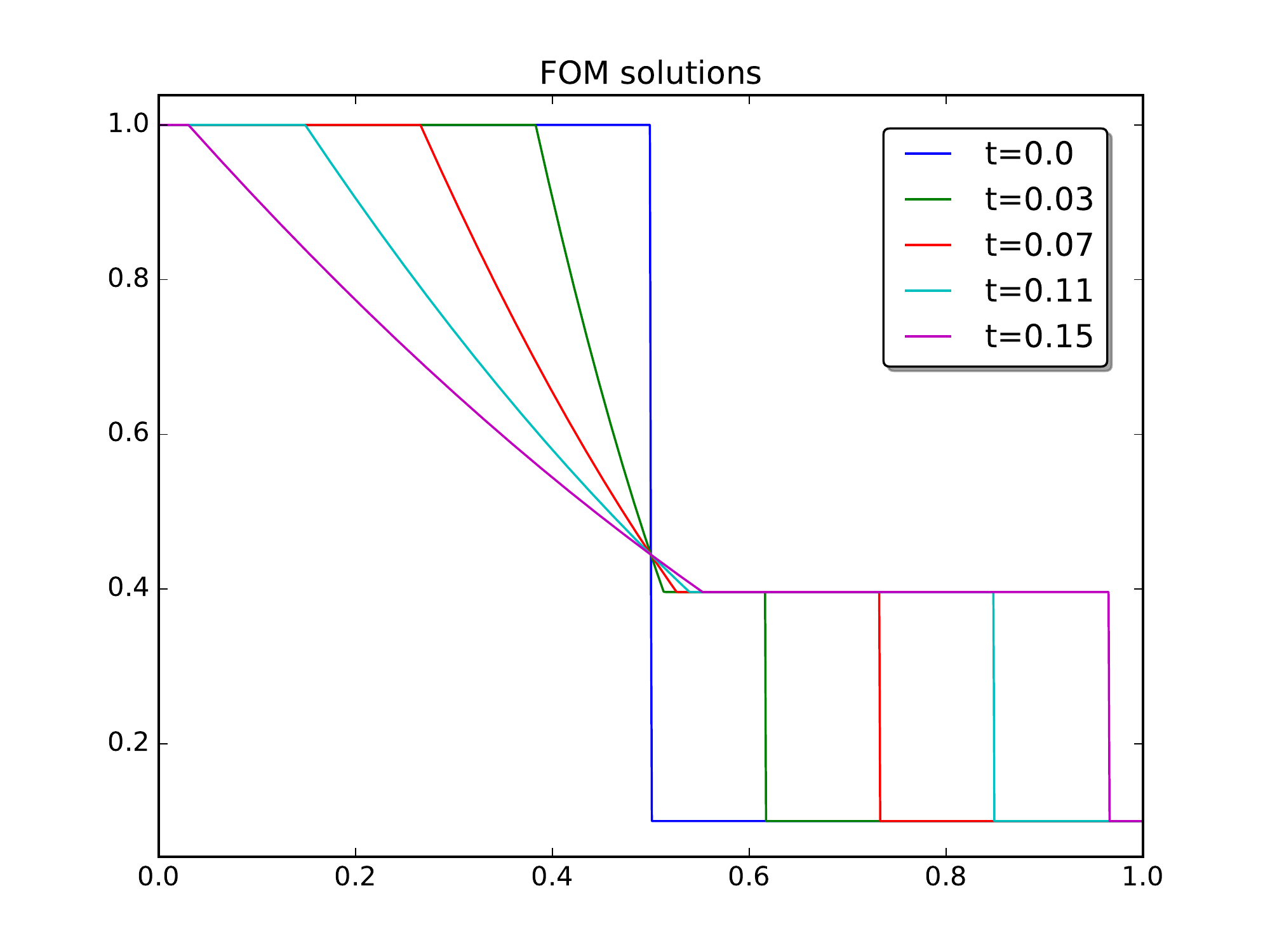}}
\end{subfigure}
\begin{subfigure}[POD modes\label{fig:POD_SW_no_calib}]{\includegraphics[width=0.45\textwidth, trim={30 30 30 20},clip]{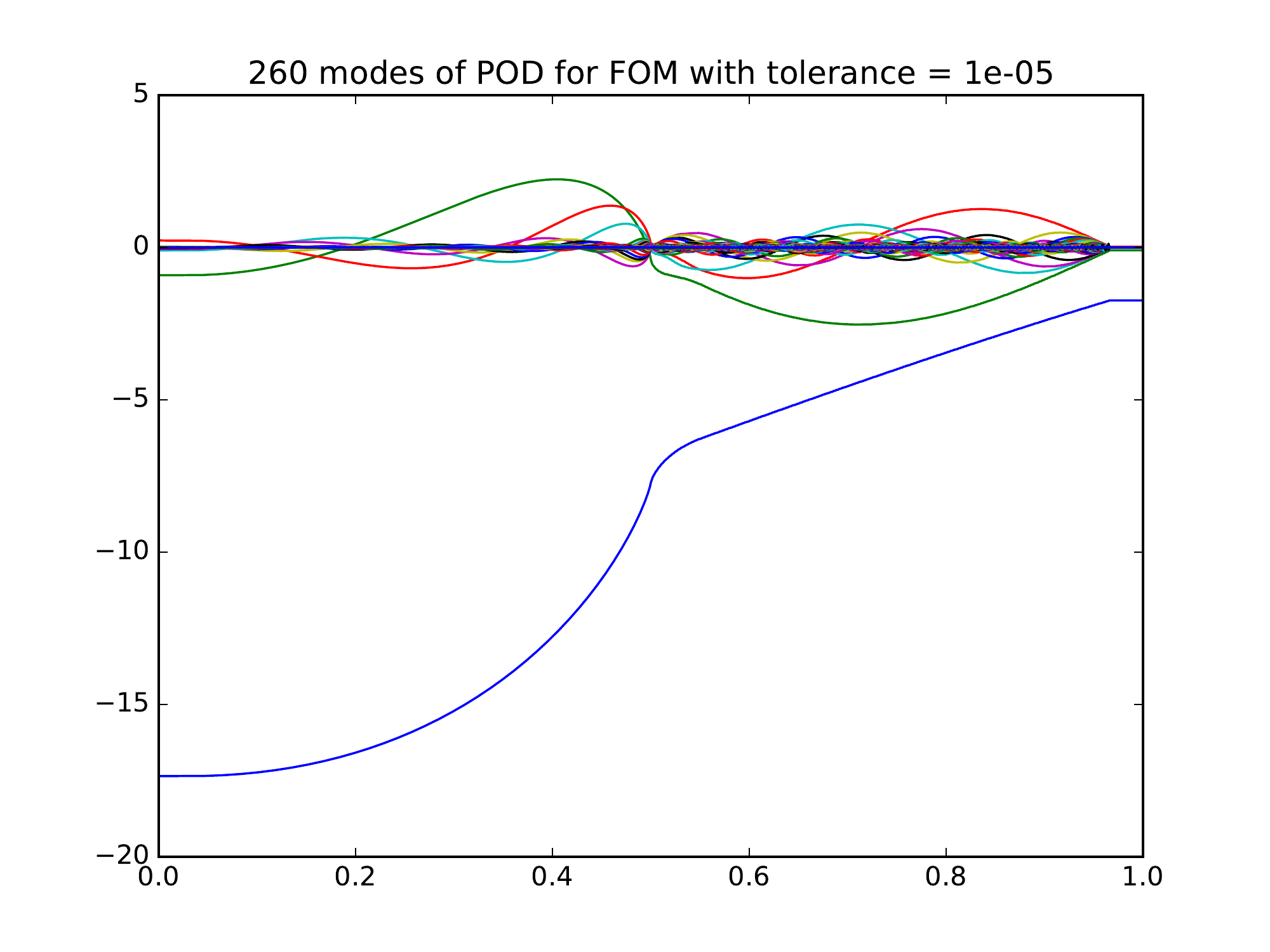}}
\end{subfigure}
\end{center}
\caption{Shallow water equations}
\end{figure} 
We can see in \cref{fig:SW_no_calib} and in \cref{fig:POD_SW_no_calib} again similar results. The dealing of multiple waves problems is more problematic than the previous ones and we will not consider them in this work. In future research they will be studied in more details. 

The common problematic of these example is the slow decay of the Kolmogorov $N$-width. It is defined as
\begin{equation}\label{def:kolmogorov}
d_N(\mathcal{S},\VFE):= \inf\limits_{\VRB \subset \VFE} \sup\limits_{f \in \mathcal{S}} \inf\limits_{g\in \VRB} ||f-g||,
\end{equation}
where $\mathcal{S}\subset \VFE$ is the manifold that we want to represent, namely, the set of all the solutions. The first infimum is taken over all the \textit{linear} subspaces of $\VFE$ with dimension $\NRB$. In other words, the Kolmogorov $N$-width is the worst error we can make given the best linear subspace $\VRB\subset \VFE$ for a fixed dimension $\NRB$. If this error $d_N(\mathcal{S},\VFE)$ is quickly decaying as $N\to \infty$, the POD will be successful. When this is not the case, linear subspaces of $\VFE$ are not the best way of compressing these information.

\section{Arbitrary Lagrangian--Eulerian Framework for MOR}\label{sec:ALE}
During the last years many works developed new non--linear techniques to approximate the manifold of the solutions $\mathcal{S}$. 
As Taddei is underlining \cite{taddei19registration}, there are, in general, two types of approaches for this problem. 
The Eulerian approaches, where the projection $\proj :\VFE \to \VRB$ may even be a nonlinear operator and the solutions are sought with different techniques, inter alia, Grassmannian learning \cite{zimmermann18onlineAdaptive}, convolutional auto-encoders \cite{carlberg18convolutionalAutoencoders}, transported/transformed snapshots \cite{sPOD, Cagniart2019}, displacement interpolation \cite{rim18displacementInterpolation} and local adaptivity \cite{peherstorfer18adaptiveBases}. 
On the other hand, there are fewer Lagrangian works where the map $\proj:\VFE \to \VRB$ is sought as a linear mapping composed with a transformation bijection, which maps the solution into a reference one, letting the transformation solve the nonlinearities of the problem. We can find some examples in \cite{iollo2014advection, mojgani17aleRB, RB_freezing, taddei19registration}. 

Moreover, some of these works make use of a transformation map to move the solutions or the reduced basis functions onto a reference domain, see \cite{sPOD, RB_freezing, trasport_greedy, iollo2014advection, Cagniart2019,taddei19registration, nair19transportedSnapshots, mojgani17aleRB}. 

What we aim to do in this work is to align the discontinuities or some features of the solution for every parameter and time step. 
To do so, we will use a transformation (bijection) of the domain into a reference one and we will rewrite the whole equation into an arbitrary Lagrangian--Eulerian (ALE) setting, where the speed of the mesh is given by the derivative in time of the transformation. 
Putting together these two features, we will be able to adapt the PODEI--Greedy algorithm to this problem and apply it directly on the reference domain, where the  decay of the Kolmogorov $N$--width will be faster.

Suppose we have a transformation map $T$ from a reference domain $\RefDom\subset \R$ to the original one $\Omega$, $T:\Theta \times \RefDom \to \Omega$, where $\Theta\subset \R^c$ is the space of the calibration parameters, that will be chosen according to time and parameter. 
Suppose that the map $T$ \textit{aligns} our solution in a way that the POD on the aligned solutions will be more effective. In practice, we define a calibration map $\theta:[0,t_f]\times \P\to \Omega$ on the spirit of \cite{Cagniart2019}, which represent a interesting feature that we want to align in all the solutions, i.\,e., $\uFE(T(\theta(t,\bmu),y),t,\bmu) \approx \bar{v}(y)$.

If we try to apply the PODEI--Greedy algorithm on the previous transformed variables, we soon realize that we do not know how to treat the transformation in the reduced equation
\begin{equation}
\begin{split}
&\sum_{i=1}^\NRB \uRBc^{k+1}(\bmu)\phiRB^i(T(\theta(t^{k+1},\bmu),y)) - \sum_{i=1}^\NRB \uRBc^{k}_i(\bmu) \phiRB^i(T(\theta(t^k,\bmu),y)) +\\
& \sum_{i=1}^\NRB \sum_{m=1}^{\NEIM} \EIMdof_m \left( \EFE (F(\uFE\left(T\left(\theta \left(t^{k+1},\bmu \right),y\right),t^{k},\bmu\right) ,\bmu)) \right) \proj_i(\EIMfun_m)\phiRB^i(T(\theta(t^k,\bmu),y))=0.
\end{split}
\end{equation}
If we keep an Eulerian approach, the bases should move with the transformation and they would depend hence on the transformation. This does not allow us to compute collocation methods, like the EIM, because it would indicate DoFs which vary for the bases functions according to transformation. 
Indeed, if we fix the EIM DoFs on the original domain $\Omega$, they will corresponds to points which have variable importance across different parameters and times and they would be meaningless for the sake of reduction. 
Many works that use Eulerian approaches do not have an \textit{online} phase because of these conflicts.

Hence, we have to adapt to an arbitrary Lagrangian--Eulerian (ALE) approach, where everything, in particular the collocation methods, is computed on the reference domain $\RefDom$. 
\subsection{Arbitrary Lagrangian--Eulerian (ALE) Framework}\label{sec:ALEframework}
Inspired by the transformations of the works of \cite{taddei19registration, RB_freezing, sPOD, Cagniart2019, Cagniart2019}, let $T:\Theta \times \RefDom \to \Omega$ be a map such that the function $T(\theta ,\cdot): \RefDom \to \Omega$ is a bijection for every $\theta \in \Theta$ and 
\begin{itemize}
\item $T(\cdot, \cdot) \in \mathcal{C}^1(\Theta\times \RefDom, \Omega ), $
\item $\exists \, T^{-1} : \Theta \times \Omega \to \RefDom$ such that $T^{-1}(\theta, T(\theta,y))=y$ for $y\in\RefDom$ and $T(\theta, T^{-1}(\theta,x))=x$ for $x \in \Omega$,
\item $T^{-1}(\cdot, \cdot) \in \mathcal{C}^1(\Theta\times \Omega , \RefDom )$.
\end{itemize}
Moreover, suppose that there exists a calibration map $\theta: \P \times [0,t_f] \to \Theta$ such that
\begin{itemize}
\item $\theta (\cdot,\bmu) \in \mathcal{C}^1([0,t_f],\Theta)$ for all $\bmu \in \P$ ,
\item $\uFE(T(\theta(t,\bmu), y),t ,\bmu) \approx \bar{v}(y) ,\quad \forall \bmu \in \P, \, t\in [0,t_f], y\in \mathcal{R},$
\end{itemize}
where the last condition expresses the way we want to align the solutions and will be explained more carefully in \cref{sec:transport_map}. There, we will also give some examples of maps that suit our typical problems.

Given this map, and a solution $\uFE(x,t,\bmu)$ of the equation \eqref{eq:cons_law}, we want to describe the behavior of the calibrated solution $\vFE(y,t,\bmu) := \uFE(T(\theta(t,\bmu),y),t,\bmu), $ through another PDE.

If we try to compute the time derivative of the calibrated solution $\vFE$, setting $x:=T(\theta(t,\bmu),y)$, we get
\begin{align}
\frac{d}{dt}\vFE(y,t,\bmu) =&  \frac{d}{dt} \uFE(T(\theta(t,\bmu),y),t,\bmu) \\
=& \partial_t \uFE(x,t,\bmu) + \partial_x \uFE(x,t,\bmu) \frac{d T(\theta(t,\bmu),y)}{dt}  \\
=& -\frac{d}{dx} F(\uFE(x,t,\bmu), \bmu) + \frac{d}{dx} \uFE(x,t,\bmu) \frac{dT(\theta(t,\bmu),y) }{dt} \\
=& -\frac{dy}{dx} \frac{d}{dy}F(\vFE(y,t,\bmu),\bmu) +\frac{dy}{dx}\frac{d}{dy} \vFE(y,t,\bmu) \frac{dT(\theta(t,\bmu),y)}{dt}.
\end{align}
So, we can write the PDE for the reference unknown $\vFE$ as
\begin{align}
\frac{d}{dt}\vFE(y,t,\bmu) +&  \frac{dT^{-1}}{dx}\frac{d}{dy} F(\vFE(y,t,\bmu),\bmu) -  \frac{dT^{-1}}{dx}\frac{d}{dy} \vFE(y,t,\bmu) 	\frac{dT(\theta(t,\bmu),y)}{dt} =0 .
\end{align}
If we remove the dependence of all the variable to simplify the notation, we obtain
\begin{align}
\frac{d}{dt}\vFE +& \frac{dT^{-1}}{dx}\frac{d}{dy} F(\vFE) -   \frac{dT^{-1}}{dx}\frac{d}{dy} \vFE	\frac{dT}{dt} =0. \label{eq:ALE}
\end{align}
Here, $ \frac{dT^{-1}}{dx}$ is the Jacobian of the inverse transformation $T^{-1}(\theta,\cdot )$ and the time derivative $\frac{dT(\theta(t,\bmu),y)}{dt} $ is also called \textit{the grid speed} in ALE context.

Now, it is clear why we need the hypotheses on the transformation and on the calibration map to be satisfied. We are using the transformation, its inverse and their derivatives in time and in space. In particular, when we compute $\frac{dT}{dt}$ we mean $\frac{dT(\theta(t,\bmu),y)}{dt}= \partial_\theta T(\theta(t,\bmu),y) \frac{d \theta(t,\bmu)}{dt}$.

The generalization to multidimensional spaces and systems of type
\begin{equation}\label{eq:conservation_law_MD}
\partial_t \uFE^s +\sum_{i =1}^d\frac{d}{dx_i}F^{s}_i(\uFE)=0, \quad \forall s=1,\dots,S,
\end{equation} is straightforward and it reads
\begin{equation}
\frac{d}{dt}\vFE^s +  \sum_{i =1}^d \sum_{j=1}^d \frac{d(T^{-1})_j}{dx_i}\frac{d}{dy_j} F^{s}_i(\vFE) -   \sum_{i =1}^d \sum_{j=1}^d \frac{d(T^{-1})_j}{dx_i}\frac{d}{dy_j} \vFE^s \frac{dT_i}{dt} =0, \quad \forall s.
\end{equation}
Nevertheless, in this work we consider only scalar and 1D problems, to avoid more technicalities that arise when more shocks or more complicated structures move at different speeds.

\subsection{MOR for ALE}
It is crucial to write the MOR algorihtm and the RB space for the reference variables $\vFE$ on the reference domain $\RefDom$. Otherwise, we would have troubles in performing the reduction and in the application of collocation methods. So, we can write
\begin{equation}\label{eq:RB_ALE}
\begin{split}
&\sum_{i=1}^\NRB (\vRBc^{k+1}-\vRBc^{k})(\bmu)\phiRB^i(y) + \sum_{i=1}^\NRB \sum_{m=1}^{\NEIM} \EIMdof_m \left( \tilde{\EFE} \left(\vRB(y),\frac{dT^{-1}}{dx}, \frac{dT}{dt} , \bmu\right) \right) \proj_i(\EIMfun_m)\phiRB^i(y)=0,
\end{split}
\end{equation}
where the new evolution operator is defined on the reference flux, following the formula \eqref{eq:ALE}, i.\,e.,
\begin{equation}\label{eq:ALE_evolution_operator}
\tilde{\EFE} \left(\vRB(y),\frac{dT^{-1}}{dx}, \frac{dT}{dt} , \bmu \right) :=
\frac{dT^{-1}}{dx}\EFE \left( F(\vRB) \right)+  \frac{dT^{-1}}{dx} \frac{dT}{dt} 
 \EFE \left( \vRB
\right) .
\end{equation}

With the ALE formulation for the RB algorithm \eqref{eq:RB_ALE}, we notice a couple of major differences with respect to the original RB formulation \eqref{eq:RB_problem}. 
First of all, we have to compute new terms regarding the transformation $\frac{dT^{-1}}{dx}$ and $\frac{dT}{dt}$, which must be easily computable, in a way not to affect the computational costs in the online phase. 
Then, the evolution scheme must be applied not only on the flux $F(\vRB)$ but also on $\vRB$ itself. 
Considering a compact stencil scheme, this can affect the computational costs of around a factor of 2.

The hope is that the reduced ALE model will need much less basis functions both in the EIM space and in the RB space. This reduction should strongly compensate the extra time we need in each computation, as we will see in the simulations of \cref{sec:results}.

\begin{oss}[Error indicator]
	The error indicator introduced in \eqref{eq:error_estimator} can be used exactly as it is in the new framework. Indeed, we need simply to substitute the evolution operator $\EFE$ with the ALE evolution operator $\tilde{\EFE}$. The considerations done in \cref{sec:error} hold for this error indicator. In particular, it is not always guaranteed to be an error bound, but it shows good behaviors in the experiments.
\end{oss}

\section{Transport Map and Learning of the Speed}\label{sec:transport_map}
\begin{figure}[h]
\begin{subfigure}[Translation\label{fig:translation}]{\includegraphics[width=0.45\textwidth, trim={30 30 20 20},clip]{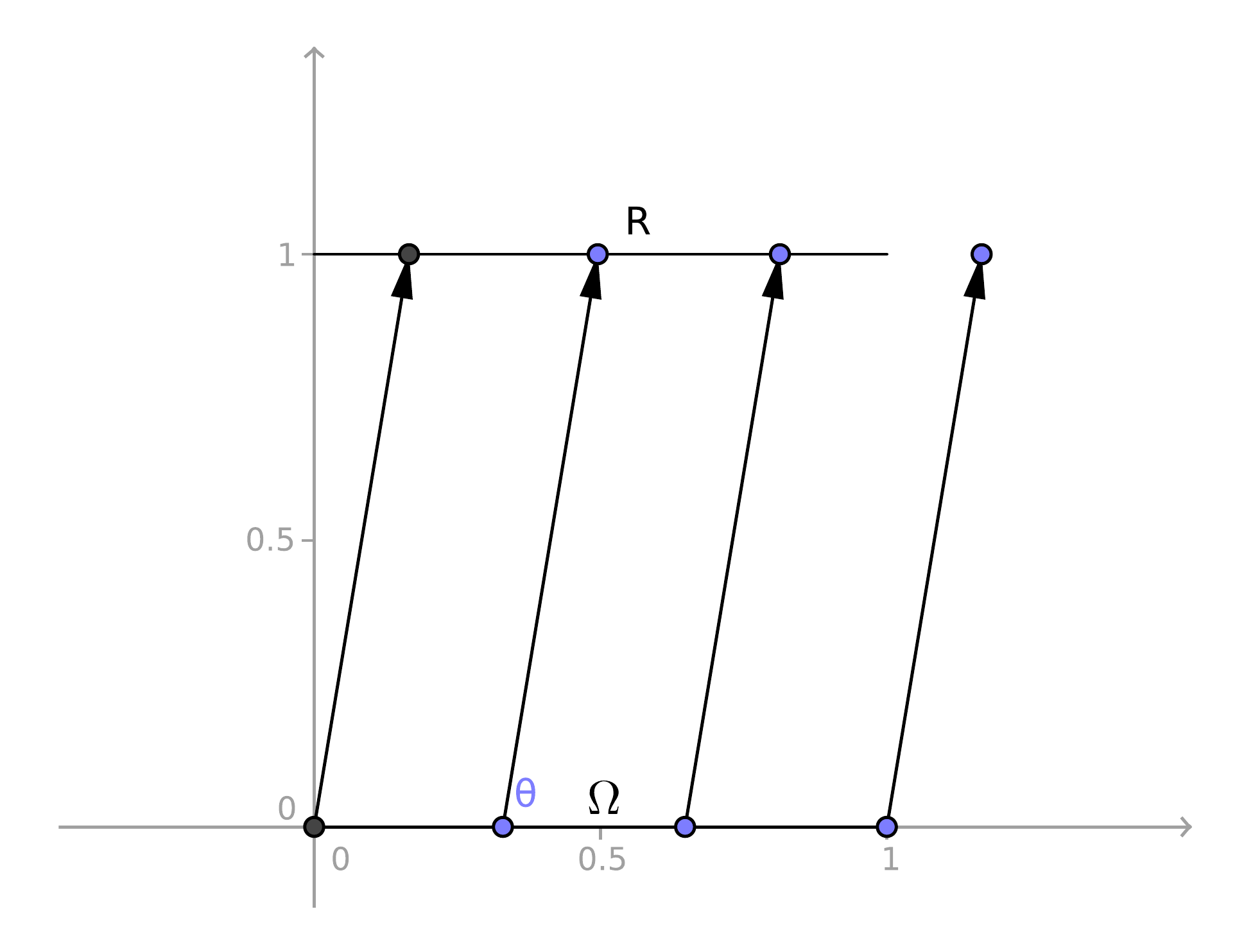}}
\end{subfigure}
\begin{subfigure}[Smooth dilatation\label{fig:dilatation}]{\includegraphics[width=0.45\textwidth, trim={30 30 20 20},clip]{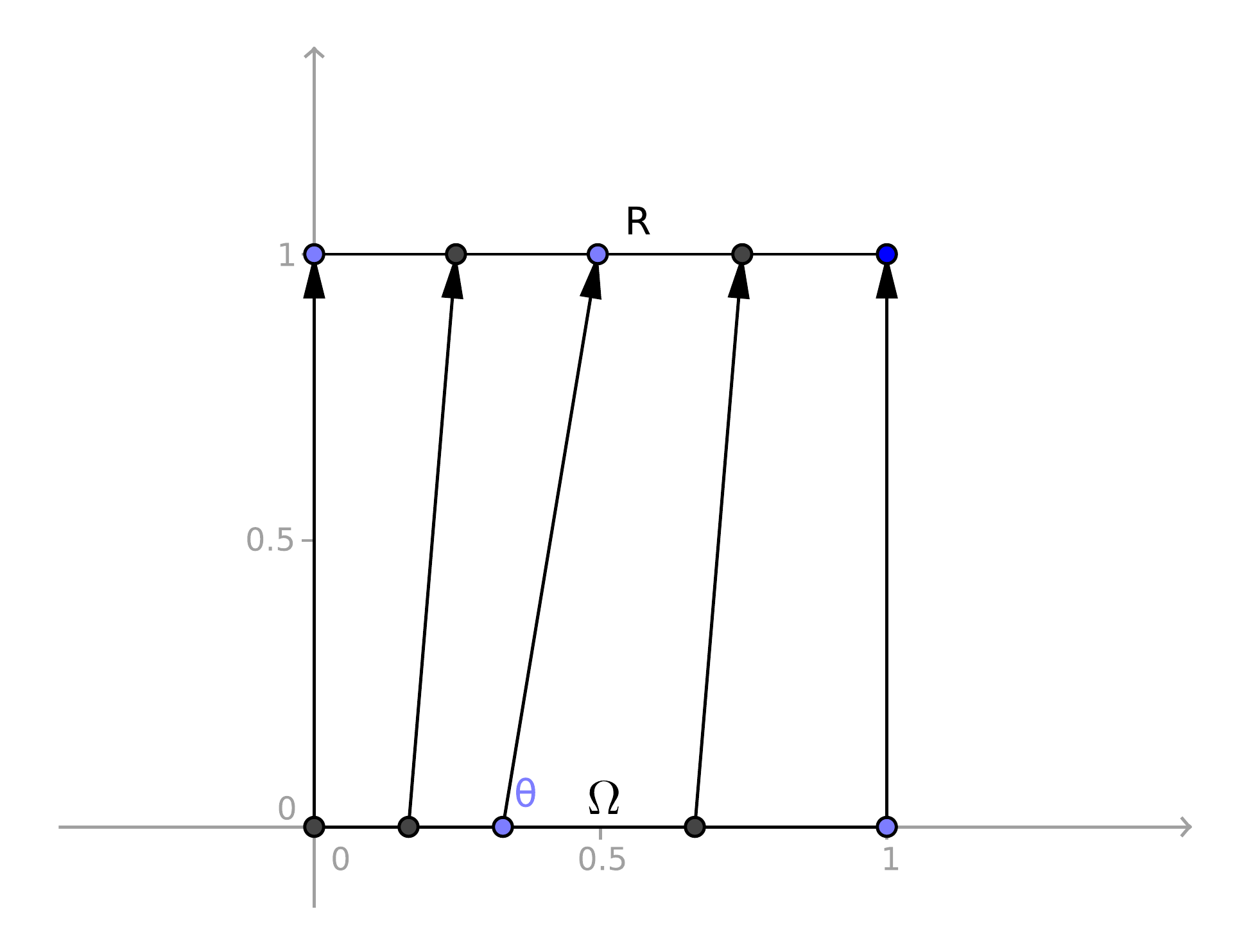}}
\end{subfigure}\caption{Examples of transport map $T^{-1}(\cdot, \theta):\Omega \to \RefDom$}
\end{figure}

To compute the reference solutions, we need a transformation map which respects the following properties:
\begin{enumerate}
\item $T$ should be written in a parameteric form $T(\theta(t,\bmu),y)$, where $T:\Theta \times \RefDom \to \Omega$, $\theta : \P   \times [0,t_f] \to \Theta$, in order to easily align solutions through calibration parameters $\theta$;
\item $T$ and $T^{-1}$ should be smooth as prescribed in \cref{sec:ALEframework};
\item $\uFE(T(\theta(t,\bmu), y),t ,\bmu) \approx \bar{v}(y) ,\, \forall \bmu \in \P, \, t\in [0,t_f],\, y\in \mathcal{R}$, in order to align the solution and gain more reduction in the RB space.
\end{enumerate} 
To address the first two point, we give a couple of examples of possible transformations that one can use. The following are the one that we will utilize in the simulations of this work.
\begin{ex}[Translation map]\label{ex:transformation_wave}
Consider the map 
\begin{equation}\label{eq:translation}
T(\theta,y) = y+\theta
\end{equation} for the traveling wave example in \cref{fig:advection_wave_no_calib}, where the domains are $\Omega=\RefDom =[ 0,1 ]$ with periodic boundary conditions, see \cref{fig:translation}.
Consider the solutions to the parametric equation $\partial_t u + \mu_1 \partial_x u=0$, where the initial conditions are $u_0(x,\bmu) = \exp(-10\sin^2(\pi(x-\mu_2)))= \mathcal{G}(x-\mu_2)$.
We know that the exact solutions are $u(x,t,\bmu)= u_0(x-\mu_1 t,\bmu)$. If we detect the maximum of each solution in $\theta(t,\bmu)=\mu_2 + \mu_1 t$ and apply the translation, we get the reference solutions 
\begin{equation}
v(y,t, \bmu)=u(T(\theta(t,\bmu),y),t,\bmu) = u(y+\mu_2 +\mu_1 t,t,\bmu)= u_0(y+\mu_2,\bmu)=\mathcal{G}(y).
\end{equation}
\end{ex}

\begin{ex}[Dilatation map]\label{ex:transformation_shock}
Consider the map $T$ and its inverse $T^{-1}$ defined by
\begin{equation}\label{eq:dilatation}
T(\theta,y) = y\frac{\theta}{(2\theta-1)y+1-\theta}, \quad T^{-1}(\theta,x) = x\frac{\theta -1}{(2\theta -1)x - \theta},
\end{equation}
for the traveling shock example in \cref{fig:advection_shock_no_calib}, where the domains are $\Omega=\RefDom =[ 0,1 ]$ with Dirichlet boundary conditions. The maps are smooth for $\theta \in (0,1)$ and $T(\theta,0)=0, \, T(\theta,1)=1$ and $T(\theta,0.5)=\theta$, see \cref{fig:dilatation}.
For this case, the exact solutions to the equation $\partial_t u + \mu_1 \partial_x u=0$, where the initial conditions are $u_0(x,\bmu) =  \mathbbm{1}_{x<\mu_2}$ are $u(x,t,\bmu)= u_0(x-\mu_1 t,\bmu)=\mathbbm{1}_{x-\mu_1 t<\mu_2}$. 
If we detect the steepest point of each solution in $\theta(t,\bmu)=\mu_2 + \mu_1 t$ and apply the transformation, we get the reference solutions 
\begin{equation}
v(y,t, \bmu)=\mathbbm{1}_{y<0.5}.
\end{equation}
\end{ex}

The alignment process, i.\,e., how we find the map $\theta(t,\bmu)$, is more challenging. 
A first possibility would be to use the information of the system to obtain a speed, in the classical ALE sense, and use this speed to compute how the transformation should behave, in order to align some features, like shocks or waves, as done in \cite{mojgani17aleRB}. 
This way does not provide a feasible method to detect the initial calibrations $\theta(t^0,\bmu)$ for different parameters $\bmu\in \P$. 
Another possibility is given by the registration procedure \cite{taddei19registration} which applies optimization techniques, given a set of parametric steady solutions.

What we will use here is a more na\"ive approach, where we detect a feature (a peak of the solution, a shock, a change in sign) and we track this feature along time and parameter domains. First, we define this map for few snapshots in a training set of Eulerian solutions $\lbrace \uFE(t,\bmu) \rbrace_{\bmu\in \P_{train}}$, then, we extend it using regression/machine learning techniques to the whole parameter and time domain as presented in \cref{algo:LearningTheta}.

\begin{algorithm}
	\fontsize{10pt}{10pt}\selectfont
	\caption{Transformation-Learning} 
	\begin{algorithmic}[1]
		{ 	\REQUIRE A training set of Eulerian snapshots: $\lbrace \uFE(t^k,\bmu): \bmu \in \P_{train}, k=0,\dots,K \rbrace $, a test set of Eulerian snapshots: $\lbrace \uFE(t^k,\bmu): \bmu \in \P_{test}, k=0,\dots,K \rbrace $, hyperparameter $s$ \\
			\STATE Detect the calibration parameter for all the training set $\lbrace \theta(t^k,\bmu): \bmu \in \P_{train}, k=0,\dots,K \rbrace $
			\STATE Test the map $\hat{\theta}$ on the test set and provide an error estimation through the test set $e \approx \max_{\bmu \in \P_{test} }\Vert\hat{\theta} -\theta\Vert$ on $[0,t_f]\times \P$
			\RETURN  Approximation map $\hat{\theta}$ and the error approximation $e$
		}
	\end{algorithmic}\label{algo:LearningTheta}
\end{algorithm}

\subsection{Learning the map $\theta$}\label{sec:learning}
Now, we present three possible regression processes that we tested in our simulations.

\subsubsection{Piecewise Linear Interpolation}
The first and most empirical method consists of a piecewise linear interpolation of the parameters of the training set (all the time steps are always spanned). 
It is straightforward to be applied when the parameters are chosen on a grid, but, according to the utilized sampling algorithm, they can be put in irregular points. In those situations, the linear interpolation becomes harder.
What we propose in this situation, given a $\bmu^*$ of which we want to compute $\hat{\theta}$, is to 
\begin{enumerate}
\item sort the parameters $\bmu \in \P_{train}$ in ascending order according to the distance $\Vert \bmu-\bmu^*\Vert$;
\item pick the first $p+1$ parameters $\lbrace \bmu^j\rbrace _{j=1}^{p+1}$;
\item write $\mu^* =\sum_{j=1}^{p+1} \alpha_j \bmu^j$, where $\sum_{j=1}^{p+1} \alpha_j =1$;
\item define $\hat{\theta} (t,\bmu^*) = \sum_{j=1}^{p+1} \alpha_j \theta (t,\bmu^j)$.
\end{enumerate}

This interpolation has some important drawbacks. First of all, the \textit{online} interpolation scales as the number of the training sample parameters, which may be many in case of low tolerance. Secondly, as the dimension of the parameter space increases, the less probable is to have a convex combination of points. This leads to instabilities in the computation of the linear combination.

A positive aspect of this interpolation that has been observed in the simulations is that few training parameters are often enough to catch a general (simple) behavior of the function. When we are facing something close to linear in the parameter space, few parameters are enough.

\subsubsection{Polynomial Regression}
A second option would be a polynomial representation of the map. This choice can also be justified by the typical examples of calibrations that we have observed in \cref{ex:transformation_wave} and in \cref{ex:transformation_shock}. Given a maximum degree $s$
\begin{equation}
\hat{\theta}(t,\bmu) = \sum_{\gamma : \Vert\gamma\Vert_{\ell^\infty} \leq s} \beta_\gamma t^{\gamma_0} \prod_{i=1}^p \mu_i^{\gamma_i},
\end{equation}
where $\gamma$ are multi--indexes of size $p+1$ and the coefficients $\beta_\gamma$ can be found through a least--square method on the training set.

Here, the hyperparameter $s$ must be carefully chosen. We can easily see that the number of parameters $\beta_\gamma $ involved in this regression are $\begin{pmatrix}
p+s+1\\p
\end{pmatrix}$. This means that the number of parameters grows exponentially with the dimension of the parameter space and the degree of the polynomials. It is really easy to ends up in overfitting phenomena when the training set is not big enough.
On the other side, a small degree $s$ may not be enough to capture the physical behavior of the calibration points.
 Even if this tuning may seem complicated, we will see that the regression given by this model is the closest to the training set. A quick hyperparameter analysis can be done on the training set before the offline phase of the MOR algorithm.

\subsubsection{Multilayer Perceptron}
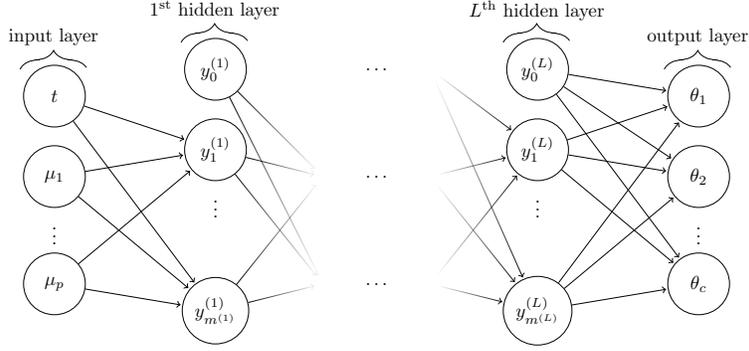
\begin{figure}[h]
	\centering
	\begin{adjustbox}{width=10cm, height=6cm, keepaspectratio}
		\begin{tikzpicture}[shorten >=1pt]
		\tikzstyle{unit}=[draw,shape=circle,minimum size=1.15cm]
		\tikzstyle{hidden}=[draw,shape=circle,minimum size=1.15cm]

		\node[unit](x0) at (0,3.5){$t$};
		\node[unit](x1) at (0,2){$\mu_1$};
		\node at (0,1){\vdots};
		\node[unit](xd) at (0,0){$\mu_p$};

		\node[hidden](h10) at (3,4){$y_0^{(1)}$};
		\node[hidden](h11) at (3,2.5){$y_1^{(1)}$};
		\node at (3,1.5){\vdots};
		\node[hidden](h1m) at (3,-0.5){$y_{m^{(1)}}^{(1)}$};

		\node(h22) at (5,0){};
		\node(h21) at (5,2){};
		\node(h20) at (5,4){};
		
		\node(d3) at (6,0){$\ldots$};
		\node(d2) at (6,2){$\ldots$};
		\node(d1) at (6,4){$\ldots$};

		\node(hL12) at (7,0){};
		\node(hL11) at (7,2){};
		\node(hL10) at (7,4){};
		
		\node[hidden](hL0) at (9,4){$y_0^{(L)}$};
		\node[hidden](hL1) at (9,2.5){$y_1^{(L)}$};
		\node at (9,1.5){\vdots};
		\node[hidden](hLm) at (9,-0.5){$y_{m^{(L)}}^{(L)}$};

		\node[unit](y1) at (12,3.5){$\theta_1$};
		\node[unit](y2) at (12,2){$\theta_2$};
		\node at (12,1){\vdots};	
		\node[unit](yc) at (12,0){$\theta_c$};

		\draw[->] (x0) -- (h11);
		\draw[->] (x0) -- (h1m);

		\draw[->] (x1) -- (h11);
		\draw[->] (x1) -- (h1m);

		\draw[->] (xd) -- (h11);
		\draw[->] (xd) -- (h1m);

		\draw[->] (hL0) -- (y1);
		\draw[->] (hL0) -- (yc);
		\draw[->] (hL0) -- (y2);

		\draw[->] (hL1) -- (y1);
		\draw[->] (hL1) -- (yc);
		\draw[->] (hL1) -- (y2);

		\draw[->] (hLm) -- (y1);
		\draw[->] (hLm) -- (y2);
		\draw[->] (hLm) -- (yc);

		\draw[->,path fading=east] (h10) -- (h21);
		\draw[->,path fading=east] (h10) -- (h22);
		
		\draw[->,path fading=east] (h11) -- (h21);
		\draw[->,path fading=east] (h11) -- (h22);
		
		\draw[->,path fading=east] (h1m) -- (h21);
		\draw[->,path fading=east] (h1m) -- (h22);
		
		\draw[->,path fading=west] (hL10) -- (hL1);
		\draw[->,path fading=west] (hL11) -- (hL1);
dd		\draw[->,path fading=west] (hL12) -- (hL1);
		
		\draw[->,path fading=west] (hL10) -- (hLm);
		\draw[->,path fading=west] (hL11) -- (hLm);
		\draw[->,path fading=west] (hL12) -- (hLm);
		
		\draw [decorate,decoration={brace,amplitude=10pt},xshift=-4pt,yshift=0pt] (-0.5,4) -- (0.75,4) node [black,midway,yshift=+0.6cm]{input layer};
		\draw [decorate,decoration={brace,amplitude=10pt},xshift=-4pt,yshift=0pt] (2.5,4.5) -- (3.75,4.5) node [black,midway,yshift=+0.6cm]{$1^{\text{st}}$ hidden layer};
		\draw [decorate,decoration={brace,amplitude=10pt},xshift=-4pt,yshift=0pt] (8.5,4.5) -- (9.75,4.5) node [black,midway,yshift=+0.6cm]{$L^{\text{th}}$ hidden layer};
		\draw [decorate,decoration={brace,amplitude=10pt},xshift=-4pt,yshift=0pt] (11.5,4) -- (12.75,4) node [black,midway,yshift=+0.6cm]{output layer};
	\end{tikzpicture}
	
\end{adjustbox}
\caption{Multilayer Perceptron architecture}\label{fig:MLper}
\end{figure}

In this section we describe the specifications for one artificial neural network (ANN) that can be used to learn the calibraton map. More details about ANN can be found in \cite{Goodfellow16deepLearning}. 
A multilayer perceptron is an ANN composed of several layers: an input layer $(t,\bmu)\in \mathcal{Y}^{(0)} = [0,t_f] \times \P$, where we pass the parameters of our problem, $L$ hidden layers $y^{(k)}\in \mathcal{Y}^{(k)} \subset \R^{m^{(k)}}$ for $k=1,\dots, L$ and an output layer $\theta \in \Theta$, where we receive the prediction of the calibration parameter. 
In \cref{fig:MLper} one can observe the architecture of this ANN. Each layer is connected to the following and the previous ones through \textit{weights}, i.e.,  affine maps $\delta^{(k)}:\mathcal{Y}^{(k)}\to \mathcal{Y}^{(k+1)}$, which are represented by arrows in \cref{fig:MLper}. 
In every node of the hidden and output layers a nonlinear activation function $\zeta:\R\to \R$ is performed component--wise. We denote with $\tilde{\zeta}$ the component--wise extension of the scalar map $\zeta$.
Overall, the multilayer perceptron map is defined as
\begin{equation}
\hat{\theta}(t,\bmu)  = \tilde{\zeta}(\delta^{(L)} (\tilde{\zeta} ( \dots (\delta^{(0)}(t,\bmu))))).
\end{equation} 

Using a training and a validation set, the learning process changes the weights $\delta^{(k)}$ based on the error of the output with respect to the exact values. The supervised learning is carried out through the backpropagation of the error combined with a stochastic gradient descent algorithm, which is an extension of the least mean squares method. Details about this algorithm can be found in \cite{Goodfellow16deepLearning}. We make use of the Keras package \cite{chollet2015keras} in Python to build and learn the ANN. 

We choose $\zeta(x):=\tanh (x)$ as activation function, because we are looking for a smooth transformation.
After different tests, where we have not observed large variation in the results, we choose to have $L=4$ hidden layers with $8$ nodes each. 

The stochastic gradient descent prevents overfitting phenomena even if the number of parameters $\sum_{k=0}^L (m^{(k)}+1)m^{(k+1)}$ can be large. To obtain a reasonable validation error for this algorithm, we have to use a large training set, as we will see in the simulations. This is probably due to the fact that the shape of the function is far away from the exact one that we want to reach.

After a hyperanalysis on the number of nodes and layers, we fix them for all the tests to 4 hidden layers and 8 nodes each layer. 
The hyperanalysis shows that more layers/nodes are too difficult to be trained by the samples that we can produce with the full solver and fewer layers/nodes are producing worse approximation results.

Improvements on this algorithms and on other learning algorithms are anyway under investigation.

\subsection{Final Algorithm}\label{sec:total_algo}
To complete the algorithm, we glue together the different pieces as specified in \cref{algo:ALE-PODEI-Greedy-Offline}. 
First of all, we compute the Eulerian solutions of the training and validation sets $\P_{train}, \P_{valid}$. Using the calibration procedure, we obtain the calibration points for these sets. 
As remarked in \cref{sec:learning}, we run different training processes on the different tests to obtain optimal hyperparameters of the methods. 
In particular, we will run different tests with different methods to compare the different results.
We check on the validation set that the error of the calibration process is smaller than a tolerance (something related to the discretization scale, we chose, for example, $5\Delta x$). 
Thus, we use the approximated calibration map $\hat{\theta}$ to compute the PODEIM--Greedy algorithm on the ALE solutions.

\begin{algorithm}
	\fontsize{10pt}{10pt}\selectfont
	\caption{ALE PODEI Greedy with learning -- Offline Phase} 
	\begin{algorithmic}[1]
		{ 	\REQUIRE Eulerian solver, Lagrangian solver, a training set $\P_{train}$, a test set $\P_{valid}$, hyperparameter for learning $s$, $\epsilon^{tol}_\theta$, $\epsilon^{tol}_{\RB}$, $\NRB_{max}$ \\
			\STATE Generate the training and test set with the Eulerian solver $\mathcal{M}_{train}=\lbrace \uFE(t^k,\bmu) : \bmu \in \P_{train}, \, \forall k \rbrace $ and $\lbrace \uFE(t^k,\bmu) : \bmu \in \P_{valid}, \,\forall k \rbrace $ 
			\STATE Compute the calibration points $\lbrace \theta(t^k,\bmu) : \bmu \in \P_{train}, \, \forall k \rbrace $ and $\lbrace \theta(t^k,\bmu) : \bmu \in \P_{valid}, \, \forall k \rbrace $
			\STATE $\hat{\theta}$, $\texttt{err}_\theta$=Transformation learning($\lbrace \theta(t^k,\bmu) : \bmu \in \P_{train}, \, \forall k \rbrace $, $\lbrace \theta(t^k,\bmu) : \bmu \in \P_{test}, \, \forall k \rbrace $, $s$) 
			\STATE Check that $\texttt{err}_\theta<5 \Delta x$
			\STATE $\RB, \EIM$ = PODEIM--Greedy($ \P_{train}$, $\epsilon^{tol}$, $\NRB_{max}$, $\tilde{\EFE}$, $\hat{\theta}$) with the ALE flux  $\tilde{\EFE}$.
			\RETURN  $\RB, \EIM, \hat{\theta}$
		}
	\end{algorithmic}\label{algo:ALE-PODEI-Greedy-Offline}
\end{algorithm}
\subsection{Online Phase}
The \textit{online} of the ALE PODEI-Greedy algorithm is given by the simple formula \eqref{eq:RB_ALE}, where the reduced evolution operator is computed with the EIM algorithm as in \eqref{eq:ALE_evolution_operator} and the total number of flux evaluation computed are at most $2\NEIM \cdot \NRB$ instead of $\NFE$.
Moreover, the evaluation of the reduced evolution operator necessitates of the map $\hat{\theta}$ in the ALE framework. The final solution on the original domain can be quickly recomputed through the maps $T^{-1}$.
The computation costs of the calibration map and of the reconstruction are negligible with respect to the computation of the solution.

With the ALE strategy we wish to strongly decrease the dimensions $\NEIM $ and $ \NRB$, in order to gain computational advantages in the \textit{online} phase. This is what we show in the next section of simulations.

\section{Results}\label{sec:results}
In this section we present some tests and their hyperanalysis.
We will study just simple scalar 1D hyperbolic problems, but extension to systems and 2D problems is straightforward when we still deal with one speed features. More complicated structures and more shocks will be object of future works.

We will consider the linear advection equation, the Burgers' equation and the Buckley--Leverett equation, respectively,
\begin{align} \label{eq:linAdv}
&\partial_t u + a \partial_x u =0,\\
&\partial_t u + a \partial_x \frac{u^2}{2} =0, \label{eq:burgers}\\
&\partial_t u + \partial_x \frac{u^2}{u^2+a(1-u^2)}=0. \label{eq:buckley}
\end{align}

In order to run all the simulations with the same time steps to facilitate the online phase, we fix \begin{equation}
	\Delta t : = \text{CFL} \min_{\bmu \in \P_{train}, x\in \Omega} \frac{\Delta x}{|J_uF(\uFE(x,0,\bmu),\bmu)|},
\end{equation} with CFL=0.25. This guarantees us a reasonable security of not incurring into oscillations. For all the simulations we used 1000 nodal points in the domain $\Omega$. In all the algorithms we use a very stable Rusanov scheme as space discretization and forward Euler in time.

We proceed with the different cases showing first of all the training process of the regression of the calibration maps, plotting the error on the validation set with respect to the dimension $N_{train}$ of the training set. 
Afterwards, we motivate the choice of the regressions we use in the \textit{offline} phase. 

For the \textit{offline }phases, we show the plot of the error decay of the PODEI--Greedy process only for few tests but we store all the resulting data in \cref{tab:RBdata}. In these plots we show the maximum of the error indicator for all the parameters, the maximum of the error and the average error with respect to the dimension of the RB space. On the error indicator we print the dimension of the EIM space. In a couple of simulations we also plot the EIM error decay with respect to the dimension of the EIM space. 
Sometimes the error goes up as the dimension increases, this means that a new solution has been considered for the extension of the EIM space, see \cref{algo:PODEIM_greedy,algo:EIM}.

Then, we plot some simulations of the \textit{online} phase of both Eulerian and ALE approaches for one parameter in the range and few time steps.

\subsection{Advection of Solitary Wave}\label{sec:advSol} 
In this test, we consider a solitary wave with different amplitudes and centers traveling at different speeds. The domain is $\Omega=[0,1]$ and the initial conditions are 
\begin{equation}\label{eq:advSol}
u_0(x,\bmu)=e^{-(100+500\mu_1)(x-0.2+0.1\mu_2)^2}, \qquad \mu_1,\mu_2 \in [-1,1],
\end{equation}
and we solve the linear transport equation \eqref{eq:linAdv} with $a=\mu_0 \in [0,2]$. We use periodic boundary conditions and the translation \eqref{eq:translation} as calibration map. The calibration point is chosen to be the maximum value point.
\begin{table}
	\begin{center}
		\begin{tabular}{c||c|c|c|c|c|c}
			Test & Dim RB & Dim EIM & FOM time\footnote[1]{The computations are performed with a Intel(R) Xeon(R) CPU E7-2850  @ 2.00GHz.} & RB time\footnotemark[1] & Ratio & Online error\footnote[2]{This error refers to the simulation run for one parameter.}\\ \hline
			\eqref{eq:advSol} ALE Poly2 & 4 & 7 & 516s & 18s & 3 \% & $5.2\, 10^{-4}$  \\ 
			\eqref{eq:advSol} ALE ANN & 12 & 20 & 516s & 38s & 7 \% & $1.7\, 10^{-4}$ \\ 
			\eqref{eq:advSol} Eulerian & 52 & 54 & 191s & 24s & 12 \% & $2.4\, 10^{-4}$  \\ \hline
			\eqref{eq:advShock} ALE Poly2 & 17 & 22 & 125s & 6s & 5 \% & $7.6\, 10^{-5}$ \\ 
			\eqref{eq:advShock} Eulerian & 64 & 124 & 49s & 9s & 18 \% & $5.3\, 10^{-4}$ \\ \hline
			\eqref{eq:BurOsc} ALE Poly3 & 50 & 60 & 314s & 35s & 11 \% & $2.9\, 10^{-4}$ \\ 
			\eqref{eq:BurOsc} ALE ANN & 139 & 167 & 298s & 66s & 22 \% & $6.4\, 10^{-4}$\\ 
			\eqref{eq:BurOsc} Eulerian\footnote[3]{For this test the algorithm does not reach the requested tolerance $10^{-3}$.} & 153\footnotemark[3] & 335\footnotemark[3] & 119s & 50s\footnotemark[3] & 42\%\footnotemark[3] & $1.2\, 10^{-3}$  \\ \hline
			\eqref{eq:BurSin} ALE Poly4 & 19 & 41 & 444s & 53s & 11\% & $3.8\, 10^{-4}$ \\ 
			\eqref{eq:BurSin} Eulerian & failed & $>600$ & 167s & $\infty$ & $\infty$ &$\infty$ \\ \hline
			\eqref{eq:BuckIC} ALE pwL & 25 & 45 & 462s & 79s & 17\% & $5.5\, 10^{-4}$ \\ 
			\eqref{eq:BuckIC} Eulerian\footnotemark[3] & 16\footnotemark[3] & 270\footnotemark[3] & 190s & 69s\footnotemark[3] & 36\%\footnotemark[3] & $9.2\, 10^{-3}$ \\ \hline
		\end{tabular}
	\end{center}\caption{Times and dimensions of the tests}\label{tab:RBdata}
\end{table}

\begin{figure}[h]
	\begin{center}
\begin{subfigure}[Offline error in ALE framework with Poly2\label{fig:test0:offlineALE}]{\includegraphics[width=0.45\textwidth, trim={0 0 20 40},clip]{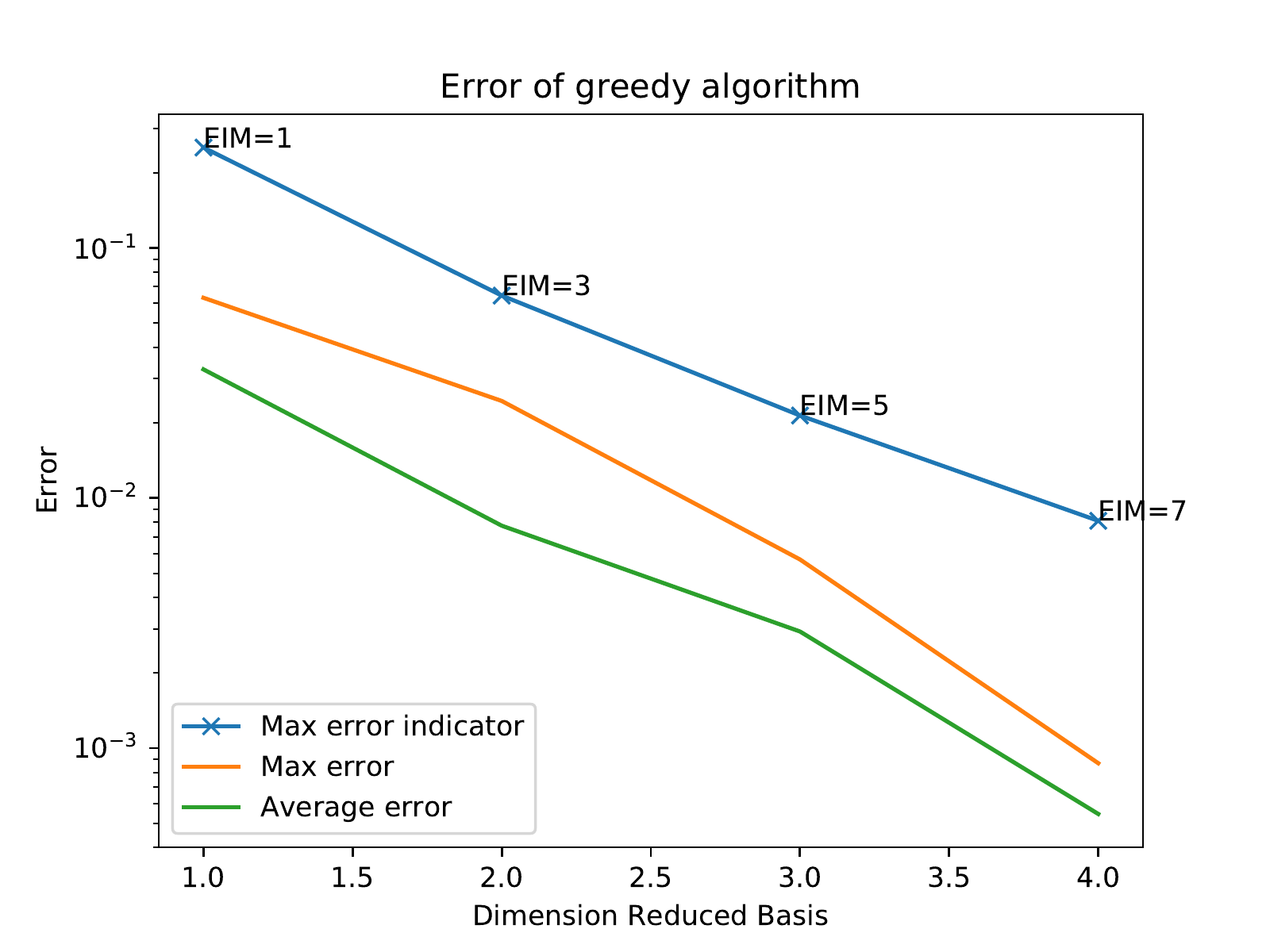}}
\end{subfigure}
\begin{subfigure}[Training of calibration\label{fig:test0:training}]{\includegraphics[width=0.45\textwidth, trim={0 0 20 40},clip]{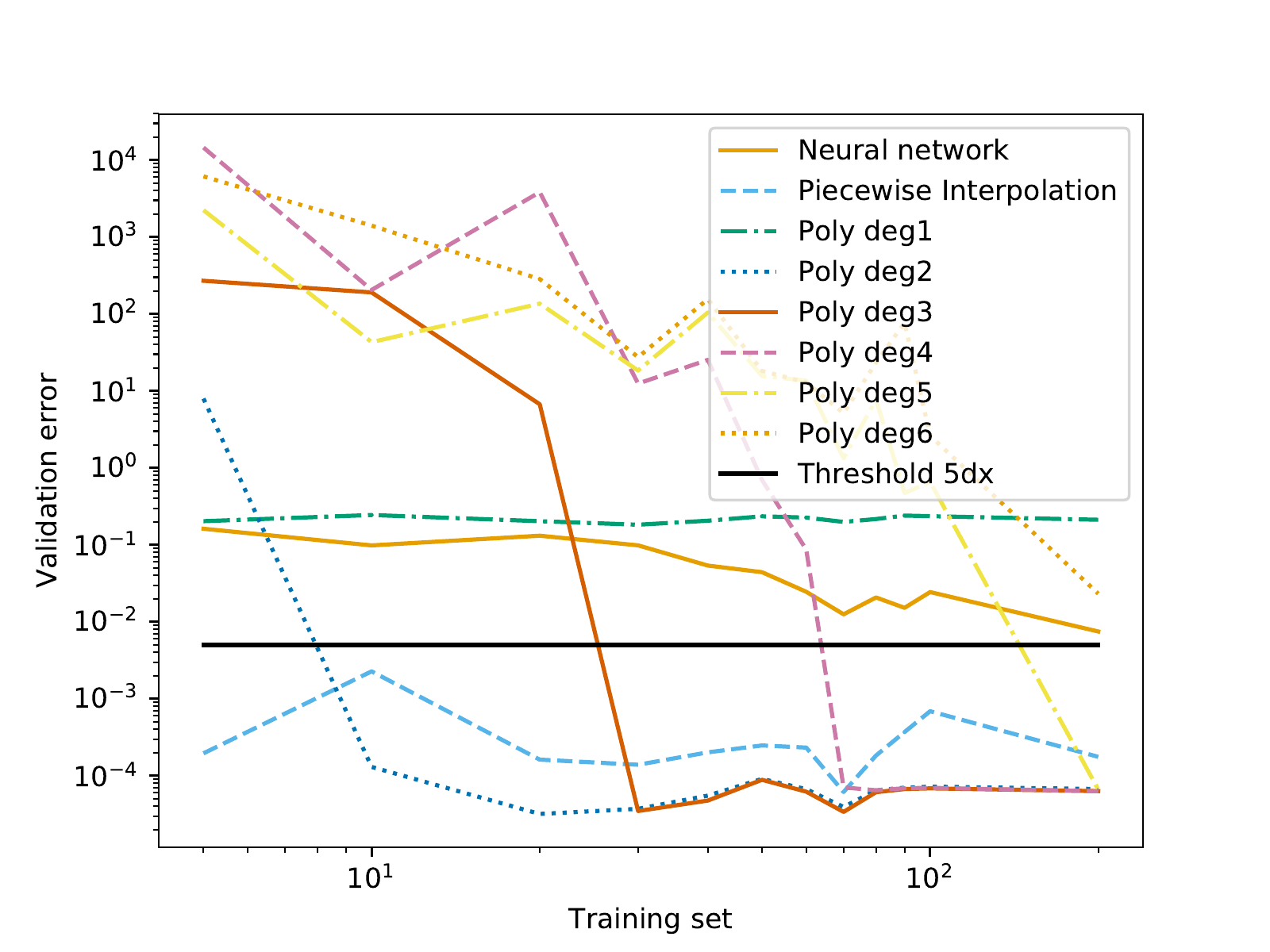}}
\end{subfigure}
\\
\begin{subfigure}[Offline error in Eulerian framework\label{fig:test0:offlineEu}]{\includegraphics[width=0.45\textwidth, trim={0 0 20 39},clip]{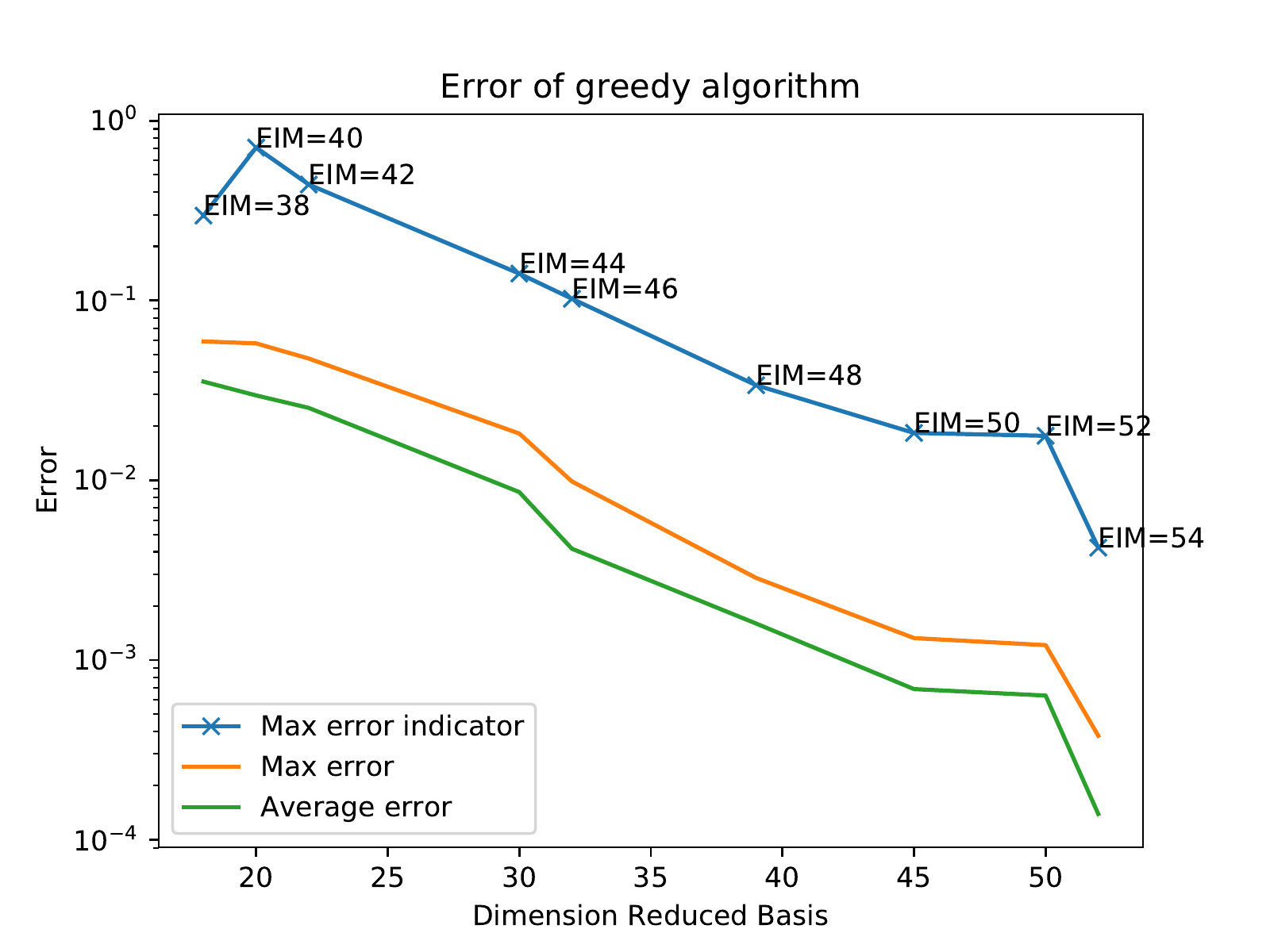}}
\end{subfigure}
\begin{subfigure}[EIM offline error in Eulerian framework\label{fig:test0:EIMEu}]{\includegraphics[width=0.45\textwidth, trim={0 0 20 40},clip]{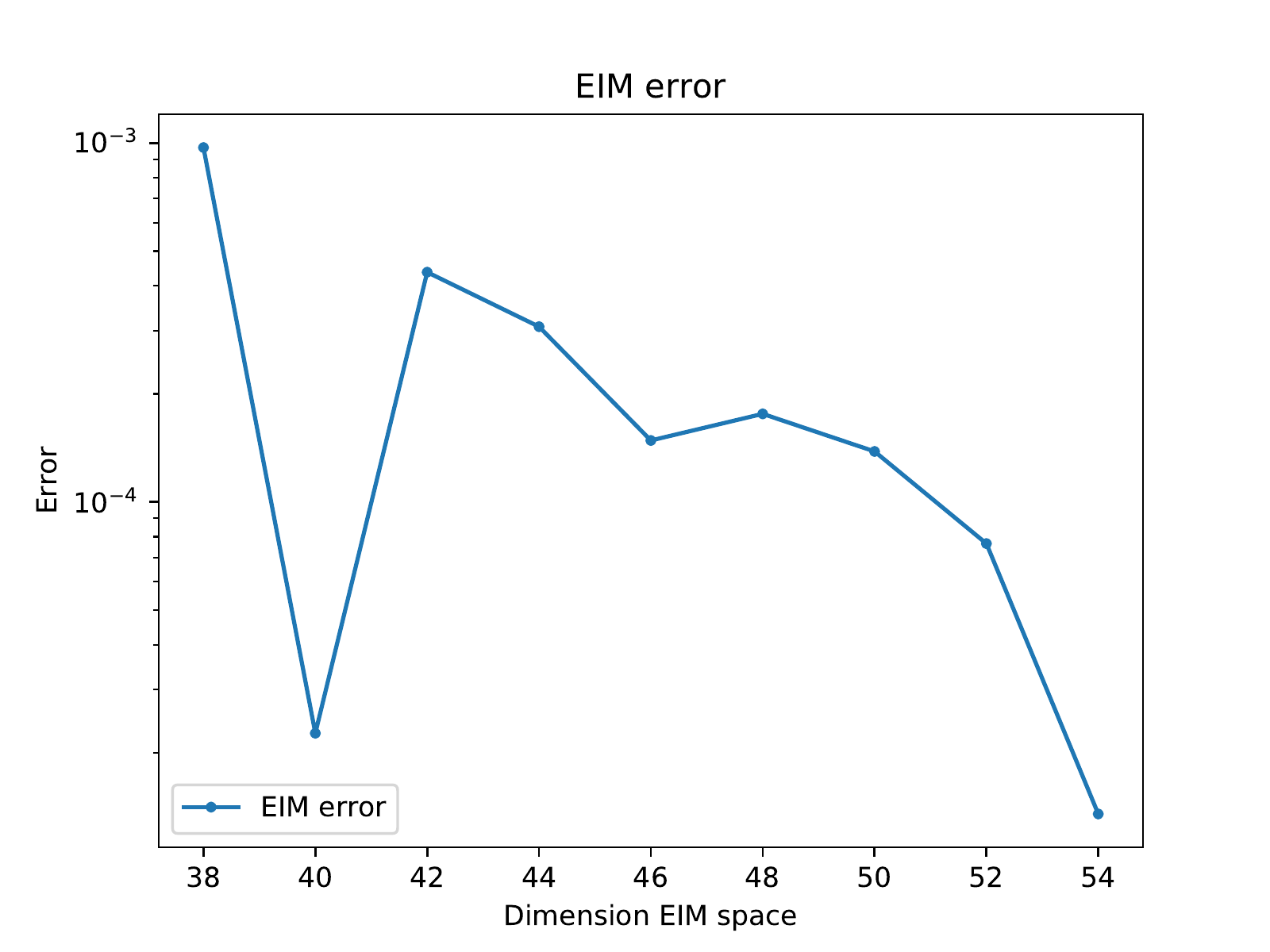}}
\end{subfigure}
\caption{Offline phases of advection of solitary wave \label{fig:test0:offline}}
\end{center}
\end{figure}
This problem travels at constant speed $\mu_0$, i.\,e. $\theta \approx t\mu_0+\mu_2$, hence it is straightforward to obtain an almost perfect map with polynomial interpolation of degree at least 2. In \cref{fig:test0:training} we see how the polynomial of degree 2 obtains a very good regression map with just 10 samples. This is true also for the piecewise linear interpolation. The higher order polynomials need more samples to regularize their coefficients and the neural network has an even slower decay of the error. 

In \cref{fig:test0:offlineEu,fig:test0:offlineALE,fig:test0:EIMEu} we compare different \textit{offline} phases for the ALE with second order polynomial regression map and the Eulerian algorithm with a tolerance of $10^{-3}$.

The Eulerian algorithm produced RB and EIM spaces of dimensions (52, 54), while the second order polynomial ALE has dimensions (4, 7) and the ANN ALE has (12, 20) as one can see in \cref{tab:RBdata}. We already see the big improvement, in particular for the case where the learned regression map behaves well. 

The RB space of the Eulerian framework does not behave particularly bad in this situation because the scheme we used is very diffusive and the simulated shape functions are smooth enough.
In computational time we can appreciate the advantage of the ALE algorithm only relatively. Indeed, in \cref{tab:RBdata} we see that even if the time for computing a ALE solution is roughly the double of the Eulerian solution, we still gain something with this approach, because we obtain very small reduced basis spaces.

\begin{figure}[h!]
	\begin{center}
\begin{subfigure}[ALE FOM and RB solutions for Poly2\label{fig:test0:onlineALE}]{\includegraphics[width=0.45\textwidth, trim={0 0 20 40},clip]{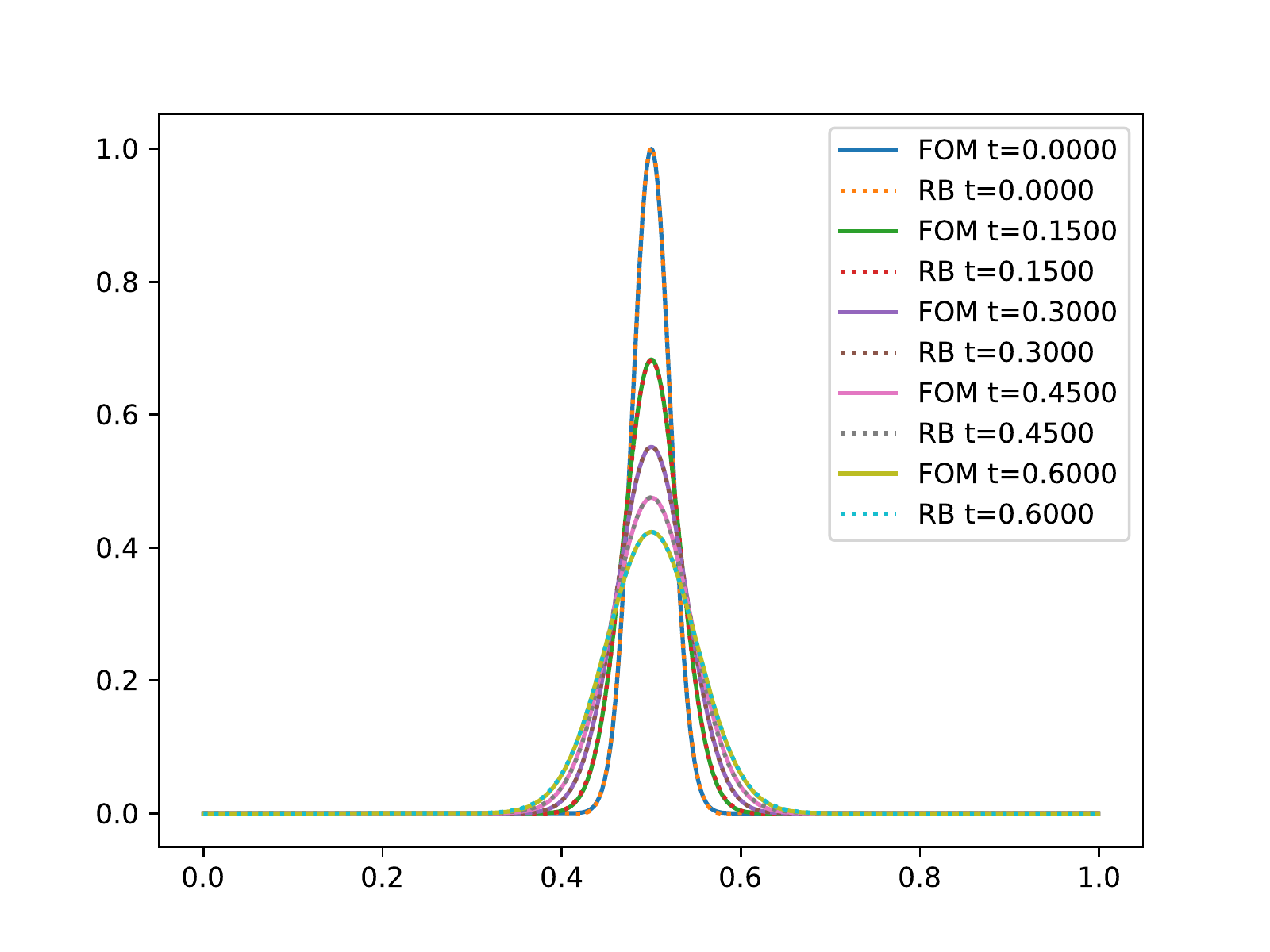}}
\end{subfigure}
\begin{subfigure}[ALE FOM and RB solutions for ANN\label{fig:test0:onlineALENN}]{\includegraphics[width=0.45\textwidth, trim={0 0 20 40},clip]{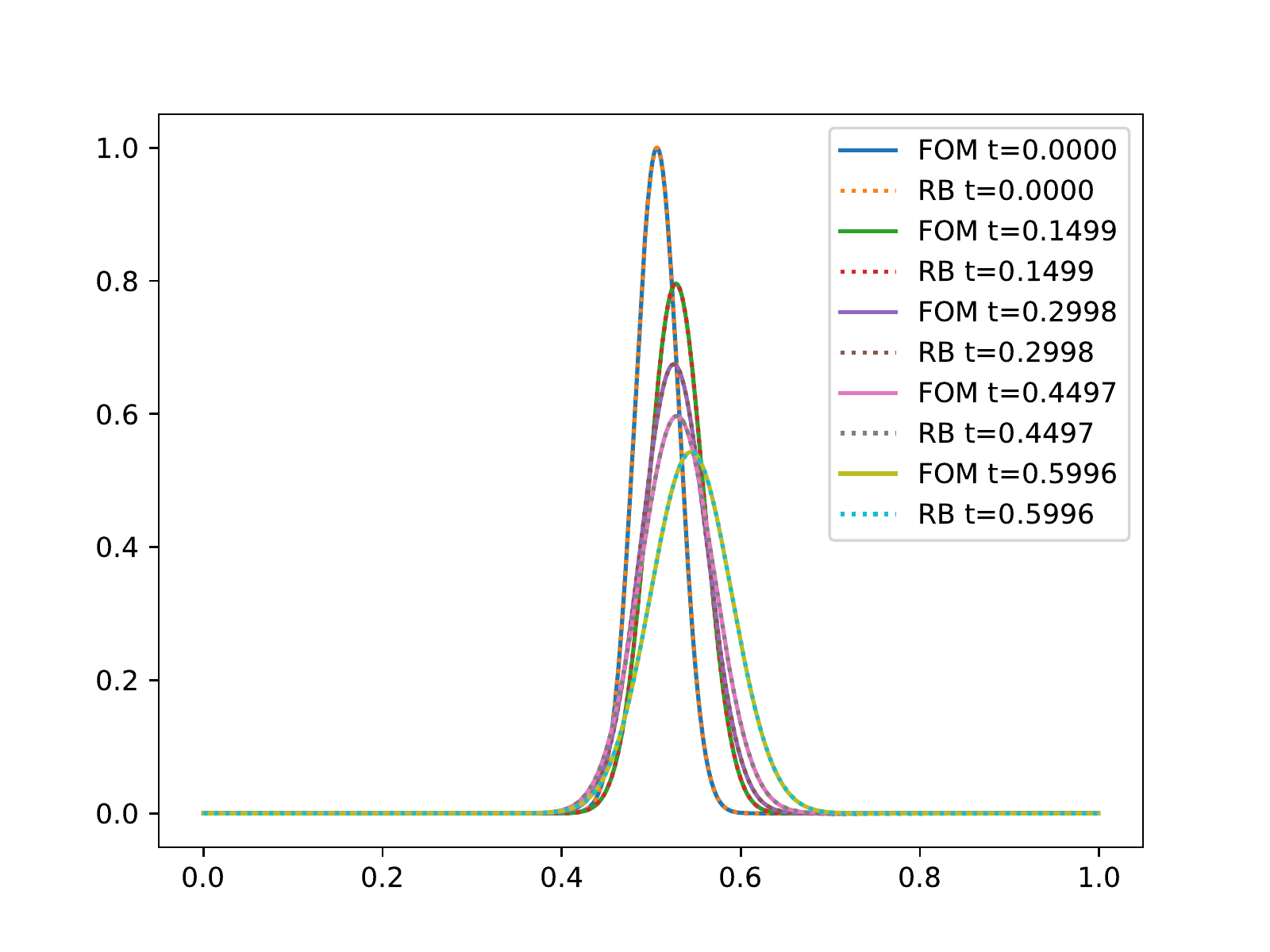}}
\end{subfigure}\caption{Online phase for ALE algorithms\label{fig:test0:ALEonline}}
\end{center}
\end{figure}
In \cref{fig:test0:onlineALE,fig:test0:onlineALENN} we see the quality of the solutions for one parameter of the domain with the two regressions. We see the troubles that the RB algorithm with the ANN has to tackle due to the not precise alignment of the selected feature, but the quality of the solutions still remains accurate.

\subsection{Advection of Shock Wave} 
In this test, we consider a shock traveling at a parametric speed with a random initial position. The domain is $\Omega=[0,1]$ and the initial conditions are 
\begin{equation}\label{eq:advShock}
u_0(x,\bmu)=\begin{cases}
\mu_1\quad &\text{if } x<0.35+0.05\mu_2, \\
0 \quad &\text{else},
\end{cases} \qquad \mu_1,\mu_2 \in [-1,1].
\end{equation}
We solve again the linear transport equation \eqref{eq:linAdv} with $a=\mu_0 \in [0,2]$ till final time $T=1.5$. The boundary conditions are inflow on the left and outflow on the right. We use the hyperbolic dilatation \eqref{eq:dilatation} as calibration transformation. We select the steepest point of the solution as calibration point.

\begin{figure}[h!]
	\begin{center}
	\begin{subfigure}[Offline error in ALE framework\label{fig:test91:offlineALE}]{\includegraphics[width=0.45\textwidth, trim={0 0 20 40},clip]{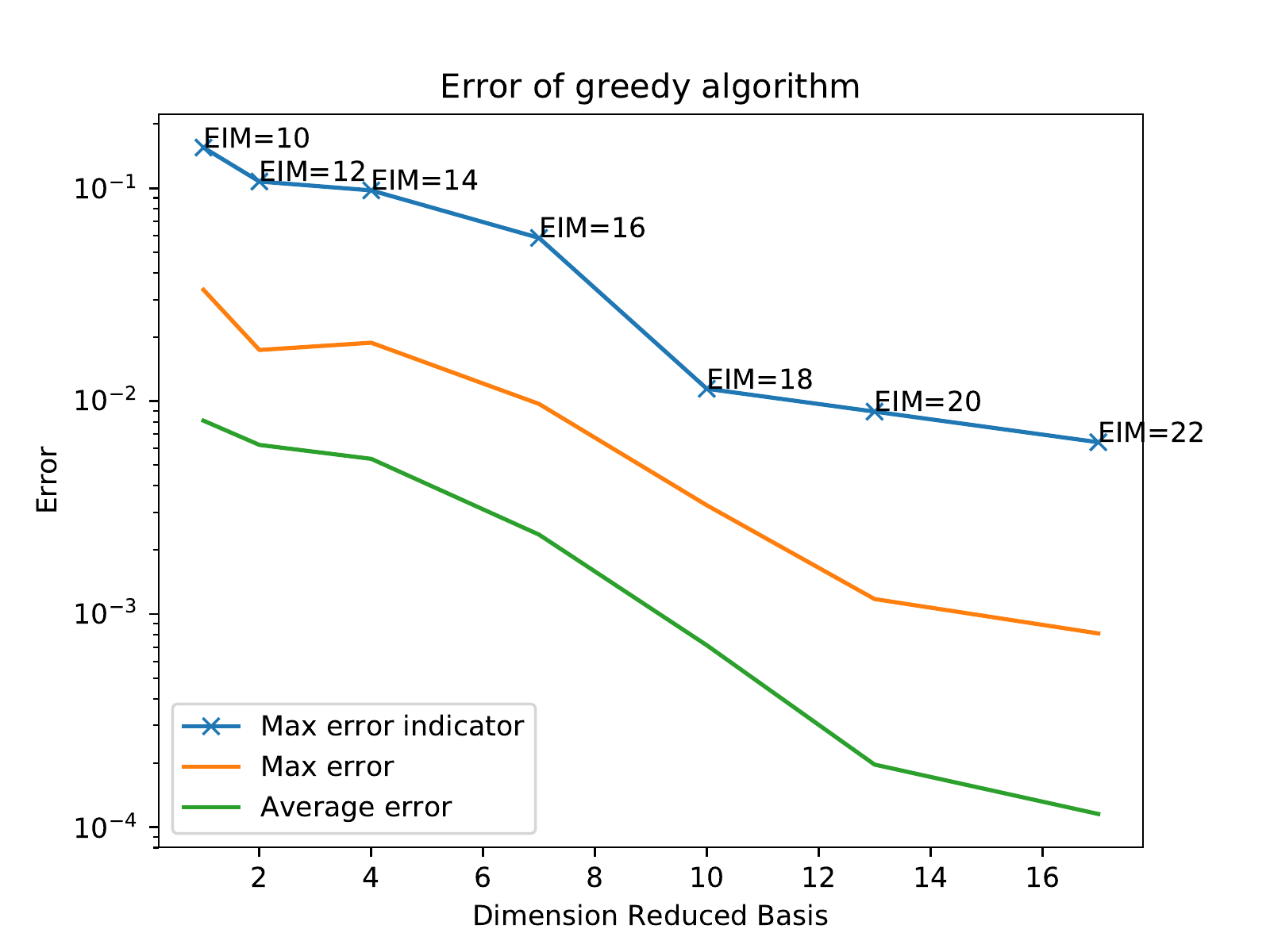}}
\end{subfigure}
	\begin{subfigure}[Training of calibration\label{fig:test91:training}]{\includegraphics[width=0.45\textwidth, trim={0 0 20 40},clip]{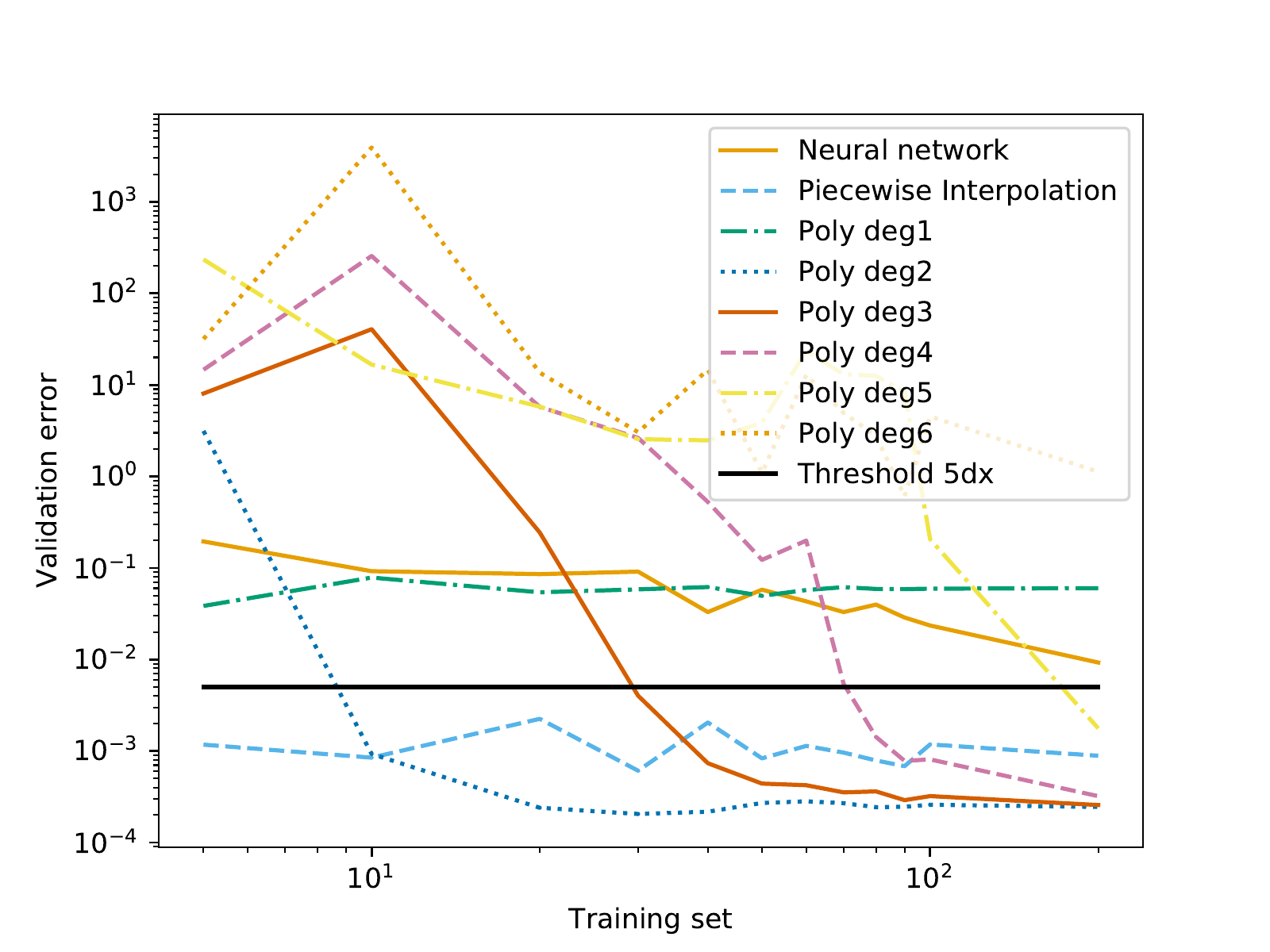}}
\end{subfigure}
\\
	\begin{subfigure}[Offline error in Eulerian framework\label{fig:test91:offlineEu}]{\includegraphics[width=0.45\textwidth, trim={0 0 20 40},clip]{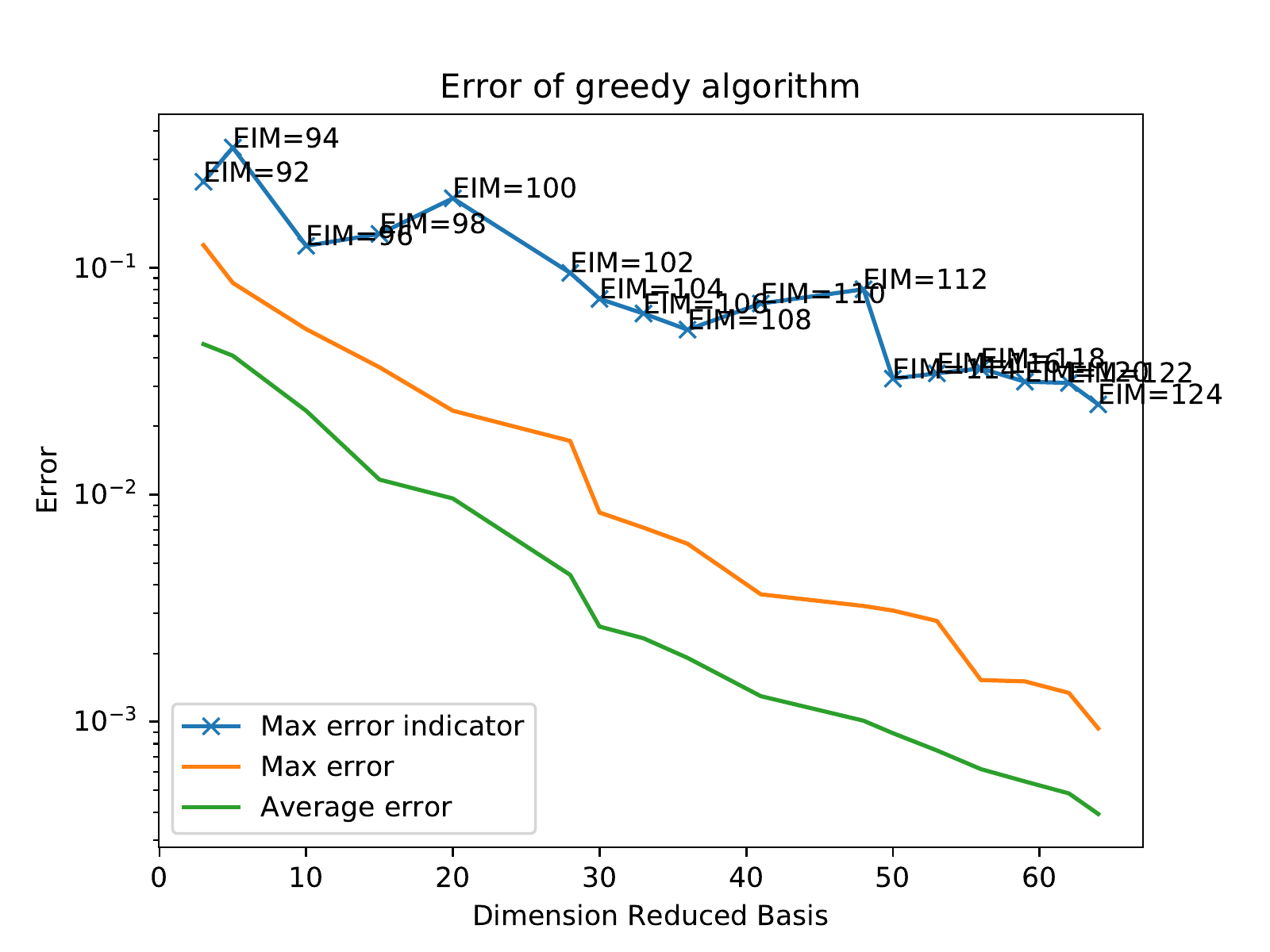}}
\end{subfigure}
\begin{subfigure}[EIM offline error in Eulerian framework\label{fig:test91:EIMEu}]{\includegraphics[width=0.45\textwidth, trim={0 0 20 37},clip]{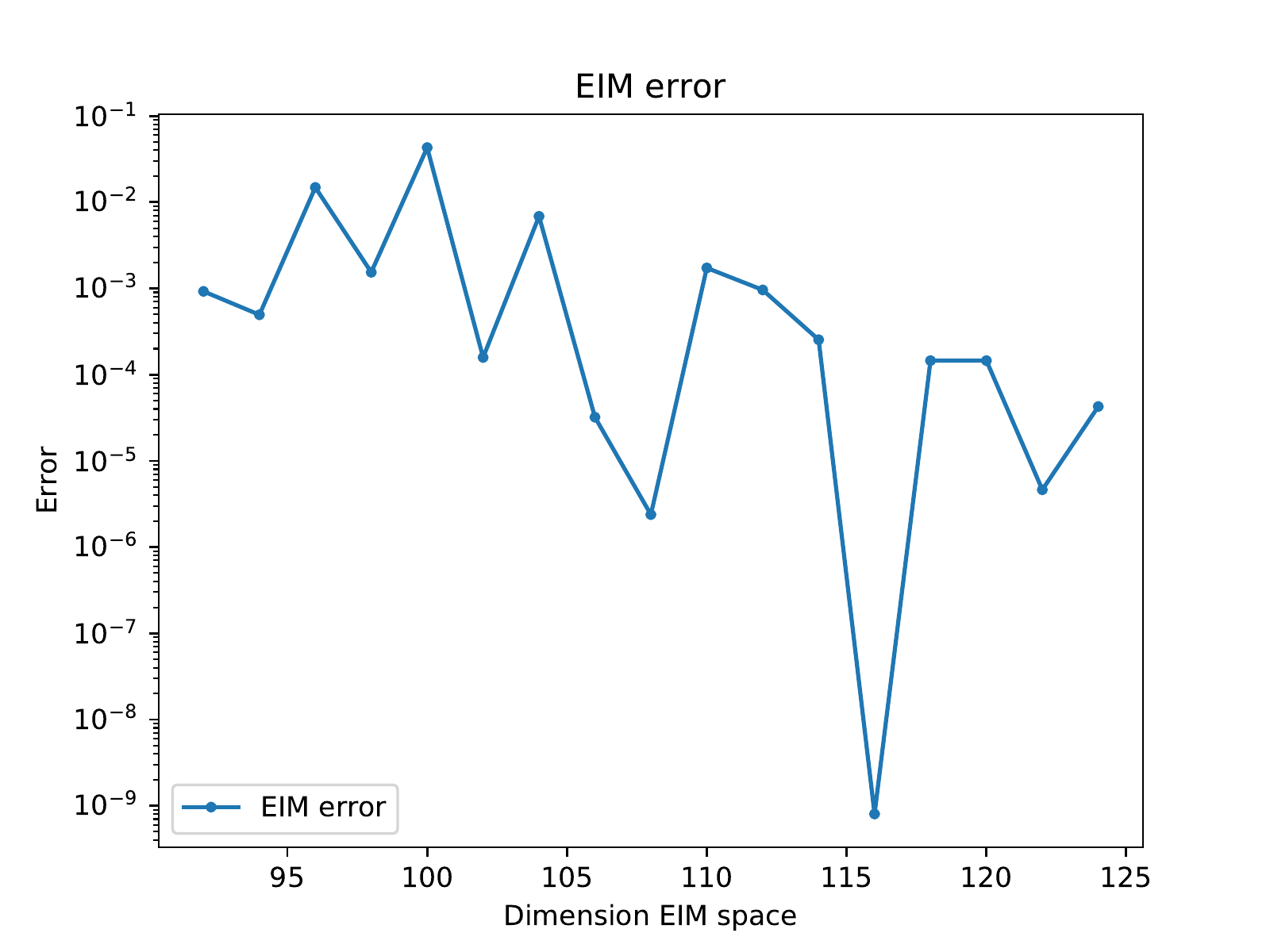}}
\end{subfigure}\\
	\begin{subfigure}[FOM and RB solutions for ALE sovlers \label{fig:test91:onlineALE}]{\includegraphics[width=0.45\textwidth, trim={0 0 20 40},clip]{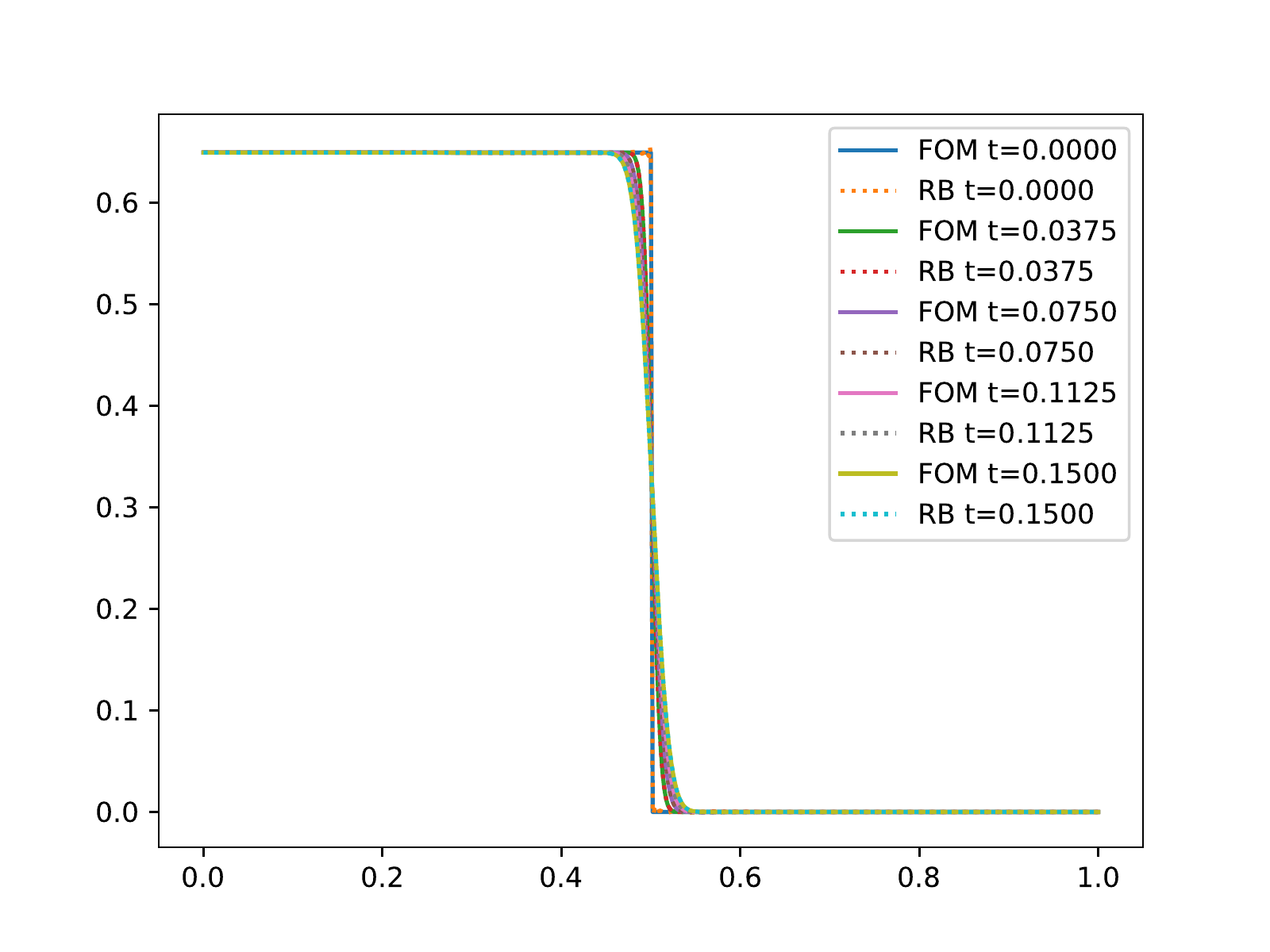}}
\end{subfigure} 
\begin{subfigure}[FOM and online RB solutions for Euler solvers \label{fig:test91:onlineEulRB}]{\includegraphics[width=0.45\textwidth, trim={0 0 20 40},clip]{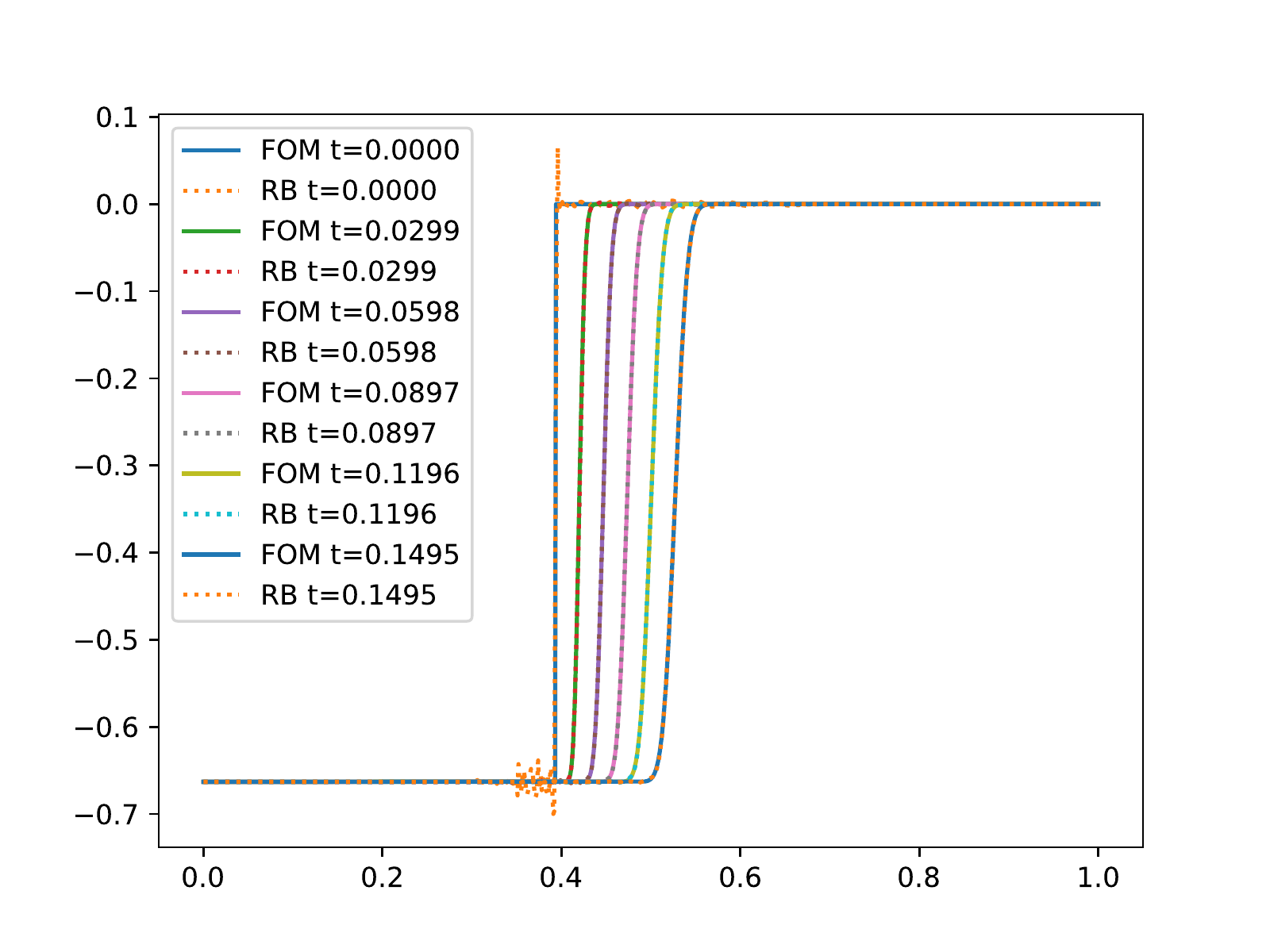}}
\end{subfigure}
\caption{Shock wave tests \label{fig:test91:offline}}
\end{center}
\end{figure}
As for the previous test, this problem travels at constant speed, hence polynomials of degree 2 are the sought regression map. In \cref{fig:test91:training} we see which regression maps are quickly converging. Similar considerations to the previous test hold for this one. 

In \cref{fig:test91:offlineALE,fig:test91:offlineEu,fig:test91:EIMEu} we compare, for a tolerance of $10^{-3}$, the classical PODEI--Greedy errors, the classical EIM errors and the new ALE PODEI Greedy with the polynomial regression of degree 2. 
More than before, we see how it is easy to catch the behavior of the solutions with few basis functions in RB and EIM spaces in ALE framework (17, 22), while the algorithm struggles in representing the right solutions in the Eulerian framework where many basis functions are needed (64, 124). 
In \cref{tab:RBdata} we compare again the dimensions and the computational times, where we observe a strong advantage in the new methodology.

In \cref{fig:test91:onlineALE,fig:test91:onlineEulRB} we can observe the main reason why the new methodology defeats the classical one. In \cref{fig:test91:onlineALE} the ALE--PODEI--Greedy produces high quality reduced solution, while in \cref{fig:test91:onlineEulRB} the Eulerian framework obtains strong spurious oscillations that are very dangerous in many physical application, where they can represent, for example, negative density or pressure.

\subsection{Burgers Oscillation} 
In this test, we solve the Burgers' equation \eqref{eq:burgers} on the domain $\Omega=[0,1]$, with as initial conditions an oscillation dumped at the boundaries. This problem can develop in finite time a shock, accordingly to parameters, and it can travel in both directions with nonlinear speed. 
The initial conditions are, more precisely,
\begin{equation}\label{eq:BurOsc}
u_0(x,\bmu)= \sin(2\pi (x+0.1\mu_1))e^{-(60+20\mu_2)(x-0.5)^2}(1+0.5\mu_3x),\qquad \mu_1 \in [0,1],\, \mu_2,\mu_3 \in [-1,1],
\end{equation}
with $a=\mu_0\in [0,2]$ till final time $T=0.6$. We use homogeneous Dirichlet  boundary conditions and the dilatation \eqref{eq:dilatation} as calibration map. The detection criterion for the calibration point is chosen checking the point of the function that crosses the x--axis.
\begin{figure}[h!]
	\begin{center}
		\begin{subfigure}[Training of regression maps for  calibration\label{fig:test6:training}]{\includegraphics[width=0.45\textwidth, trim={0 0 20 40},clip]{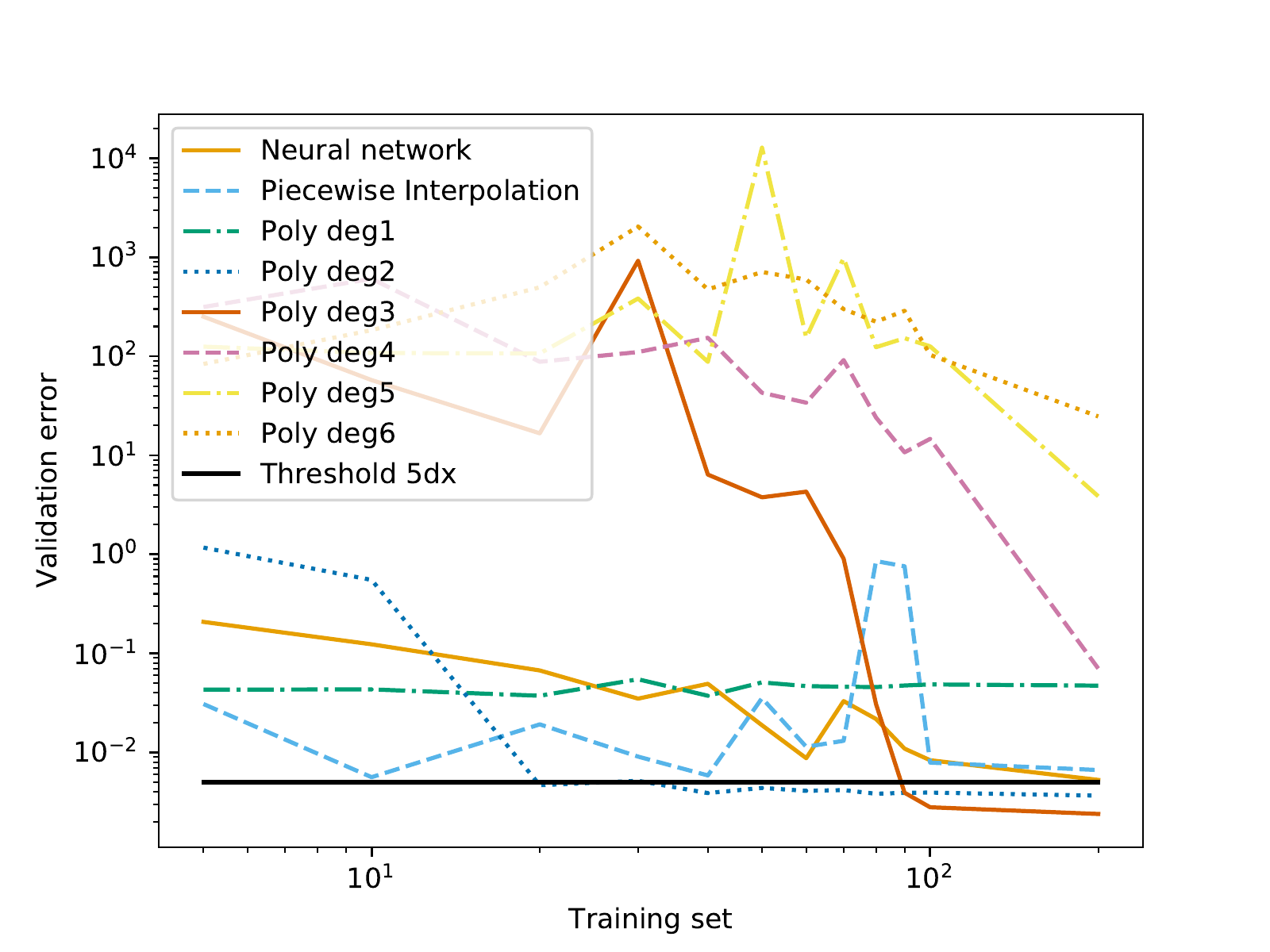}}
\end{subfigure}
\begin{subfigure}[Offline error of Greedy in Eulerian framework\label{fig:test6:offlineEu}]{\includegraphics[width=0.45\textwidth, trim={0 0 10 40},clip]{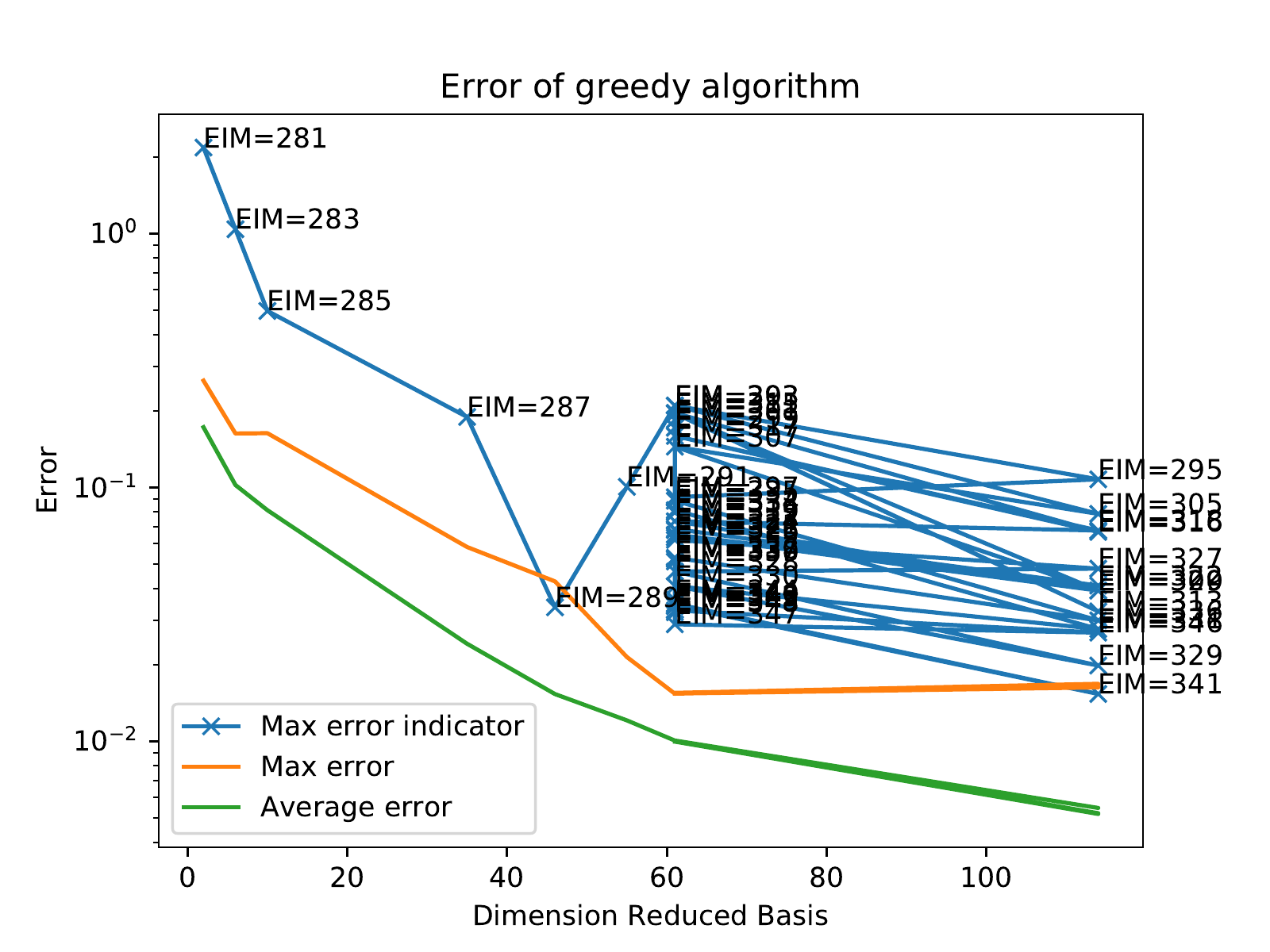}}
\end{subfigure}\\
	\begin{subfigure}[Original solutions \label{fig:test6:onlineFOM}]{\includegraphics[width=0.45\textwidth, trim={0 0 20 40},clip]{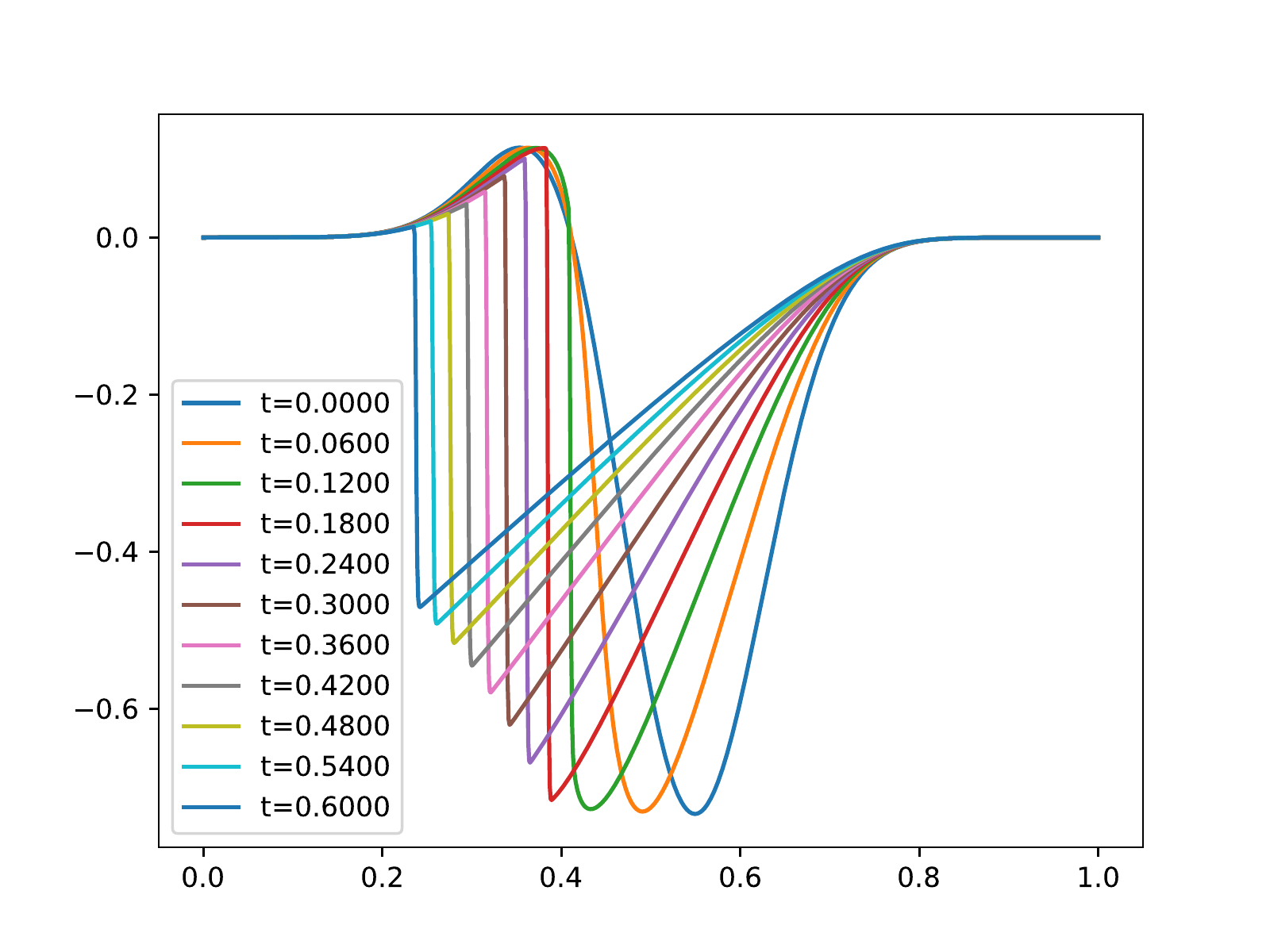}}
	\end{subfigure}
	\begin{subfigure}[ALE FOM and RB solutions\label{fig:test6:onlineALE}]{\includegraphics[width=0.45\textwidth, trim={0 0 20 40},clip]{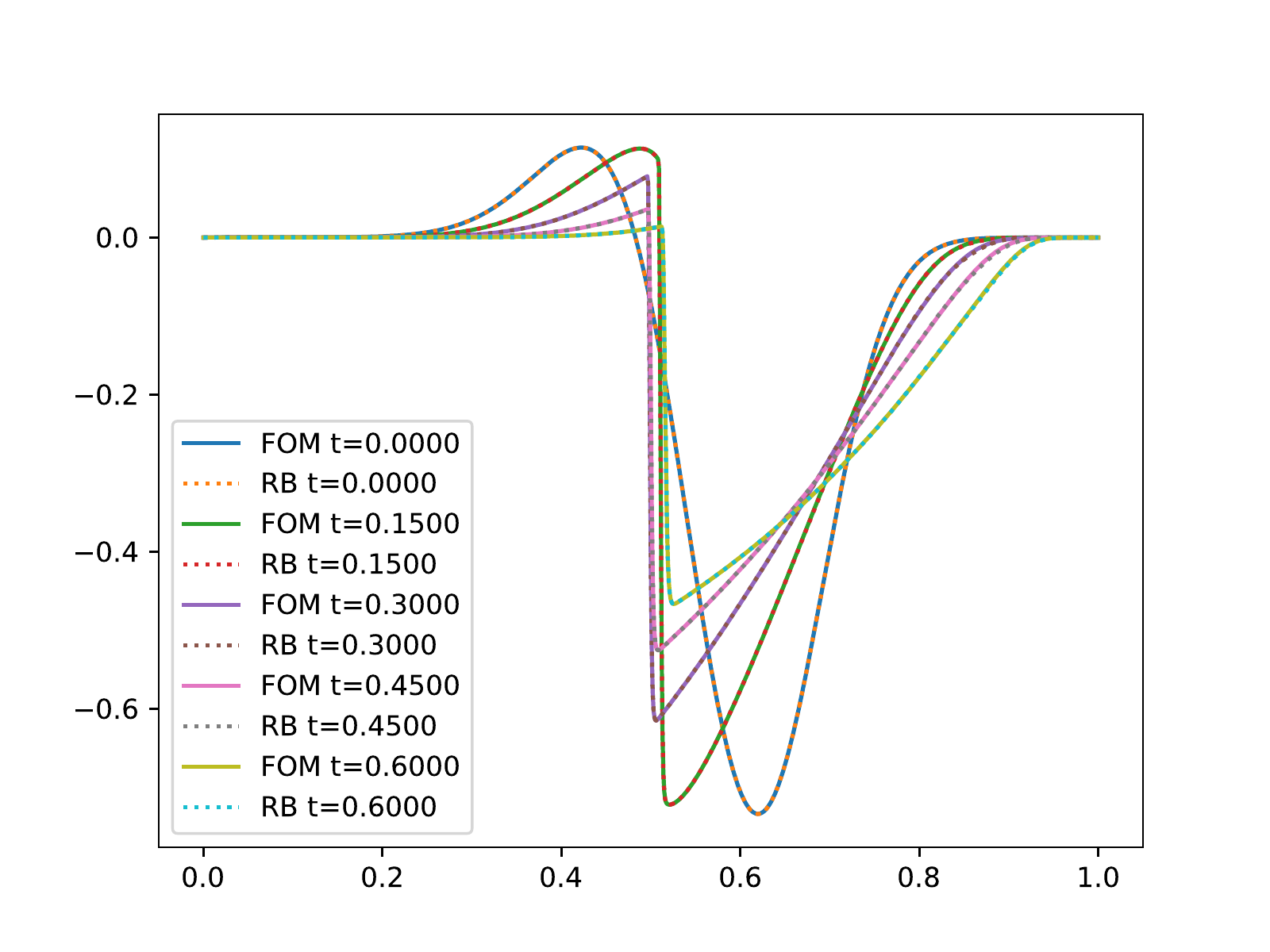}}
	\end{subfigure} \\
	\begin{subfigure}[Eulerian FOM and RB solutions\label{fig:test6:onlineEulRB}]{\includegraphics[width=0.45\textwidth, trim={0 0 20 40},clip]{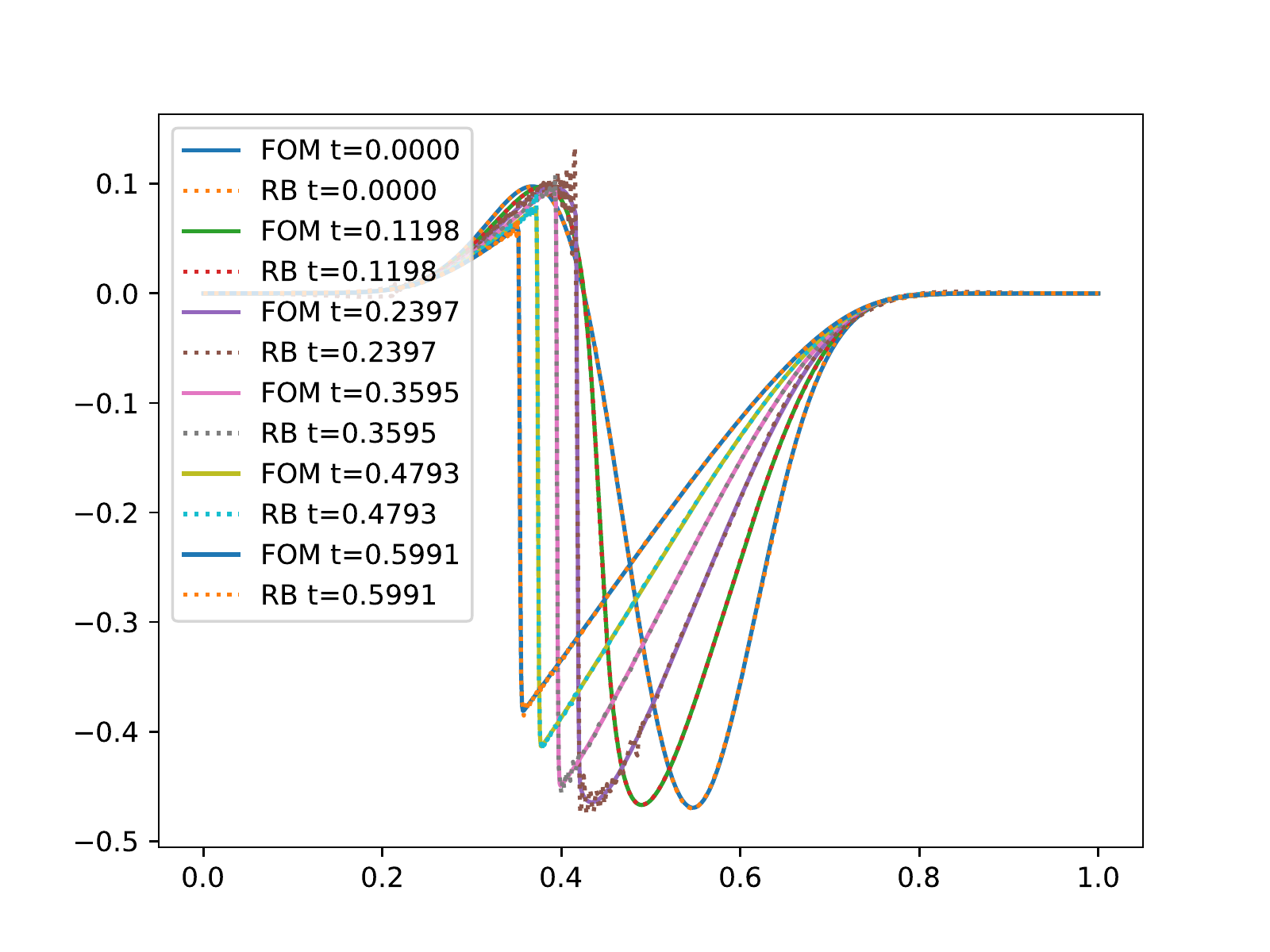}}
	\end{subfigure}
	\begin{subfigure}[Eulerian FOM and RB solutions (zoom)\label{fig:test6:onlineEulRBzoom}]{\includegraphics[width=0.45\textwidth, trim={0 0 20 40},clip]{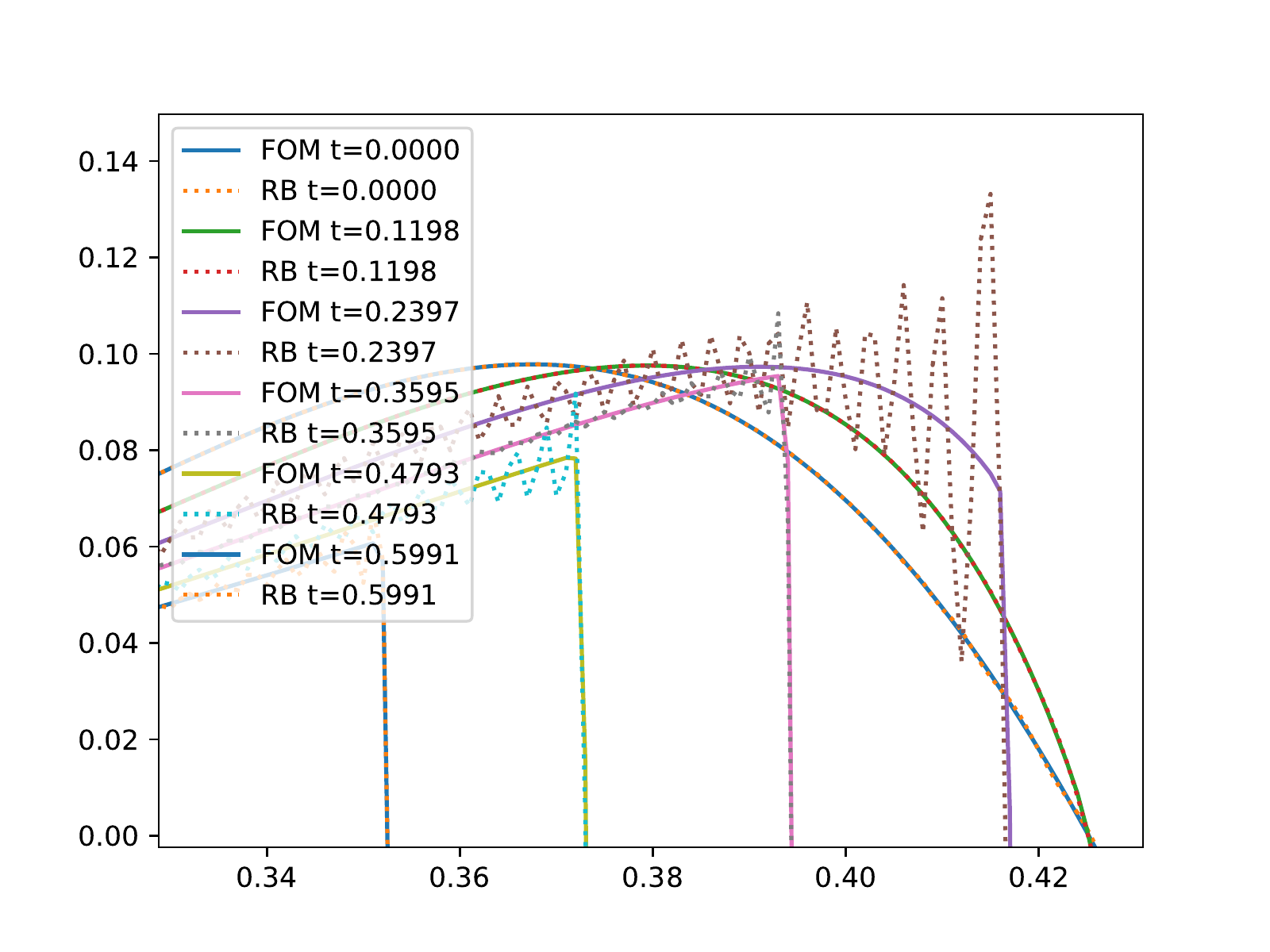}}
	\end{subfigure}\caption{Burgers oscillation test\label{fig:test6:ALEonline}}
\end{center}
\end{figure}

In this test, the calibration parameters are much harder to be found. 
In \cref{fig:test6:training} we see that all the regression maps just barely touch the threshold line. This is also due to the larger number of parameters of the problem. Nevertheless, we test the third order polynomials, which seem able to capture the behavior of the calibration curve. For our simulations we choose also the ANN to compare the results.


In \cref{fig:test6:onlineFOM,fig:test6:onlineALE,fig:test6:onlineEulRB,fig:test6:onlineEulRBzoom}  we plot the offline phase for the classical Eulerian algorithm and we see that it really needs many EIM basis functions before the error can decrease and we do not even reach the tolerance of $10^{-3}$. In \cref{tab:RBdata} are stored the other values for the \textit{offline } phases of the Poly3 ALE--PODEI--Greedy and the ANN ALE--PODEI--Greedy. 
For this simulation, in all algorithms the convergence is slower, since the problem is more involved. 
In the ALE framework, nevertheless, with 50 RB functions and 60 EIM bases we obtain the tolerance sought, while with the Eulerian algorithm the EIM space must span more than one third of the FOM space, leading to very high dimensional EIM and RB spaces and long computational times.
In \cref{tab:RBdata} we observe a strong advantage in the computational time of the new methodology. 

In \cref{fig:test6:onlineALE} we see the quality of the solutions in the ALE framework, while in \cref{fig:test6:onlineEulRBzoom} in the Eulerian framework we notice again the even higher spurious oscillations.

\subsection{Burgers Sine} 
In this test, we solve the Burgers' equation \eqref{eq:burgers} on the domain $\Omega=[0,\pi]$, with the absolute value of the sine as initial conditions. This problem develop in finite time a shock and it travels with a speed that depends nonlinearly on the parameters of the problem. The initial conditions are, more precisely,
\begin{equation}\label{eq:BurSin}
u_0(x,\bmu)=\lvert \sin(x+\mu_1) \rvert + 0.1,\qquad \mu_1 \in [0,\pi],
\end{equation}
with the coefficient of the Burgers' equation $a=\mu_0 \in [0,2]$ till final time $T=0.15$. We use periodic boundary conditions and the translation \eqref{eq:translation} as calibration map. We select the calibration point as the minimum value of the solution.

\begin{figure}[h]
	\begin{center}
	\begin{subfigure}[Training of regression maps for  calibration\label{fig:test1:training}]{\includegraphics[width=0.45\textwidth, trim={0 0 20 40},clip]{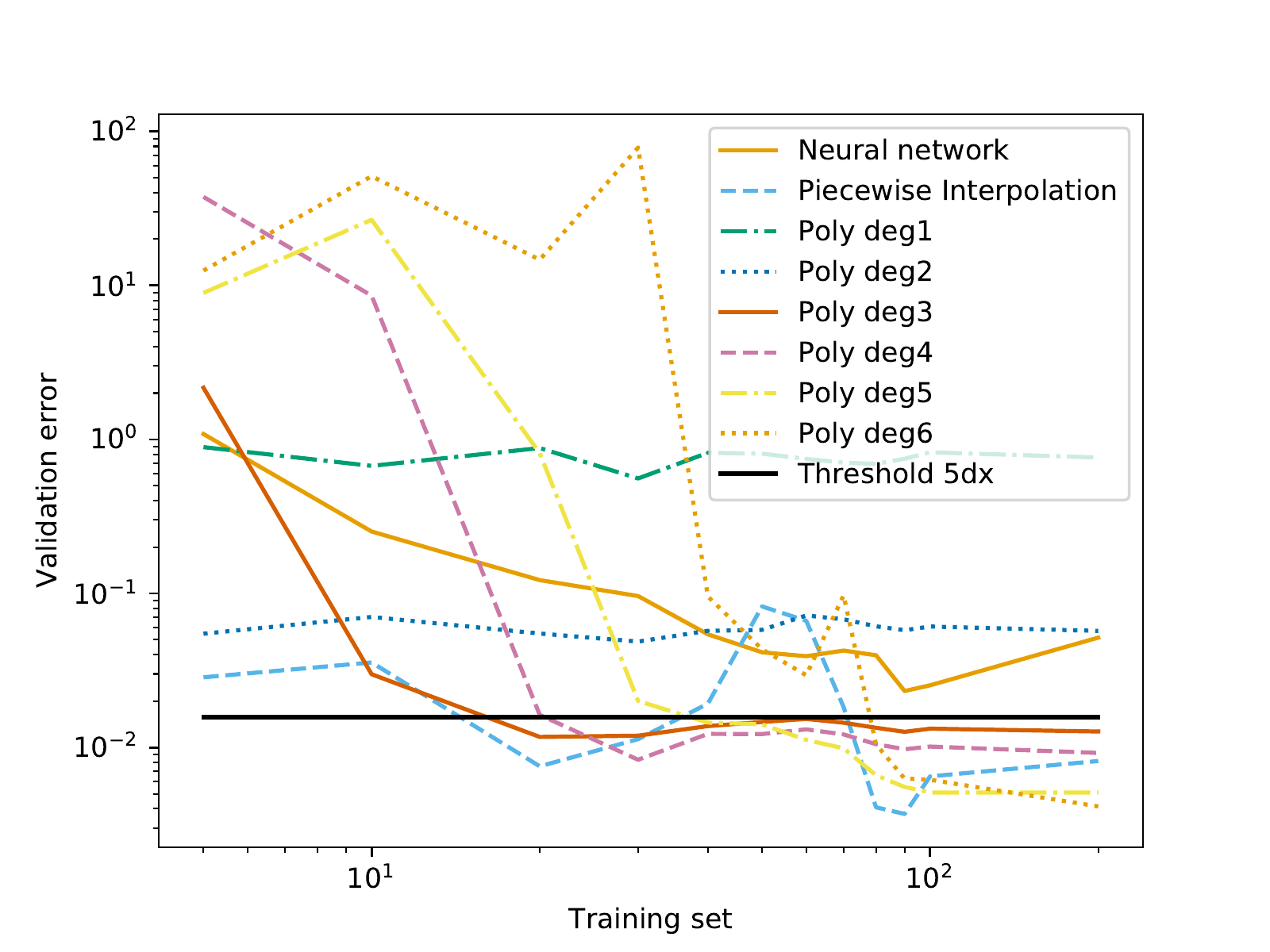}}
	\end{subfigure}
	\begin{subfigure}[Offline error of Greedy in ALE framework with Poly4\label{fig:test1:offlineEu}]{\includegraphics[width=0.45\textwidth, trim={0 0 20 40},clip]{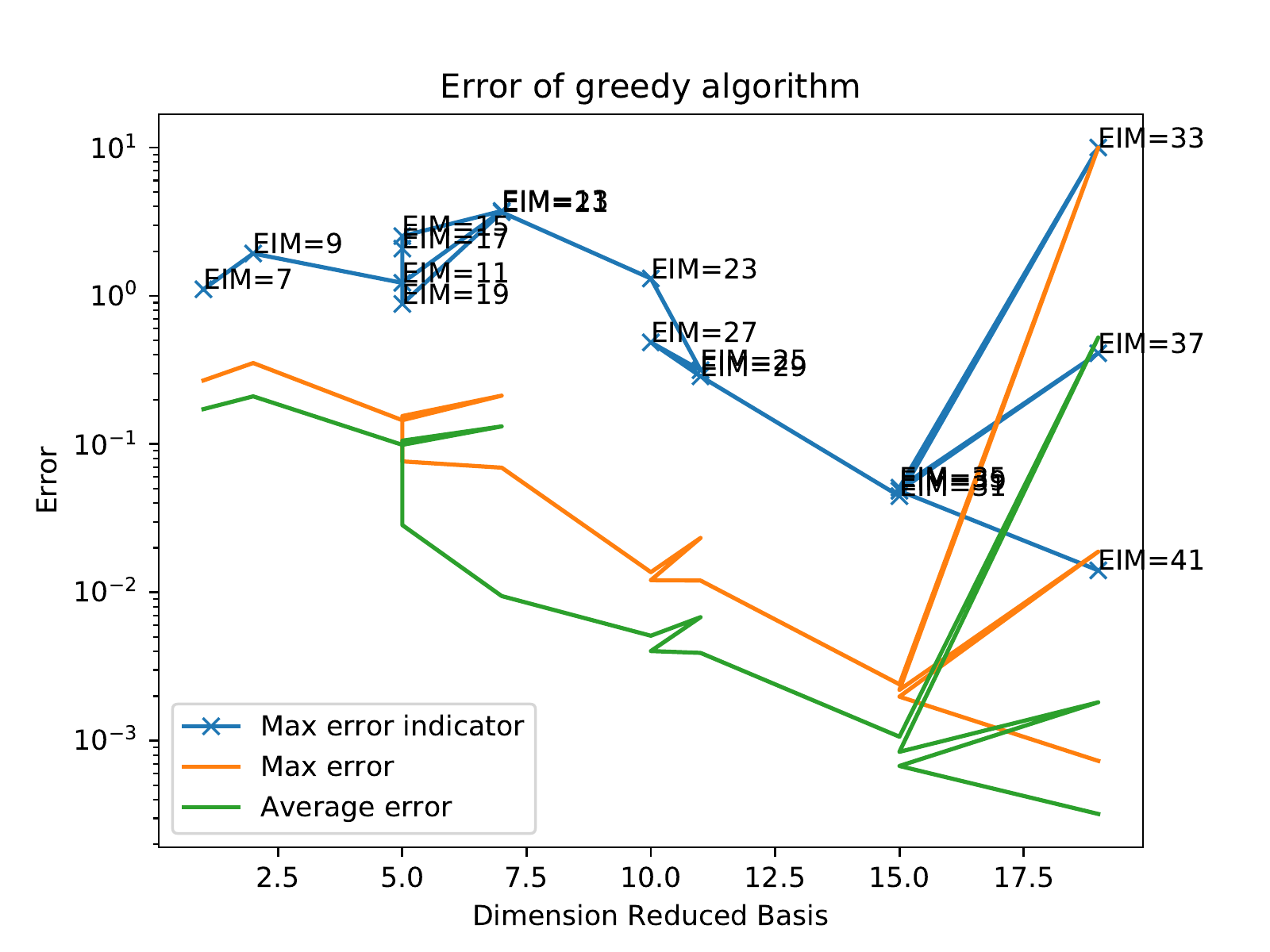}}
	\end{subfigure}\\
	\begin{subfigure}[Original solutions \label{fig:test1:onlineFOM}]{\includegraphics[width=0.45\textwidth, trim={0 0 20 40},clip]{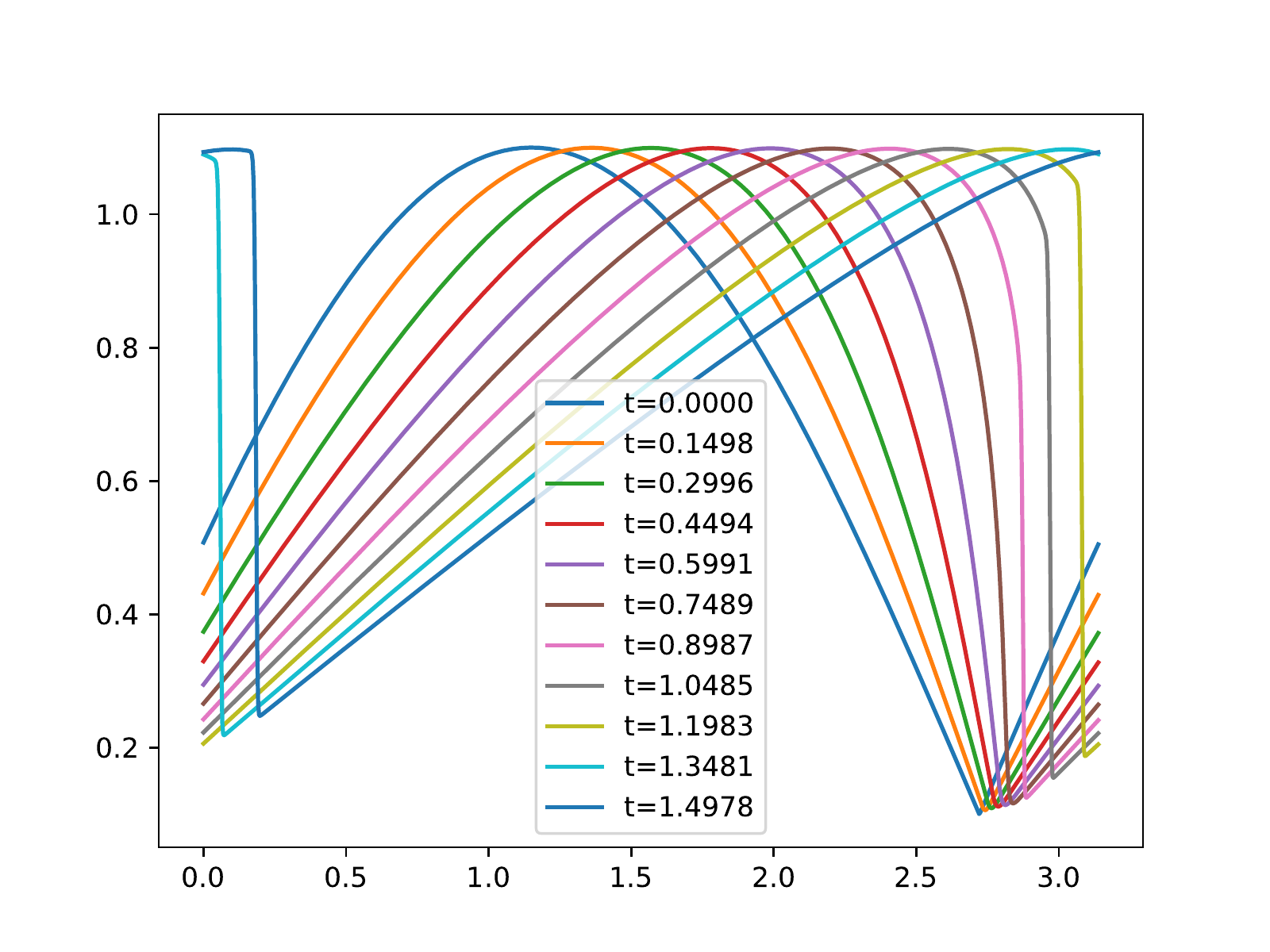}}
	\end{subfigure}
	\begin{subfigure}[ALE FOM and RB solutions\label{fig:test1:onlineALE}]{\includegraphics[width=0.45\textwidth, trim={0 0 20 40},clip]{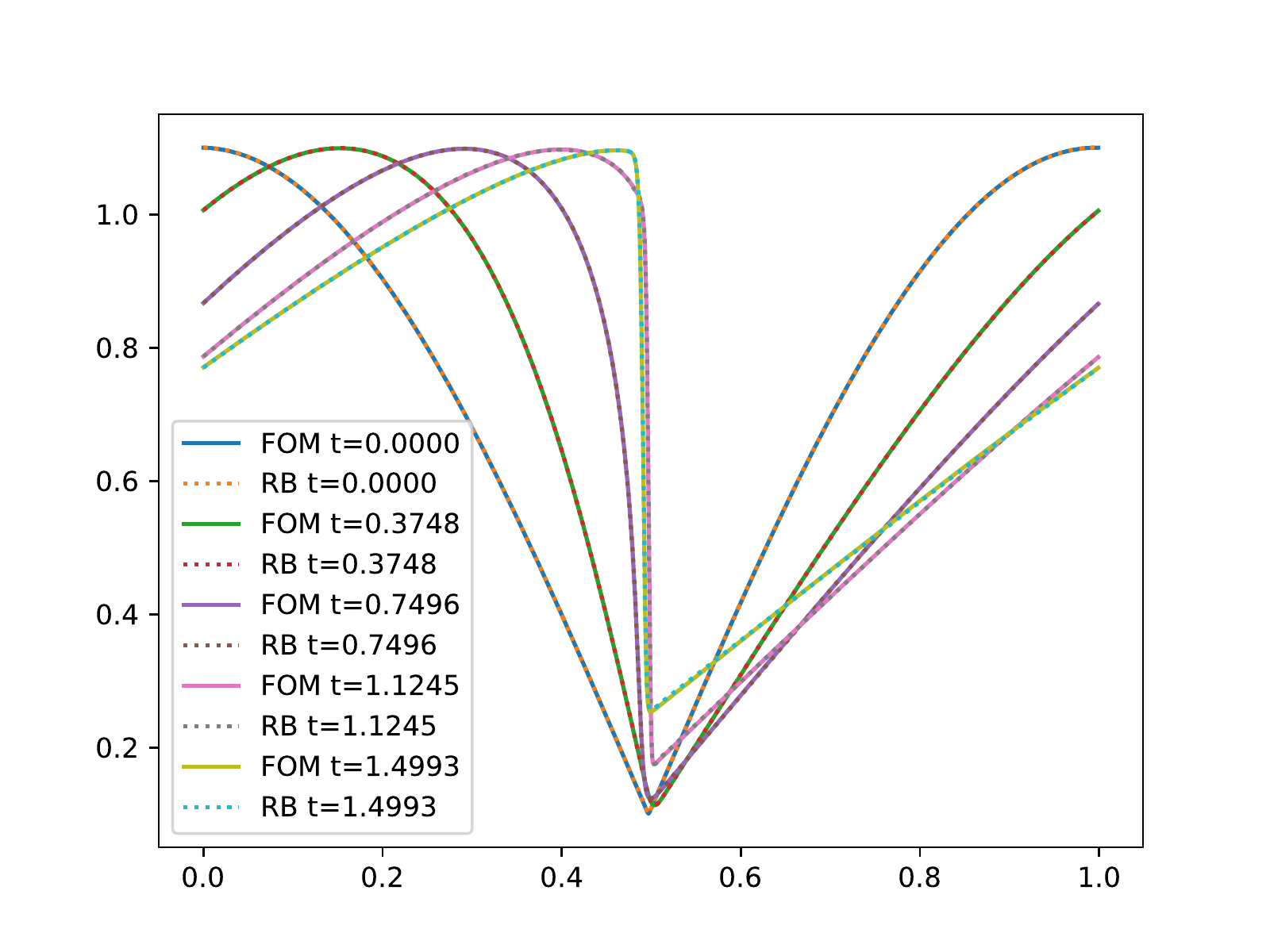}}
\end{subfigure}
\end{center}
	\caption{Burgers sine test\label{fig:test1:regression}}
\end{figure}

\begin{figure}[h!]
	\begin{center}
	\begin{subfigure}[Training of regression maps for  calibration\label{fig:test10:training}]{\includegraphics[width=0.45\textwidth, trim={0 0 20 40},clip]{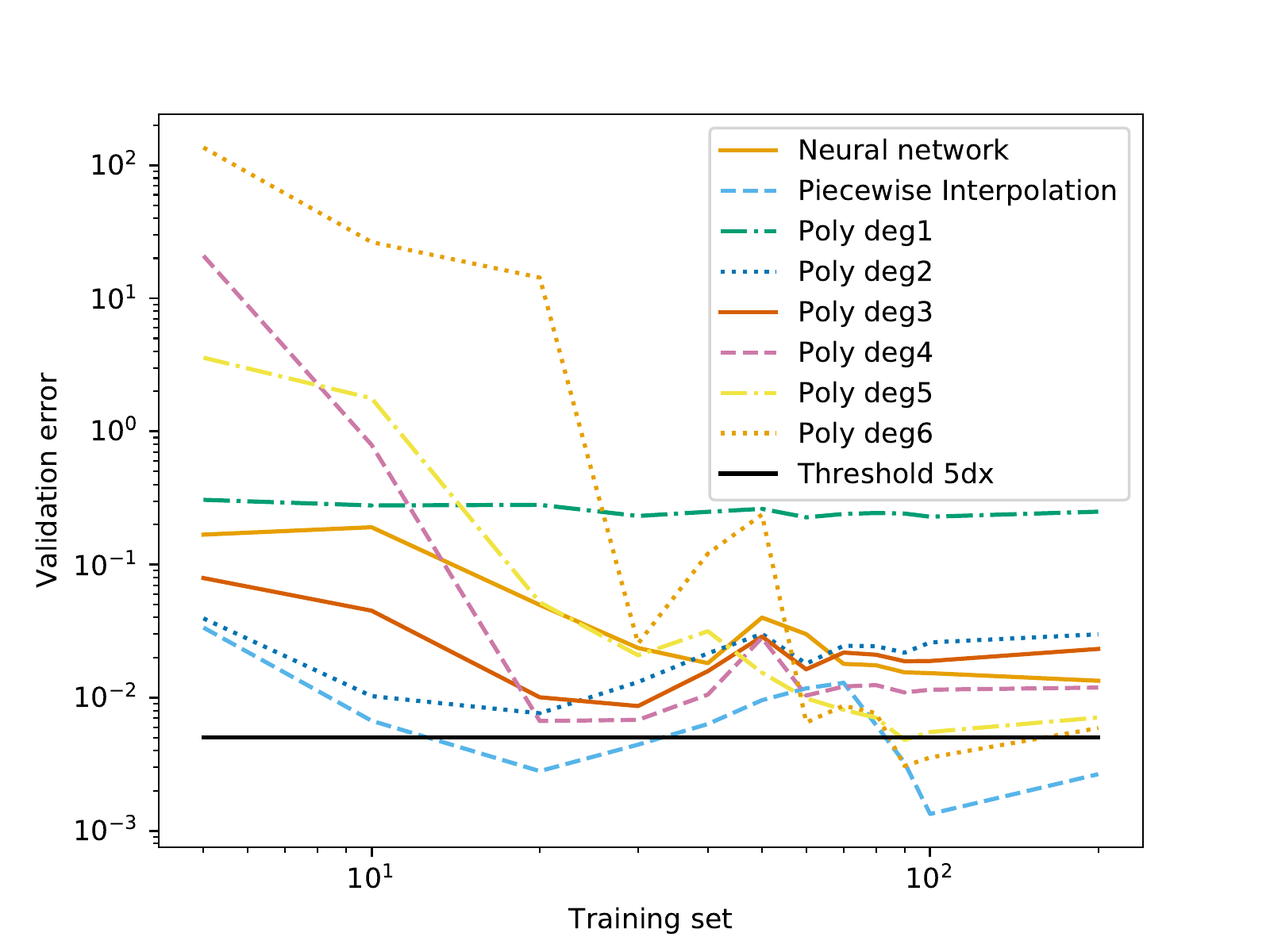}}
	\end{subfigure}
	\begin{subfigure}[Original solutions \label{fig:test10:onlineFOM}]{\includegraphics[width=0.45\textwidth, trim={0 0 20 40},clip]{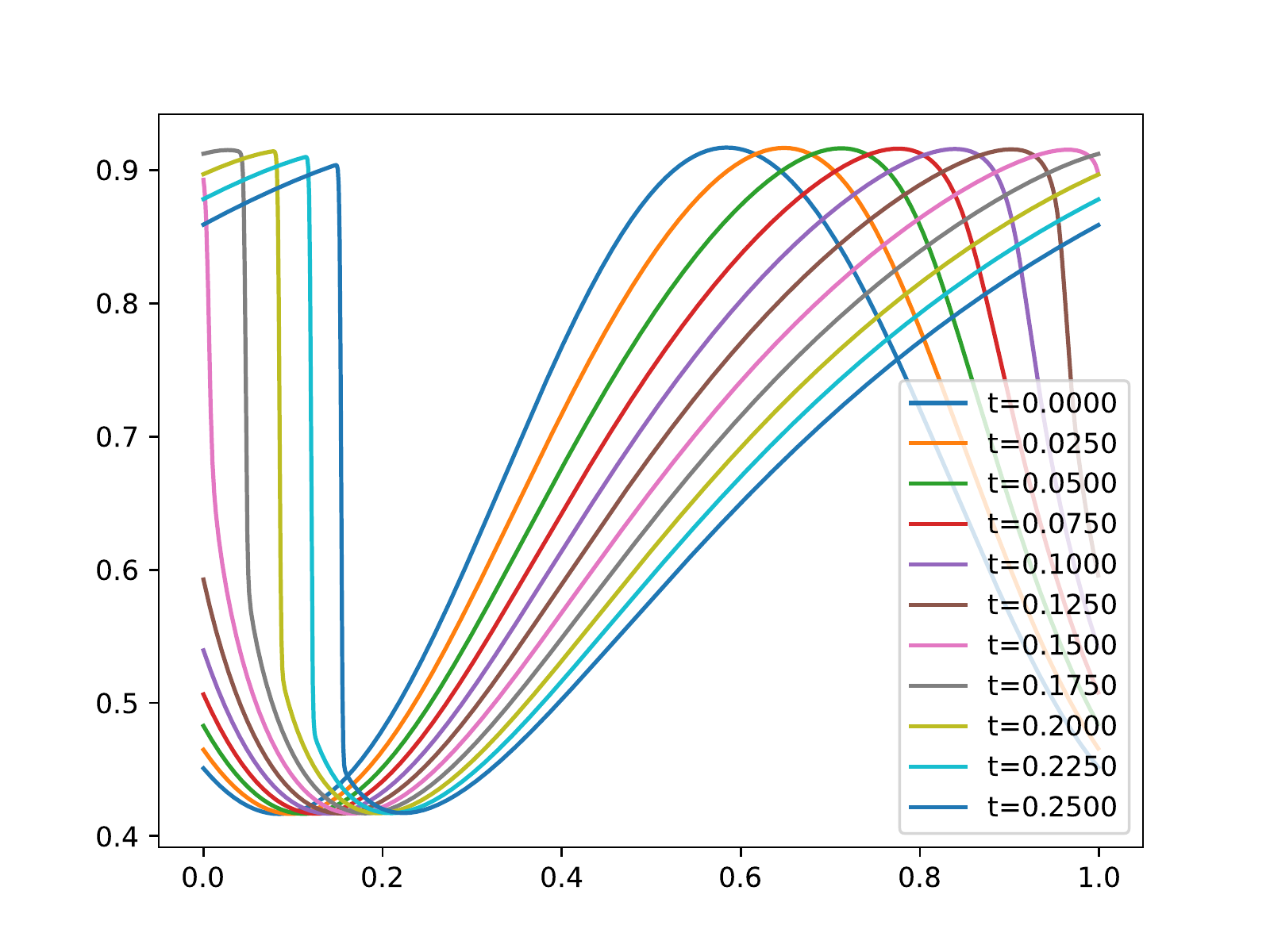}}
		\end{subfigure}\\
	\begin{subfigure}[Offline error of Greedy in Eulerian framework\label{fig:test10:offlineEu}]{\includegraphics[width=0.45\textwidth, trim={0 0 19 40},clip]{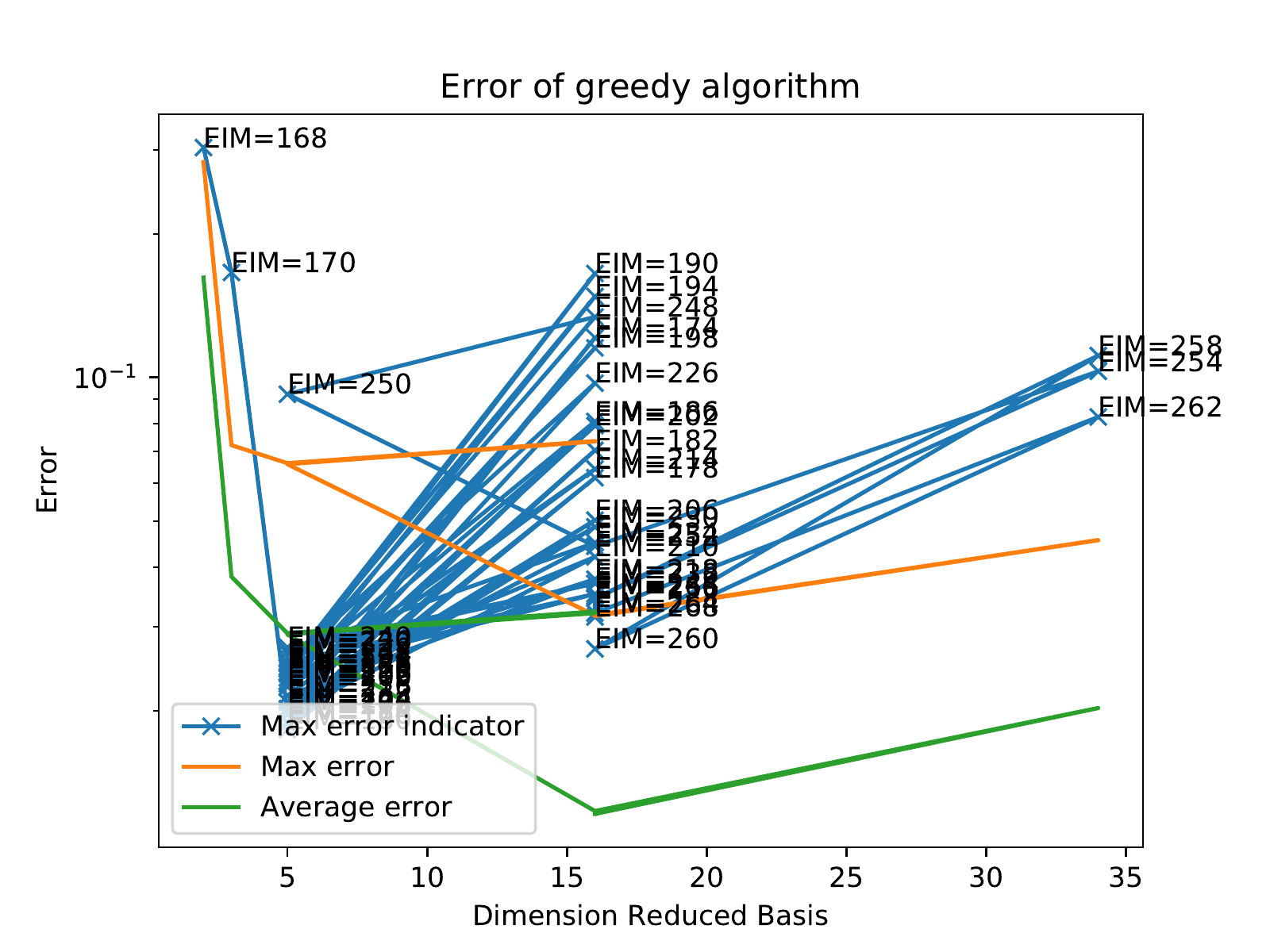}}
\end{subfigure}
	\begin{subfigure}[Offline error of Greedy in ALE framework\label{fig:test10:offlineALE}]{\includegraphics[width=0.45\textwidth, trim={0 0 20 40},clip]{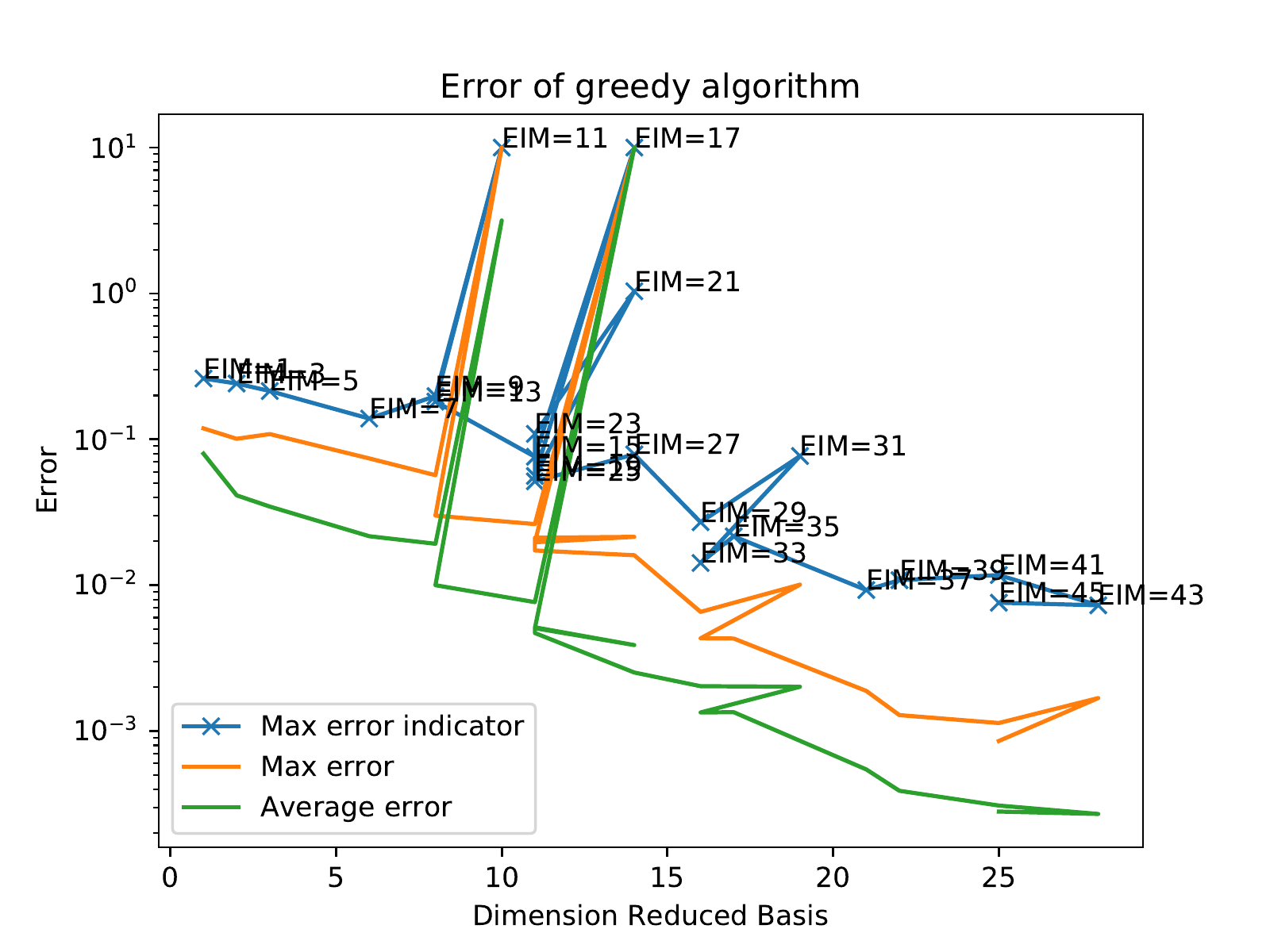}}
\end{subfigure}\\
	\begin{subfigure}[Eulerian FOM and RB solutions \label{fig:test10:onlineEulRB}]{\includegraphics[width=0.45\textwidth, trim={0 0 20 40},clip]{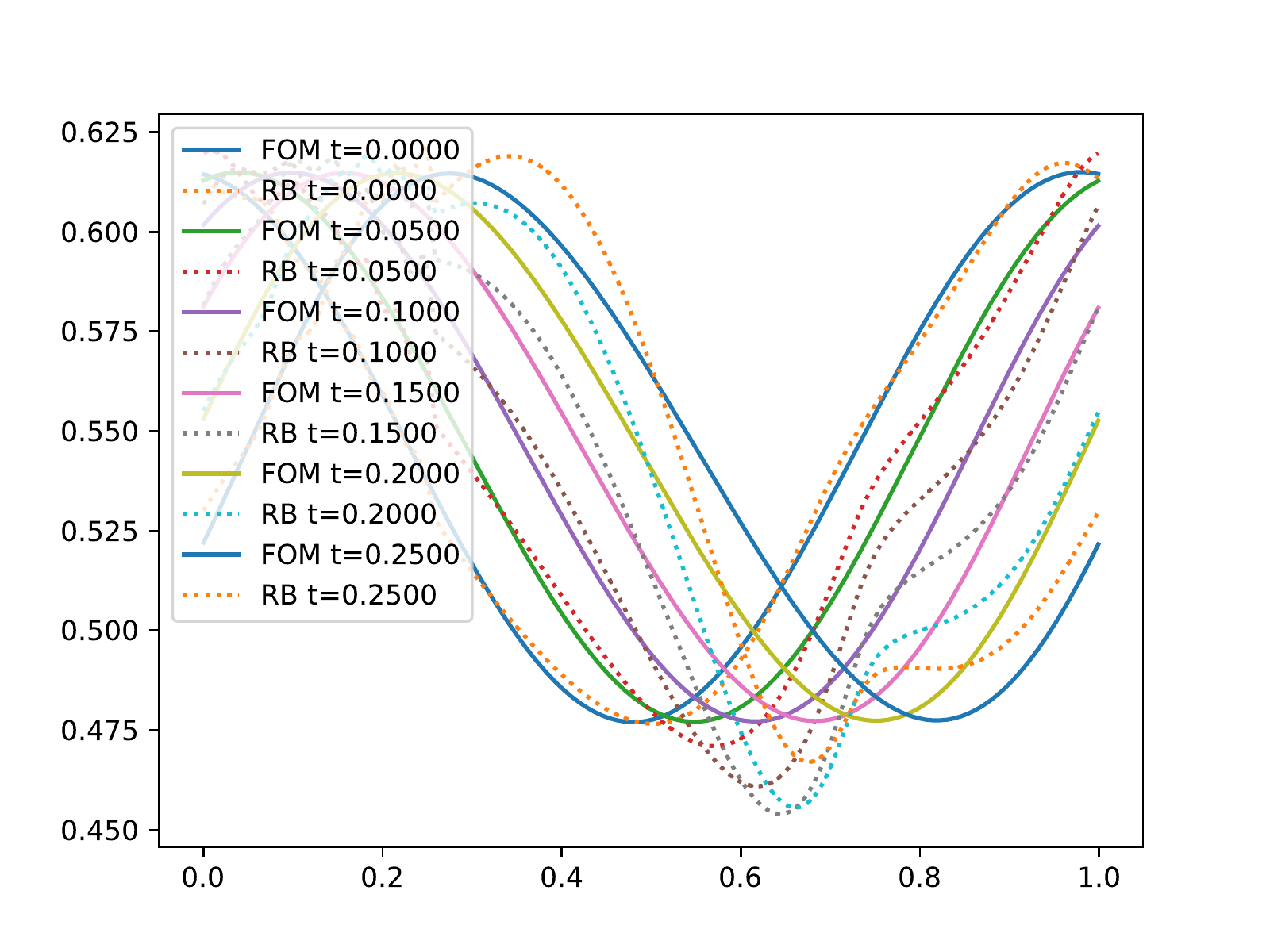}}
	\end{subfigure}
	\begin{subfigure}[ALE FOM and RB solutions\label{fig:test10:onlineALE}]{\includegraphics[width=0.45\textwidth, trim={0 0 20 40},clip]{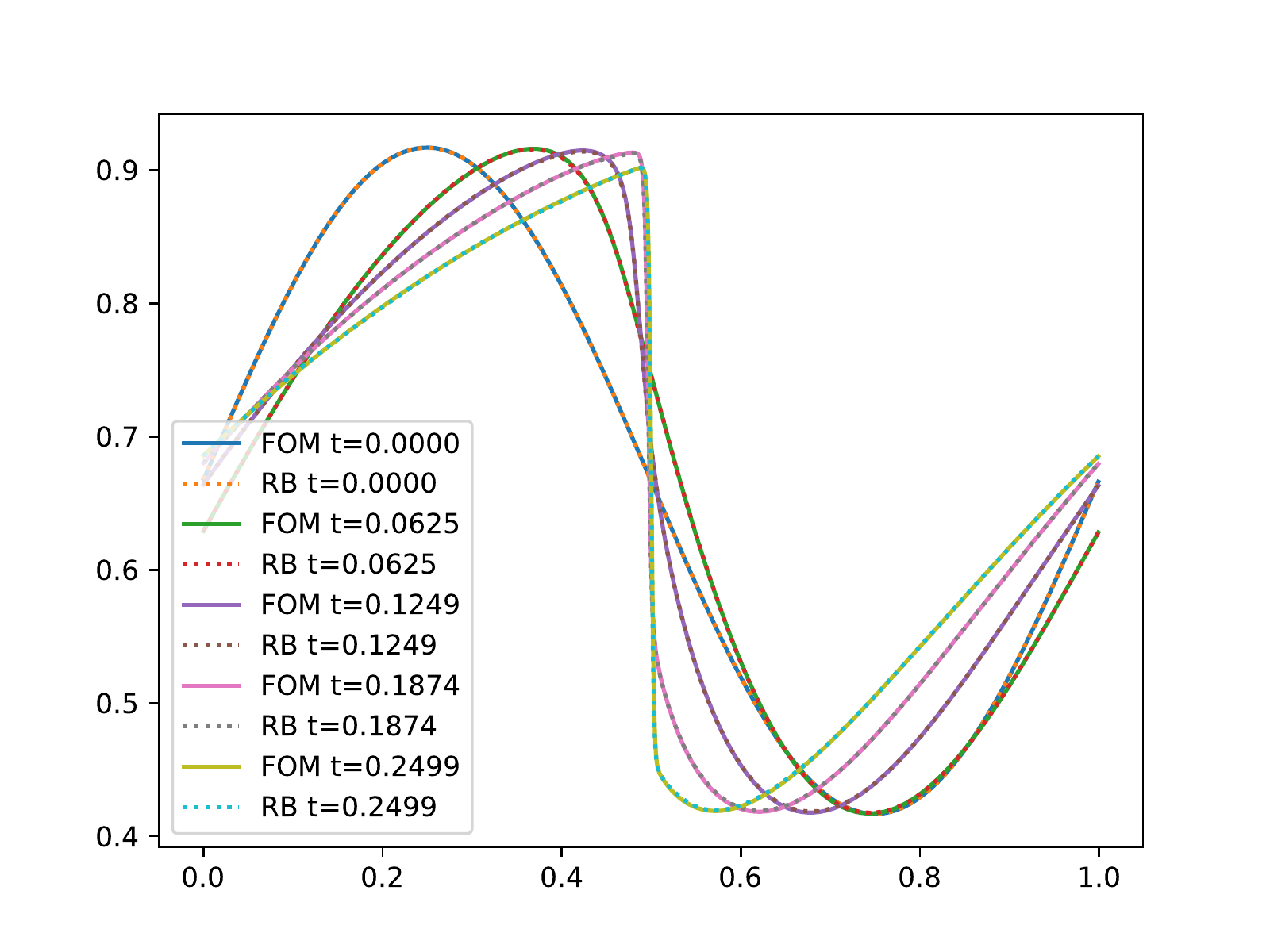}}
\end{subfigure}
\caption{Offline and online phases of Buckley simulations\label{fig:test10:ALEonline}}
\end{center}
\end{figure}

This problem is particularly tough to be analyzed. First of all, the shock is moving very quickly along the domain, secondly the speed is not a linear function of time and this makes the regression of the calibration map particularly hard. 
In \cref{fig:test1:training} we see that both polynomials of degree 1 and 2 can not at all represent this map, while third order and fourth order polynomials need more parameters to get a good enough regression map. 
We choose to use polynomials of fourth order in this example.

In \cref{fig:test1:offlineEu} we see that the error decay reaches with relatively few bases (19,41) the threshold, even if the algorithm needs to refine the EIM space a bit longer before having meaningful results.
If compared with the failure of the algorithm in the Eulerian framework, it is impressive how few bases we need.
The simulation in \cref{fig:test1:onlineFOM} shows the qualitatively behavior of a solution in the Eulerian framework and in \cref{fig:test1:onlineALE} we see how the solutions are aligned in the ALE framework and the quality of the reduced solution.

\subsection{Buckley Equation}
In this test, we solve the Buckley--Leverett equation \eqref{eq:buckley} on the domain $\Omega=[0,1]$, with a sine wave as initial conditions. This problem can develop in finite time a shock, accordingly to parameters. The initial conditions are, more precisely,
\begin{equation}\label{eq:BuckIC}
u_0(x,\bmu)= 0.5 + 0.2 \mu_1 + 0.3 \mu_1 \sin(2\pi (x-\mu_1-0.5)),\qquad \mu_1 \in [0.1,1],
\end{equation}
with $a=\mu_0 \in [0.001,2]$ till final time $T=0.25$. We use periodic boundary conditions and the translation \eqref{eq:translation} as calibration map. We select the steepest descending point as calibration point.

In \cref{fig:test10:training} we observe that also in this test the calibration map is not linear with respect to time. Indeed, all the polynomials struggle in obtaining a good approximation. Even the error of  the ANN is decaying too slowly to be used for the simulations. 
So, we pick the piecewise linear transformation to perform the ALE simulations since it can reach bigger area of the domain if we consider a bigger training set $N_{train}=100$.
In \cref{fig:test10:offlineEu} the \textit{offline} phase of the classical Eulerian is shown and we see that the algorithm fails in reaching the tolerance error even with $270$ EIM basis functions. 
We used, anyway, the resulting reduced spaces to perform an online phase in \cref{fig:test10:onlineEulRB}, where we see the difficulties in catching the right behavior of the solution. The data of these \textit{offline} and \textit{online} performances are stored in \cref{tab:RBdata}.
In \cref{fig:test10:offlineALE} the error of the ALE--PODEI--Greedy algorithm decays mush faster, overcoming some problems and resulting in final RB and EIM spaces of dimensions (25,45). The related simulation in the online phase in \cref{fig:test10:onlineALE} shows a better resolved solution which almost coincides with the exact one.

\section{Conclusion and Outlook}\label{sec:conclusion}
We have presented a new MOR technique that is able to effectively reduce the dimensions of the solution manifold of many advection dominated problems and is able to solve them in an \textit{online }phase gaining computational time.
It is applicable when there is one feature traveling in the domain that depends on time and parameters. 
The algorithm is capable of aligning all these features through a calibration map which is learned with different techniques, e.g. polynomial regression or artificial neural networks.
The \textit{offline} and the \textit{online} phases of the MOR algorithm are run on an ALE framework which modifies the original equation accordingly to the calibration process.
The preformed tests show that this algorithm is more robust than classical ones and obtains bigger reductions.
It is also general enough to deal with many nonlinear fluxes and different types of boundary conditions.

At the moment the ALE PODEI--Greedy has been tested on 1D scalar problem, but the generalization to systems with one speed is straightforward and will be object of future studies, even if these type of problems are not common at all. 
Also the generalization to 2D problems with one leading speed is possible and it can be done with different calibration maps. For example, the Gordon--Hall maps proposed in \cite{cagniartThesis} can serve the purpose. There, the calibration map would be a parametrized curve that follows the feature of interests, e.g. the shock or the wave.
It is possible to use this algorithm even for different advection dominated equations with few changes, according to the used scheme.

This approach still does not give an answer to more complicated problems, for example Riemann problems for systems of hyperbolic equations, where more shocks propagate from the same point.
The algorithm here proposed would fail in the inversion of the calibration map when the shocks meet and a singularity arises in it. 
The techniques used in other works for this situation do not meet the requirements of the proposed \textit{online} phase.

In future, we plan to extend this work to this type of problems, introducing a classification of different regimes or an usage of different nonlinear tools.

\section*{Acknowledgments}
Davide Torlo is supported by ITN ModCompShock project funded by the European Union’s Horizon 2020 research and innovation program under the Marie Sklodowska-Curie grant agreement No 642768 and by the SNF grant No 200020\_175784.\\
For this work, we also thank Tommaso Taddei for the scientific discussion over the Lagrangian and Eulerian methods, Maria Han Veiga and Fatemeh Mojarrad for the ideas and the trials on deep learning techniques.

\bibliographystyle{siam}
\bibliography{biblio_calib,biblio}

\end{document}